\pdfoutput=1
\RequirePackage{ifpdf}
\ifpdf 
\documentclass[pdftex]{sigma}
\else
\documentclass{sigma}
\fi

\usepackage{stmaryrd}
\newcommand\overt{\mathbin{\text{\rotatebox[origin=c]{-90}{$\ominus$}}}}

\usepackage{multirow}
\usepackage{phaistos}
\usepackage{accents}
\usepackage{tikz}
\usepackage{tikz-cd}

\newcommand\INT[1]{\operatorname{int}(#1)}

\newcommand{\interval}[1]{\llbracket #1 \rrbracket}
\renewcommand{\restriction}[2]{{#1\mathbin{\upharpoonright}_{#2}}}
\newcommand{\translucidation}[2]{{{#1}_{#2\leadsto0}}}
\newcommand{\linearspan}[1]{{\CC #1}}
\newcommand{\vertset}{{\overt}}
\def\uvert{{\mathcal{E}}^{\overt}}
\def\nc{\mathrm{NC}}
\def\ip{\mathrm{I}}
\def\monp{\mathrm{M}}

\def\ibp{\mathrm{I}}
\def\bnc{\mathrm{NC}}

\def\mbnc{\mathrm{M}}

\newcommand{\stdl}{\lessdot}
\newcommand{\stdg}{\gtrdot}

\newcommand{\stdgeq}{\,{\mathop{\geq\hspace{-10.5pt}\raisebox{1pt}[0pt][0pt]{$\cdot$}\hspace{3pt}}\nolimits}\,}
\newcommand{\stdgeqs}{{\mathop{\geq\hspace{-6pt}\raisebox{0.3pt}[0pt][0pt]{$\cdot$}\hspace{3pt}}\nolimits}}
\newcommand{\stdleq}{{\mathop{\leq\hspace{-7.5pt}\raisebox{1pt}[0pt][0pt]{$\cdot$}\hspace{1pt}}\nolimits}}
\newcommand{\bigwedgedot}{{\mathop{\bigwedge\hspace{-7.8pt}\raisebox{-1pt}[0pt][0pt]{$\cdot$}\hspace{3pt}}\nolimits}}
\newcommand{\bigwedgedots}{{\mathop{\bigwedge\hspace{-6.5pt}\raisebox{-1.5pt}[0pt][0pt]{$\cdot$}\hspace{2pt}}\nolimits}}
\newcommand{\bigveedot}{{\mathop{\bigvee\hspace{-8pt}\raisebox{1pt}[0pt][0pt]{$\cdot$}\hspace{3pt}}\nolimits}}
\newcommand{\wedgedotb}{{\mathop{\wedge\hspace{-5.2pt}\raisebox{-1pt}[0pt][0pt]{$\cdot$}\hspace{2pt}}\nolimits}}
\newcommand{\wedgedot}{{\mathop{\wedge\hspace{-4.2pt}\raisebox{-1pt}[0pt][0pt]{$\cdot$}\hspace{2pt}}\nolimits}}
\newcommand{\wedgedots}{{\mathop{\wedge\hspace{-4pt}\raisebox{-1.5pt}[0pt][0pt]{$\cdot$}\hspace{1pt}}\nolimits}}

\newcommand{\SWc}{{\mathfrak{T}\mathcal{W}}}
\newcommand\IWc{\mathcal{W}}
\newcommand{\SWint}{\SWc_{\sf I}}
\newcommand{\IWint}{{\IWc}_{\sf I}}
\newcommand{\BNCc}{\mathcal{N}\mathcal{C}}
\newcommand{\BPc}{{\mathcal{P}}}
\def\MBNCc{{\mathcal{M}}}
\def\BBc{{\mathcal{I}}}
\def\shBNCc{\BNCc_{\,{\color{red}\scalebox{0.5}{$\pmb{|}$}}}}
\newcommand{\shMBNCc}{\MBNCc_{\,{\color{red}\scalebox{0.5}{$\pmb{|}$}}}}
\DeclareMathOperator{\source}{\mathsf{source}}
\DeclareMathOperator{\target}{\mathsf{target}}
\DeclareMathOperator{\Type}{\mathsf{type}}
\newcommand{\vp}{\mathsf{m}^{\overt}}
\newcommand{\vcp}{\Delta^{\overt}}
\newcommand{\bvcp}{\bar{\Delta}^{\overt}}
\newcommand\vpS{\mathsf{m}^{\overt}_{\IWc}}
\newcommand\vpBP{\vp_{\BPc}}
\newcommand\vpBNC{\vp_{\BNCc}}
\newcommand\vcpS{\vcp_{\IWc}}
\newcommand{\bncinc}{\BNCc}
\newcommand{\bpinc}{\BPc}

\newcommand\barvcpS{\bar{\Delta}^{\overt}_{\IWc}}

\newcommand\hp{\mathsf{m}^{\ominus}}

\newcommand\hpIWc{\hp_{\IWc}}

\newcommand\mS{\mathsf{m}^{\ominus}_{\IWc}}
\newcommand\mSS{\mathsf{m}^{\ominus}_{\IWc\overt\IWc}}

\newcommand\vcpSpr{\vcp_{\prec}}
\newcommand\vcpSsu{\vcp_{\succ}}

\newcommand\LL{{\mathsf{L}}}
\newcommand\RR{{\mathsf{R}}}

\newcommand\cyan[1]{{\color{cyan}#1}}
\newcommand\sincp{shaded noncrossing bipartition}
\newcommand\sincps{shaded noncrossing bipartitions}

\newcommand\CAT[1]{\mathbf{#1}}
\newcommand{\Vect}{\CAT{Vect}}
\newcommand{\Set}{\CAT{Set}}
\newcommand\sgominus{\CAT{SemiGrp}(\ominus)}

\newcommand\icsw{ \Vect^{\SWc} }

\def\CC{\mathbb{C}}
\def\CIWc{{\linearspan{\IWc}}}

\newcommand\shgrI{{\mathsf M}_{\mathsf I}}
\newcommand\shalI{{\mathfrak m}_{\mathsf I}}

\newcommand{\shalII}{{\mathfrak m}_{\mathsf P}}

\newcommand{\shgrII}{{\mathsf M}_{\mathsf P}}

\newcommand{\Singleton}{\raisebox{-2pt}{\scalebox{0.7}{\PHcat}}}

\numberwithin{equation}{section}

\newtheorem{Theorem}{Theorem}[section]
\newtheorem{Proposition}[Theorem]{Proposition}
\newtheorem{Corollary}[Theorem]{Corollary}
\newtheorem{Observation}[Theorem]{Observation}

{ \theoremstyle{definition}
\newtheorem{Definition}[Theorem]{Definition}
\newtheorem{Notation}[Theorem]{Notation}
\newtheorem{Construction}[Theorem]{Construction}
\newtheorem{Example}[Theorem]{Example}
\newtheorem{Remark}[Theorem]{Remark}
}

\begin{document}

\allowdisplaybreaks

\newcommand{\arXivNumber}{2201.11747}

\renewcommand{\thefootnote}{}

\renewcommand{\PaperNumber}{006}

\FirstPageHeading

\ShortArticleName{Shuffle Algebras and Non-Commutative Probability for Pairs of Faces}

\ArticleName{Shuffle Algebras and Non-Commutative Probability\\ for Pairs of Faces\footnote{This paper is a~contribution to the Special Issue on Non-Commutative Algebra, Probability and Analysis in Action. The~full collection is available at \href{https://www.emis.de/journals/SIGMA/non-commutative-probability.html}{https://www.emis.de/journals/SIGMA/non-commutative-probability.html}}}

\Author{Joscha DIEHL~$^{\rm a}$, Malte GERHOLD~$^{\rm ab}$ and Nicolas GILLIERS~$^{\rm ac}$}

\AuthorNameForHeading{J.~Diehl, M.~Gerhold and N.~Gilliers}

\Address{$^{\rm a)}$~Universit\"at Greifswald, Institut f\"ur Mathematik und Informatik, \\
\hphantom{$^{\rm a)}$}~Walther-Rathenau-Str. 47, 17489 Greifswald, Germany}
\EmailD{\href{mailto:joscha.diehl@gmail.com}{joscha.diehl@gmail.com}}
\URLaddressD{\url{https://diehlj.github.io/}}

\Address{$^{\rm b)}$~Norwegian University of Science and Technology (NTNU), \\
\hphantom{$^{\rm b)}$}~Department of Mathematical Sciences, 7491 Trondheim, Norway}
\EmailD{\href{mailto:malte.gerhold@ntnu.no}{malte.gerhold@ntnu.no}}
\URLaddressD{\url{https://sites.google.com/view/malte-gerhold/}}

\Address{$^{\rm c)}$~Institut de Math\'ematiques de Toulouse, UMR5219, Universit\'e de Toulouse,\\
\hphantom{$^{\rm c)}$}~CNRS, UPS, F-31062 Toulouse, France}
\EmailD{\href{mailto:nicolas.gilliers@gmail.cm}{nicolas.gilliers@gmail.com}}
\URLaddressD{\url{http://nicolas-gilliers.github.com}}

\ArticleDates{Received March 30, 2022, in final form January 10, 2023; Published online January 31, 2023}

\Abstract{One can build an operatorial model for freeness by considering either the \emph{right-handed} or the \emph{left-handed} representation of algebras of operators acting on the free product of the underlying pointed Hilbert spaces. Considering both at the same time, that is, computing distributions of operators in the algebra generated by the left- and right-handed representations, led Voiculescu in 2013 to define and study \emph{bifreeness} and, in the sequel, triggered the development of an extension of noncommutative probability now frequently referred to as \emph{multi-faced} (two-faced in the example given above). Many examples of two-faced independences emerged these past years. Of great interest to us are \emph{biBoolean}, \emph{bifree} and \emph{type I bimonotone independences}. In this paper, we extend the preLie calculus pertaining to free, Boolean, and monotone moment-cumulant relations initiated by K.~Ebrahimi-Fard and F.~Patras to their above-mentioned two-faced equivalents.}


\Classification{46L53; 60A05; 18M05; 46L54; 16T05; 16T10; 16T30}

\renewcommand{\thefootnote}{\arabic{footnote}}
\setcounter{footnote}{0}

\section{Introduction}

\subsection{Background and motivation}
In \cite{voiculescu2014free}, Voiculescu introduced an extension of free probability; a new notion of independence, motivated by computations of joint distributions of \emph{left and right} creation and annihilation operators acting on the reduced free product of pointed Hilbert spaces. These operators prototype \emph{bifree independence} (introduced in \cite{charlesworth2015combinatorics, voiculescu2014free}) between pairs of random variables, in the same way that left (\emph{or} right) operators on their own prototype free independence. This operator algebraic root of bifreeness is supplemented by a combinatorial one~\cite{charlesworth2015combinatorics}, to which the poset of noncrossing \emph{bipartitions} and M\"obius inversion are central, extending the now well-developed combinatorial approach to freeness established by Speicher \cite{nica2006lectures, Speicher1994freecumulants}.
Since its inception, bifreeness developed into a \emph{theory of noncommutative probability for pairs of random variables}, also called \emph{two-faced} (left is one face, right is the other) random variables, with a steadily increasing set of two-faced independences.

Recall that Ben Ghorbal and Sch\"urmann \cite{BenGhorbalSchurmann} established a set of axioms defining the concept of independence in noncommutative probability. Muraki proved that there are only five such independences, namely free, monotone, anti-monotone, Boolean, and tensor independence~\cite{Muraki2003}. The axiomatic framework has been adapted to the two-faced (and more generally multi-faced-multi-state) case by Manzel and Sch\"urmann~\cite{ManzelSchurmann17}.

With the works of Skoufranis, Gu, Hasebe, Gerhold, Liu,
\cite{gerhold2017preprint,gu2020bi,gu2017bi,gu2019bi,liu2019}, several more
instances of two-faced independences besides bifreeness emerged, mixing one type of independence (free, monotone,
anti-monotone, boolean, or tensor) for left-sided random variables with another
one for right-sided random variables (Gerhold, Hasebe and Ulrich in \cite{GHU21p} and Gerhold and Var\v{s}o in \cite{GerholdVarso2023combinatorial}, see also \cite{Varso-PhD}, even establish continuous families of two-faced independences, including nontrivial bi-tensor independences).
For most of these independences, cumulants have been defined by exhibiting a
certain poset of bipartitions and applying the usual machinery of Rota's combinatorics and M\"obius inversion. Note that
neither does fixing the independence for the left and the right variables alone
determine the two-faced independence, nor do we know whether mixings for all
combinations of single-faced independences exist.

There are intimate relations between free, Boolean, and monotone cumulants.
The work of Ebrahimi-Fard and Patras~\cite{ebrahimi2015cumulants} shed new light on the interrelation of those three independences by showing that their corresponding additive convolution products originate from a certain splitting into two compatible co-operations of a coassociative coproduct on words on random variables. These (non-coassociative) co-operations satisfy the co-shuffle\footnote{Shuffle algebras have been popularized by Loday \cite{loday1763dialgebras} under the name \emph{dendriform algebras}.} relations.
Note that this coproduct is familiar to deformation theorists, it is called the \emph{double bar} of an algebra and was first defined by Baues~\cite{baues1981double}.

 Dualization of this coproduct and its splitting yields three different possibly nonassociative products satisfying the relations constitutive of \emph{a shuffle algebra}.
Subordination products, Bercovici--Pata bijection, and formulae relating the various families of cumulants can be cast in the shuffle (and the accompanying preLie) algebraic realm, as shown by Ebrahimi-Fard and Patras~\cite{EbrahimiFardPatras2019applications}.
Ebrahimi-Fard, Foissy, Kock and Patras recognized later that the discussed unshuffle structure is the restriction of an unshuffle structure on (polynomials on) noncrossing partitions, drawing a tight connection with the \emph{theory of operads}, see~\cite{ebrahimi2020operads}.
The present work takes inspiration from \cite{ebrahimi2020operads}, where an extension of the perspective on the moments-cumulants relations described in~\cite{ebrahimi2015cumulants} to operator-valued probability is proposed: both, \cite{ebrahimi2020operads} and this work, ultimately rely on an appropriate formalism for composing certain partitions, namely noncrossing partitions in \cite{ebrahimi2020operads} and noncrossing bipartitions here.
In particular, we use notions originating from higher category theory, \emph{duoidal categories} namely.

Of particular interest in this article are \emph{biBoolean independence}, \emph{type I bimonotone independence}, and \emph{bifree independence}. In those three cases, combinatorics of the moment-cumulant relations is guided by partitions that are noncrossing \emph{after} action of a specific permutation related to the two-faced structure.

In this work, we propose to study the
compositional structure of noncrossing bipartitions,
in order to exhibit a \emph{preLie algebra} connecting bifree, biBoolean and type~I bimonotone independences. This connection is similar to the single-faced case.
We build a~monoid and a~comonoid (but with respect to two different tensor products) supported by words on random variables and \emph{placeholders} of two types (left, right).
Recall that in the single-faced setting the unshuffle Hopf algebra is supported by words on words on random variables, which can be interpreted as incomplete words on random variables with placeholders, in our terminology, of only \emph{one} type.
The \emph{$($full$)$ unshuffle} coproduct of this monoid splits into two \emph{half-unshuffle} parts yielding the structure of a shuffle algebra on the dual space. We reveal a similar picture as for the single-faced setting: biBoolean moment-cumulant relations correspond to right half-shuffle exponentiation, bifree moment-cumulant relations correspond to left half-shuffle exponentiation, and type I bimonotone to full shuffle exponentiation.

One shortcoming of the established (un)shuffle algebra approach in the single-faced case is the failure to encompass classical moment-cumulant relations. To those corresponds another (commutative unshuffle) Hopf algebra whose relation to the previous one is unclear. We hope that the techniques we develop in this article will also help to shed some light in the future on the relationship between the two aforementioned unshuffle algebras.

\subsection{Outline}

\begin{itemize}\itemsep=0pt
\item In Section~\ref{sec:bckgnd} we give a brief overview of noncommutative probability for pairs of faces. We recall the definitions of the three posets of partitions at stake: interval, noncrossing and monotone bipartitions.
\item In Section~\ref{sec:categoryofwords} we introduce a certain category $\SWc$ which we use to index bipartitions and words on two-faced random variables in order to formalize their respective compositions. To that end, the \emph{incidence category} $\Vect^\SWc$ of $\SWc$-indexed collections of vector spaces is equipped with a monoidal structure, the vertical composition $\overt$.
 Additionally, we show that $\Vect^\SWc$ can be equipped with a \emph{horizontal} semigroupal product $\ominus$ satisfying a certain exchange relation with the aforementioned monoidal structure, that we call \emph{vertical}, to make a clear distinction between the two.
\item In Section~\ref{sec:words} we arrive at the $\SWc$-indexed collection of incomplete words supporting the algebraic structures we are interested in and we give our main result in Theorem~\ref{thm:unshuffle}.
\item In Section~\ref{sec:momentscumulantsrelations} we compute the three exponentials (the full shuffle exponential and the two half-shuffle exponentials), thereby recovering type I bimonotone, biBoolean, and bifree moment-cumulant relations.
\end{itemize}

\subsection{Conventions}
\begin{enumerate}\itemsep=0pt
\item We write $A\subset B$ if $A$ is a, not necessarily strict, subset of $B$.
\item For a positive integer $n$, $\llbracket n \rrbracket := \{1, \dots, n\}$.
\end{enumerate}

\section{Combinatorics of two-faced independences}\label{sec:bckgnd}

Two-faced independence is, roughly speaking, an independence relation for pairs
of noncommutative random variables. The axiomatization of the concept of
independence in noncommutative probability based on universal products
(suitable replacements for the product measure in classical probability) has
been initiated by Ben Ghorbal and Sch\"urmann in \cite{BenGhorbalSchurmann}. The
different independences focused on in the article at hand -- bifree,
bimonotone, and biBoolean independence -- are covered by a multivariate
extension of the theory of universal products, which was spelled out explicitly
by Manzel and Sch\"urmann in \cite{ManzelSchurmann17}, see also \cite[Section
3]{Gerhold21p}. Those three independences can be derived from Gu and
Skoufranis' c-bifree independence \cite{gu2017bi} much like free, monotone
and boolean independence are derived from Speicher and Bo\.zejko's c-free
independence \cite{BozejkoSpeicher1991}. We will not give details on the
definitions as we are only interested in the associated moment-cumulant
relations. Those have been established by Gu and Skoufranis  \cite{gu2019bi} for biBoolean
independence, Gu, Skoufranis and Hasebe  \cite{gu2020bi} for the bimonotone independence and Mastnak and Nica
\cite{mastnak2015double} for bifree independence, by relating the independences
to specific sets of bipartitions.

\subsection{Sets of bipartitions}
\label{sec:setsofbipartition}

A {\it partition} of a set $X$ is a set $\pi$ of subsets of $X$, called {\it blocks} of $\pi$, such that
\begin{itemize}\itemsep=0pt
	\item $V\neq\varnothing$ for all blocks $V\in \pi$,
	\item $V_1\cap V_2=\varnothing$ for all distinct blocks $V_1,V_2\in \pi$,
	\item $\displaystyle \bigcup_{V\in \pi} V=X$.
\end{itemize}

A partition $\pi$ is uniquely determined by its associated equivalence relation $\sim_\pi$ on $X$ defined by
$i \sim_\pi j$, $i,j\in X$ if and only if $i$ and $j$ are in the same block of $\pi$.

An {\it ordered partition} is a partition together with a total order on the blocks, identified as an injective labelling of the blocks with integers in $\llbracket |\pi|\rrbracket$.

From now on, we will only consider (ordered) partitions of totally ordered sets $X$, the order on the base space $X$ is fundamental to defining the following classes of partitions.

\begin{figure}[!ht]\centering
	\includegraphics[scale=1.4]{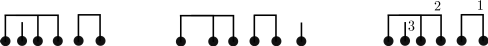}
	\caption{Example of (from left to right) a noncrossing partition, an interval partition and a monotone partition.}\label{fig:examplespartitions}
\end{figure}

\begin{Definition}
	Let $(X,<)$ be a totally ordered set.
	\begin{itemize}\itemsep=0pt
		\item A partition $\pi$ of $X$ is {\it noncrossing} if for all $a<b<c<d\in X$ such that $a\sim_\pi c$ and $b\sim_\pi d$ it follows that $a\sim_\pi b$, i.e., if the blocks of $\pi$ do not cross when drawn as an arc diagram.
		\item A partition $\pi$ of $X$ is an {\it interval partition} if for all $a<b<c\in X$ such that $a\sim_\pi c$ it follows that $a\sim_\pi b$, i.e., if all blocks of $\pi$ are intervals of $X$.
		\item An ordered partition $\pi$ of $X$ is a {\it monotone partition} if it is noncrossing and for all $a<b<c\in X$ such that $a\sim_\pi c$ it follows that the blocks $V\ni a,c$ and $W\ni b$ fulfill $V\leq W$, i.e., inner blocks are higher than outer blocks.
	\end{itemize}
	See examples in Figure~\ref{fig:examplespartitions}.
\end{Definition}
\begin{Notation}The sets of all noncrossing, interval, and monotone partitions of $(X,<)$ are denoted by $\nc(X,<)$, $\ip(X,<)$, and $\monp(X,<)$, respectively (or simply $\nc(X)$, $\ip(X)$, and $\monp(X)$ if there is no risk of confusion).
\end{Notation}

\begin{figure}[!ht]\centering
	\includegraphics[scale=0.45]{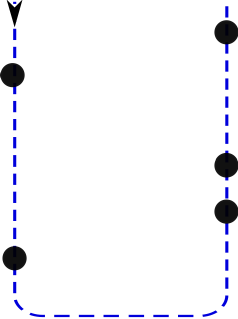}
	\caption{The ordered biset $X = \{1 < 2 < 3 < 4 < 5\}$ with left-right structure $\alpha(1)={\RR}$, $\alpha(2)={\LL}$, $\alpha(3)={\RR}$, $\alpha(4)={\RR}$, $\alpha(5)={\LL}$ and necklace order $2 \stdl 5 \stdl 4 \stdl 3 \stdl 1$.}\label{fig:biset}
\end{figure}

A {\it biset} is a set $X$ together with a map $\alpha\colon X\to\{\LL,\RR\}$ that we call \emph{left-right structure} on~$X$. The elements of $X_\LL:=\alpha^{-1}(\LL)$ are called {\it left}, the elements of $X_\RR:=\alpha^{-1}(\RR)$ are called {\it right}.
A~biset $(X,\alpha)$ with a total order on $X$ is called {\it ordered biset}. We depict ordered bisets by diagrams as in Figure~\ref{fig:biset}: starting from the top, we put in increasing order left elements on a left string and right elements on a right string. Given an ordered biset, we refer to its given total order $<$ on $X$ as the {\it natural order}. We define the {\it necklace order} $\stdl$ in the following way:
\begin{itemize}\itemsep=0pt
	\item for $x,y\in X_\LL$, $x\stdl y\iff x<y$, i.e., the order for left elements remains unchanged,
	\item for $x,y\in X_\RR$, $x\stdl y\iff x>y$, i.e., the order for right elements is reversed,
	\item for $x\in X_\LL$ and $y\in X_\RR$, $x\stdl y$, i.e., all left elements are smaller than all right elements.
\end{itemize}

\begin{Notation}\label{notation:infinity}
Given an ordered biset $(X,\stdl)$, we introduce two elements $-\infty_{\stdl}$ and $+\infty_{\stdl}$ and extend the order ${\stdl}$ to $X\cup \{-\infty_{\stdl}, +\infty_{\stdl}\}$ by
\[
-\infty_{\stdl} \stdl x \stdl +\infty_{\stdl},\qquad x \in X.
\]
\end{Notation}
The necklace order can easily be read from the corresponding diagram. Having another look at Figure~\ref{fig:biset}, we connected the vertical strings at the bottom with a horizontal link. The necklace order is the order in which the elements of $X$ appear when moving along the dotted line
of Figure~\ref{fig:biset}, starting from the top left.

\begin{Remark}\label{rk:categorybiset}
Ordered bisets form a category {\rm BiSet} with morphisms the injective set maps which respect the natural orders and the left-right structures. As a consequence, morphisms preserve also the necklace orders.
\end{Remark}
To emphasize the dependence of the necklace order $\stdl$ on the left-right structure $\alpha$, we use at times a subscript and write $\stdl_\alpha$.

\begin{Definition}	\label{def:bipartitions}
	A partition of a biset is also called {\it bipartition}. Let $X$ be an ordered biset. An (ordered) bipartition is called {\it noncrossing bipartition}, {\it interval bipartition}, or {\it monotone bipartition} if it is a noncrossing partition, an interval partition, or a monotone partition of $(X,\stdl)$, respectively.

	The sets of all noncrossing, interval, and monotone bipartitions of $X$ are denoted by~$\bnc(X)$, $\ibp(X)$, and $\mbnc(X)$ respectively.
\end{Definition}

\begin{Remark}\label{rk:finite-ordered-bisets}
	We will almost exclusively deal with finite ordered bisets. These can be easily described by words in $\{\LL,\RR\}^\star$.
	Indeed, let $\alpha=\alpha(1)\dots \alpha(n)\in \{\LL,\RR\}^\star$ be a word (note that we do not distinguish between a word $\alpha$ of length $n$ and the corresponding map $k\mapsto \alpha(k)\colon \interval{n}\to \{\LL,\RR\}$). Then $\interval{n}$ becomes an ordered biset with the natural order induced from $\mathbb N$ and the left-right structures given by $\alpha$.
 We write $X_\alpha$ for the ordered biset $(\interval{n},\alpha)$.
\end{Remark}
	Let us say that two ordered bisets are isomorphic if there is a bijection respecting the natural orders and the left-right structures. Then the prescription $\alpha\mapsto X_\alpha$ yields a bijection between words in $\{\LL,\RR\}^\star$ and isomorphism classes of finite ordered bisets.

	We will usually identify $\alpha$ with $X_\alpha$ in the following.

\begin{Notation}
Given $\alpha \in \{\sf L, \sf R\}^{\star}$, we use the notations $\bnc(\alpha)$, $\ibp(\alpha)$, and $\mbnc(\alpha)$ for the associated sets of noncrossing, interval, and monotone bipartitions. A bipartition $\pi$ of a finite ordered biset $X_\alpha$ will usually be denoted by the pair $(\pi,\alpha)$.
\end{Notation}
\begin{Remark}
For a finite ordered biset $w \in \{{\sf L,\sf R}\}^{\star}$, Definition~\ref{def:bipartitions} of a noncrossing bipartition is equivalent to the definition given in \cite{charlesworth2015combinatorics}; indeed, while they apply a permutation to the partition, we apply the inverse permutation to the partitioned set.
\end{Remark}
\begin{figure}[!htb]\centering
	\includegraphics[scale=0.65]{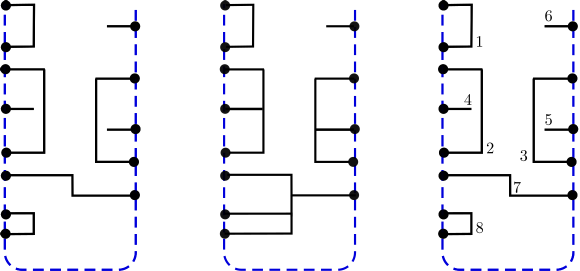}
	\caption{Example (from left to right) of a noncrossing bipartition, interval bipartition and a monotone bipartition.}\label{fig:bipartitions}
\end{figure}

\subsection{Cumulants}
\label{sec:cumulants}

Let us start with a brief historical account of classical and noncommutative cumulants and the shuffle perspective on their relation to moments. As mentioned in Hald's review article \cite{hald2000cumulants} on the early history of cumulants,
classical cumulants of a random variable were first introduced by Thorvald N.~Thiele in 1889 via an inductive formula computing the coefficients in a Gram--Charlier expansion of a density under the name of halfinvariants. Only ten years later, Thiele recognized them as the coefficients of the logarithm of the generating function of the moments of the random variable.

Almost a century later, Speed \cite{speed1983cumulants} proved a combinatorial formula for computing cumulants in terms of the moments: they are obtained by M\"obius inversion on the lattice of set-partitions of the partitioned moments of the random variable.
Around the same period of time, Voiculescu laid the foundations of the theory of free probability by introducing the notions of \emph{free independence} and \emph{free additive convolution} \cite{Voiculescu1985symmetries,Voiculescu1986addition}.

In 1993, Speicher introduced the counterpart of Thiele's cumulants in free probability theory, that is, a set of multilinear functionals linearizing free independence \cite{Speicher1994freecumulants}. The observation that moments of words on free random variables that are arranged according to a noncrossing partition
factorize (but not all moments) hints at the definition for free cumulants: replace the lattice of set-partitions by the lattice of noncrossing partitions in Speed's combinatorial formula.

Together with the emergence of new independences (Boolean, monotone, antimonotone) came the corresponding sets of cumulants, see \cite{SpeicherWoroudi1997bool} for the Boolean and \cite{HasebeSaigo2011monotone} for the {(anti\nobreakdash-)}monotone case. Lehner established a very general approach to noncommutative cumulants via so-called exchangeability systems \cite{Lehner2004exchange}. To include the monotone case, Lehner and Hasebe generalized the approach by weakening exchangeability to spreadability \cite{HasebeLehner2017p}.
A unified framework for cumulants with respect to universal product independences is given by Manzel and Sch\"urmann \cite{ManzelSchurmann17}. In this case, the moment-cumulant relations can always be expressed via a Hopf-algebraic exponential and logarithm, see also \cite{Gerhold21p}.
Formulas relating different sets of noncommutative cumulants of a random variable were proved in \cite{ArizmendiHasebeLehnerVargas2015cumulants} and are still subject of great attention. Using shuffle algebras, Ebrahimi-Fard and Patras established a single unified Hopf algebraic setting which relates free, monotone, and boolean cumulants \cite{ebrahimi2015cumulants} as well as the corresponding additive convolution products in a group-theoretic sense \cite{EbrahimiFardPatras2019applications}.

We quickly recall the moment-cumulant relations which determine the free cumulants $\kappa_n$, the monotone cumulants $K_n$ and the boolean cumulants $B_n$. (The exponential character of the formulas is most visible in the monotone case because monotone partitions are ordered partitions. In the other cases, the division by a factorial is hidden by the use of non-ordered partitions; there are $|\pi|!$ ordered partitions corresponding to the same partition $\pi$.) Let $(\mathcal A,\varphi)$ be an algebraic probability space, i.e., $\mathcal A$ is an unital algebra over $\mathbb C$ and $\varphi\colon \mathcal A\to\CC$ an unital linear functional. For a tuple $(a_1,\dots, a_n)$ and a subset $I=\{i_1<\cdots <i_k\}$, we write $\restriction{(a_1,\dots,a_n)}{I}:=(a_{i_1},\dots,a_{i_k})$. Then $\kappa_n,K_n,B_n\colon \mathcal{A}^n\to\CC$ are the unique sequences of multilinear maps such that
\begin{align}
	\label{eqn:mcrB}
	\tag{\sf MC$_{\text{B}}$} \varphi(a_1\cdots a_n) & = \sum_{\pi \in \ip(n)}\prod_{V \in \pi} B_{|V|}(\restriction{(a_1,\dots,a_n)}{V}) \\
	\label{eqn:mcmrf}
	\tag{\sf MC$_{\text{f}}$}
	 & = \sum_{\pi \in \nc(n)}\prod_{V \in \pi} \kappa_{|V|}(\restriction{(a_1,\dots,a_n)}{V}) \\
	\label{eqn:mcrmo} & \tag{\sf MC$_{\text{m}}$} = \sum_{\pi \in \monp(n)}\frac{1}{|\pi|!}\prod_{V \in \pi} K_{|V|}(\restriction{(a_1,\dots,a_n)}{V})
\end{align}
holds for all $a_1,\dots, a_n\in \mathcal A$ (with the convention that the empty product is 1; note that for $n=0$ the empty partition then yields one summand 1 in all three cases).

The two-faced extensions of those moment-cumulant relations have also been established (see \cite{CNS15:random-variables,gu2020bi, gu2019bi,mastnak2015double} for the specific cases or \cite{ManzelSchurmann17} for the general definition of cumulants with respect to a~universal product). Let $(\mathcal A,\varphi)$ be an algebraic probability space and fix subalgebras $\mathcal A^{\LL}$, $\mathcal A^{\RR}$ of \emph{left} and \emph{right} random variables (a \emph{pair of faces}).
For a word $\alpha=(\alpha_1,\dots, \alpha_n)\in \{\LL,\RR\}^\star$ we set $\mathcal A^\alpha:=\mathcal A^{\alpha_1}\times \cdots \times \mathcal A^{\alpha_n}$.
We also extend the restriction notation to words (which are also tuples, just with entries from $\{\LL,\RR\}$ instead of $\mathcal A$) in an obvious way, $\restriction{(\alpha_1,\dots,\alpha_n)}{\{i_1<\cdots<i_k\}}:=(\alpha_{i_1},\dots,\alpha_{i_k})$. Then the bifree, biBoolean and bimonotone cumulants $\kappa_\alpha,K_\alpha,B_\alpha\colon \mathcal A^\alpha\to\CC$ ($\alpha\in \{\LL,\RR\}^\star$) are the unique families of multilinear maps such that
\begin{align}
	\label{eqn:mcrbB}
	\tag{\sf MC$_{\text{bB}}$} \varphi(a_1\cdots a_n) & = \sum_{\pi \in \ibp(\alpha)}\prod_{V \in \pi} B_{\restriction{\alpha}{ V}}(\restriction{(a_1,\dots,a_n)}{V}) & & \text{\cite{gu2019bi}} \\
	\label{eqn:mcrbf}
	\tag{\sf MC$_{\text{bf}}$}
	 & = \sum_{\pi \in \bnc(\alpha)}\prod_{V \in \pi} \kappa_{\restriction{\alpha}{ V}}(\restriction{(a_1,\dots,a_n)}{V}) & & \text{\cite{CNS15:random-variables, mastnak2015double}} \\
	\label{eqn:mcrm}
	 & \tag{\sf MC$_{\text{Im}}$} = \sum_{\pi \in \mbnc(\alpha)}\frac{1}{|\pi|!}\prod_{V \in \pi} K_{\restriction{\alpha}{ V}}(\restriction{(a_1,\dots,a_n)}{V}) & & \text{\cite{gu2020bi}}
\end{align}
holds for all $(a_1,\dots, a_n)\in \mathcal A^\alpha$.

\begin{Remark}
 In the literature, the moment cumulant relations are often given for arbitrary $a_1,\dots, a_n$ in some algebraic probability space $(\mathcal B,\Phi)$. Given $a_1,\dots, a_n\in \mathcal B$ and $\alpha\in \{\LL,\RR\}^n$, the given formulas extend as follows. One can construct the algebraic probability space $(\mathcal A,\varphi):=(\mathcal A\sqcup \mathcal A, \Phi\circ\mu)$ ($\mu=\mathrm{id}\sqcup\mathrm{id}\colon \mathcal A\sqcup \mathcal A \to \mathcal A$ the canonical homomorphism), which clearly contains two copies of $\mathcal A$ as subalgebras, denoted $\mathcal A^\LL$ and $\mathcal A^\RR$, respectively. When interpreting $a_1\dots a_n$ as an element of $\mathcal A^\alpha$, the given formulas apply.
\end{Remark}

BiBoolean and bifree independence of a family of two-faced random variables are characterized by the vanishing of mixed cumulants. Concretely, two-faced random variables $\big(a^{\LL}_1, a^{\RR}_1\big)$ and $\big(a^{\LL}_2, a^{\RR}_2\big)$ are bifree (resp.\ biBoolean) independent if and only if
\[
\kappa_{\alpha}\big(a_{\varepsilon(1)}^{\alpha(1)},\dots,a^{\alpha(n)}_{\varepsilon(n)}\big) = 0, \qquad \big(\text{resp.} \ B_{\alpha}\big(a_{\varepsilon(1)}^{\alpha(1)},\dots,a^{\alpha(n)}_{\varepsilon(n)}\big)=0\big)
\]
whenever $\varepsilon\in\{1,2\}^n$ is not constant (with $\mathcal A^\LL$, $\mathcal A^\RR$ the algebras generated by the $a_k^\LL$ and $a_k^\RR$, respectively). For (bi-)monotone independence the situation is more involved because those are nonsymmetric independences, see~\cite{ManzelSchurmann17} and~\cite{HasebeLehner2017p} for details on this problem.

In \cite{ebrahimi2015cumulants}, Ebrahimi-Fard and Patras present a unified geometric approach to monotone, Boolean and free cumulants, expressing the three moment-cumulant relations as logarithmic-exponential correspondences on a shuffle algebra obtained as the dual of an unshuffle coalgebra.
In their work, cumulants correspond to a certain \emph{element of a $($pre$)$Lie algebra} $\mathfrak{g}$ and moments of a~random variable to \emph{a point in the corresponding group}~$G$.
The two half-shuffle products, $\prec$~and~$\succ$, yield \emph{two exponential maps} besides the one corresponding to the associative (full shuffle) convolution product.
These exponentials $\exp_{\prec}\colon \mathfrak{g}\to G$ and $\exp_{\succ}\colon \mathfrak{g}\to G$ are solutions to the following fixed points equations,
\begin{equation*}
	\exp_{\prec}(k) = \varepsilon + k \prec \exp_{\prec}(k),\qquad \exp_{\succ}(b) = \varepsilon + \exp_{\succ}(b) \succ b,
\end{equation*}
where $k,b \in \mathfrak g$. The above equations can be shown to be equivalent to the inductive formula computing moments from free cumulants for the first one and moments from Boolean cumulants for the second one, see \cite{EbrahimiFardPatras2019applications}. The three following sections aim at developing this algebraic interpretation for the two-faced moment-cumulant relations \eqref{eqn:mcrm}, \eqref{eqn:mcrbB} and \eqref{eqn:mcrbf}.

\section{A M\"obius category of (translucent) words}\label{sec:categoryofwords}

In this section, we introduce the categorical setting underlying the rest of the
paper. This categorical underpinning is reminiscent of the one underlying
operads, but it differs in the fact that the collections are indexed by words
in $\{\LL, \RR\}^{\star}$. This is necessary for formalizing the composition
of bipartitions and words on random variables in a two-faced probability space.
In the single-faced case, integers are sufficient to address the operad's inputs, but
here we additionally have to take the left-right structure into account.

\begin{figure} \centering
 \includegraphics[scale=0.65]{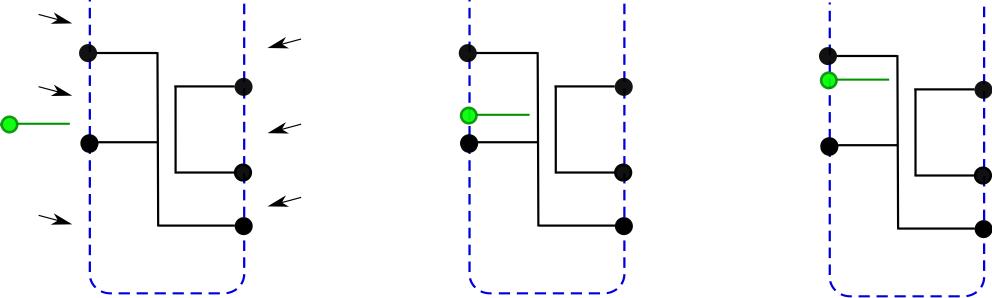}
 \caption{Ambiguity in composing noncrossing bipartition. We want to ``insert'' the singleton (drawn in green) in a space delimited by two consecutive elements (for the order defined by the left-right structure). For the chosen slot, there is an ambiguity about the linear ordering which can be lifted in two ways.} \label{fig:ambiguity}
\end{figure}

When trying to define a composition of noncrossing bipartitions, drawing their diagrams hints at a suitable set of \emph{inputs}:
given a noncrossing bipartition $(\pi, \alpha)$, each segment connecting two consecutive points for the necklace order $\stdl_{\alpha}$ can be thought
of as an input into which one can insert a certain type of noncrossing bipartition.
For example, in a segment bounded by two left points, one can certainly
insert a noncrossing partition of a set all of whose points are also on the
left, in order to obtain a new noncrossing bipartition. However, there is a
certain ambiguity in composing two noncrossing bipartitions this way that is
demonstrated in Figure~\ref{fig:ambiguity}.
In order to address the ambiguity in this example, one has to choose the relative positions of
the newly inserted points and the right-sided points that sit between the two
boundary points of the input segment (with respect to the natural order). Such
an ordering will be picked by adding \emph{placeholders} to the inserted
bipartitions, whose sole purpose is to settle the aforementioned ambiguities
in composing bipartitions.

In Section~\ref{subsec:category-of-translucent-words}, we introduce \emph{translucent words} which form the morphisms of a category $\SWc$. Collections of sets indexed by translucent words form a monoidal category introduced in Section~\ref{subsec:monoidalindexedmord-set}. The same holds for collections of vector spaces indexed by translucent words, which are discussed in Section~\ref{subsec:monoidalindexedmord-vect}.

\subsection{The category of translucent words}
\label{subsec:category-of-translucent-words}

\begin{Definition}[translucent words]	\label{def:translucentword}
 A {\it translucent word} (over the alphabet $\{\LL,\RR\}$) is a pair $t=(\alpha,b)$
	where
 $\alpha \in \{ \LL, \RR \}^\star$ is a word in the letters $\LL$, $\RR$ and $b \in \{0,1\}^\star$ is a ``Boolean word'' of equal length. We denote by $t(i)$ the pair $(\alpha(i),b(i))$. A translucent word can thus also be interpreted as a word in $(\{\LL,\RR\}\times\{0,1\})^\star$, in particular, $|t|=|\alpha|=|b|$. Given a translucent word $t=(\alpha,b)$, we call its word component $\alpha_t:=\alpha$ the {\it left-right structure} of $t$ and its Boolean component $b_t:=b$ the {\it translucent-opaque structure}.
\end{Definition}
We denote by $\SWc$ the set of all translucent words (including the {\it empty translucent word} $\varnothing := (\varnothing,\varnothing)$).

Given a translucent word $t$, we interpret $\alpha_t$ as a finite ordered biset as in Remark~\ref{rk:finite-ordered-bisets}, while $b_t$ marks the positions tagged $0$ ``translucent'' and the positions tagged $1$ ``opaque'' (the terminology will become clearer below).

In the sequel, we denote by ${\mathbf{0}}_n$ the Boolean word of length $n$ filled with $0$'s and by ${\mathbf{1}}_n$ the Boolean word with only $1$'s. Usually, the length is clear from the context, so we will also omit the index $n$. In particular, we write $(\alpha,\mathbf{1})$ for the translucent word with left-right structure $\alpha$ and all positions opaque.

 \begin{Definition}
 We call a set $X$ together with a \emph{translucent-opaque} structure (i.e., a map $X\to\{0,1\}$) a {\it translucent set} and consequently an ordered biset with a translucent-opaque structure is called {\it translucent ordered biset}.
 \end{Definition}
 In the corresponding diagrams, white dots stand for translucent points and black dots stand for opaque points. The finite ordered translucent set associated with a translucent word $t$ is sometimes denoted $X_t$, but usually, we do not distinguish between $t$ and $X_t$.

\begin{figure}[!ht]\centering
	\includegraphics[scale=0.60]{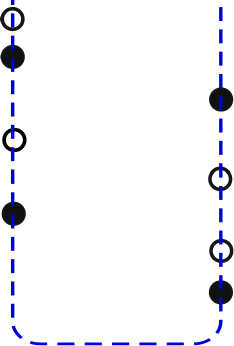}
	\caption{Graphical representation of the translucent word $({\sf LLRLRLRR}, 01100101)$.} \label{fig:incompletewords}
\end{figure}

\begin{Example}
 $({\sf LLRLRLRR}, 01100101)$
	is a translucent word. It is visualized in Figure~\ref{fig:incompletewords}.
\end{Example}

Further below we make use of a few items of notation and terminology that we introduce now. Let $t=(\alpha,b) \in \SWc$.
\begin{enumerate}\itemsep=0pt
	\item For a word $w\in A^\star$ over any alphabet and a letter $a\in A$, we write \[[w]_a:=\{j \in \interval{|w|}\colon w(j)=a\}.\]
	We extend the notation to translucent words: for $t=(\alpha,b) \in \SWc$, we set $[t]_0:=[b]_0$, $[t]_\LL:=[\alpha]_{\LL}$, etc.
	\item For a nonempty subset $I=\{i_1<\cdots<i_p\}$ of $\interval{|t|}$, the {\it restriction of $t$ to $I$}, denoted $\restriction{t}{I}=(\restriction{\alpha}{I},\restriction{b}{I})$, is the pair of subwords
	 \begin{equation*}
		 \restriction{\alpha}{I} := \alpha(i_1)\cdots \alpha(i_p),\qquad \restriction{b}{I}:=b(i_1)\cdots b(i_p).
	 \end{equation*}
	\item
	 For a subset $I$ as before, one can also define the
 {\it translucidation of $t$ at $I$}, denoted $\translucidation{t}{I}=(\translucidation{\alpha}{I},\translucidation{b}{I})$,
 by \emph{turning} all positions in $I$ translucid, i.e.,
 	 \begin{equation}
 \label{eq:tranlucidation}
		 \translucidation{\alpha}{I} := \alpha,\qquad
		 \translucidation{b}{I}(i)=
		 \begin{cases}
		 0, &i \in I, \\
		 b(i),&i\not\in I.
		 \end{cases}
	 \end{equation}
\end{enumerate}

\begin{Remark}
\label{rk:translucent-words-and-morphisms-of-ordered-bisets}
	A translucent word $t=(\alpha,b)$ defines a morphism \[ \iota_t\colon \ X_{\restriction{\alpha}{[t]_0}}=(\interval{|[t]_0|},\restriction{\alpha}{[t]_0})\to X_\alpha=(\interval{|\alpha|},\alpha)\] in the category {\rm BiSet},
	\[ \iota_t(k):=i_k,\qquad [t]_0=\{i_1<\cdots <i_{n}\}.\]
\end{Remark}

	\begin{Proposition}
	The correspondence between translucent words and morphisms between finite ordered bisets of the form $X_\alpha$ is one-to-one.
	\end{Proposition}
	\begin{proof}
 The morphism $\iota_t$ defined above is clearly strictly increasing and also preserves the left-right structure because $\restriction{\alpha}{[t]_0}(k)=\alpha(i_k)=\alpha(\iota_t(k))$ for all $k\in 1,\dots,n=|[t]_0|$.
	Conversely, if $\iota\colon X_\beta\to X_\alpha$ is a morphism of ordered bisets, then it is easy to see that $\iota=\iota_t$ for $t=(\alpha,b)$ where
	\[
	b(k)=\begin{cases}
	 0 & \text{if $k\in \iota(\interval{|\beta|})$},\\
	 1&\text{otherwise.}
	\end{cases}\tag*{\qed}
	\]\renewcommand{\qed}{}
	\end{proof}

 Instead of working directly with the category of finite ordered bisets, we prefer to concretely define a composition of translucent words, thereby turning them into the morphism class of a~(countable) category which is equivalent to the opposite of the category of finite bisets.

\begin{figure}[!ht]
	\centering
	\begin{tikzcd}[column sep = large]
		{\LL} & {\LL\LL\RR\RR} \arrow[l, "1011"'] & {\sf \LL\LL\RR\LL\RR\LL\RR\RR} \arrow[l,xshift=2mm, "{01100101}"'] \arrow[ll, "11101111"', bend left]
	\end{tikzcd}
	\caption{Composition of $s = ({\sf LLRR}, 1011)$ and $t = ({\sf LLRLRLRR}, 01100101)$ in $\SWc$. Each arrow stands for the translucent word whose left-right pattern is given by its source and whose translucent-opaque pattern is the one written above.}	\label{fig:exampleComposition}
\end{figure}
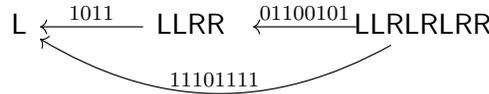

Given a translucent word $t=(\alpha, b)$, define its {\it source} to be the word $\source{t}:=\alpha\in \{\LL,\RR\}^\star$ and its {\it target} to be the word $\target{t}:=\restriction{\alpha}{[t]_0}\in \{\LL,\RR\}^\star$.
We say that two translucent words~$s$,~$t$ are {\it composable} if $\source s=\target t$, i.e., if the left-right pattern of $s$ matches the left-right pattern of the translucent points of $t$; in this case we define their {\it composition} $s\circ t:=r$,
where~$r$ is the unique translucent word with $\restriction{r}{[t]_1}=\restriction{t}{[t_1]}$ and $\restriction{r}{[t]_0}=s$, i.e., opaque points of $t$ stay opaque in the composition while translucent points of $t$ are coloured according to $s$.
Roughly speaking, the Boolean word $b_s$ overwrites the zeroes of $b_t$ to produce $b_r$, see Figure~\ref{fig:exampleComposition}.

\begin{figure}[!ht]\centering
	\includegraphics[scale=0.45]{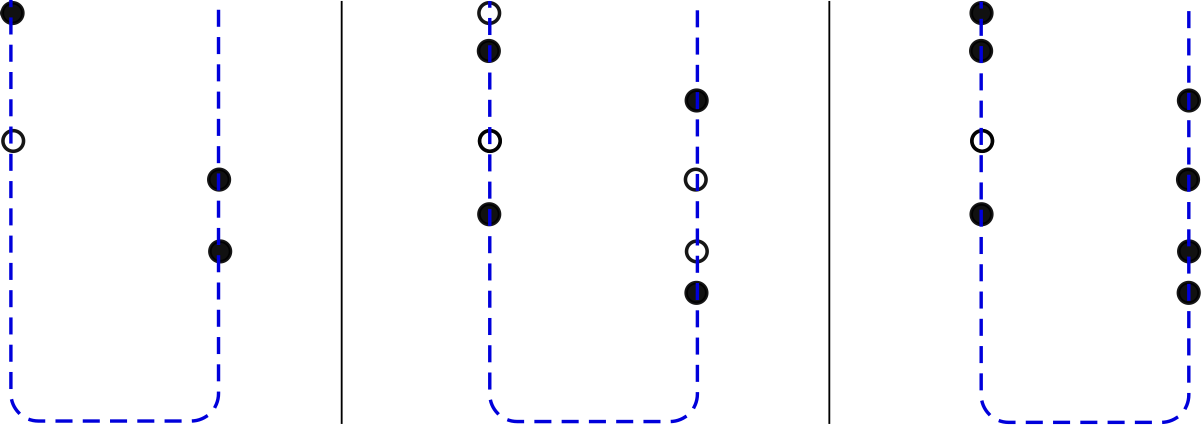}
	\caption{Graphical representation of composition of the two translucent words.
		On the left $s = ({\sf LLRR}, 1011)$, in the middle $t = ({\sf LLRLRLRR}, 01100101)$
		and on the right $s\circ t = ({\sf LLRLRLRR}, 11101111)$.}		\label{fig:exampleComposition-diagram}
\end{figure}

The composition of translucent words corresponds to the superimposition of the two diagrams representing the translucent words: if an opaque point covers a translucent point, the resulting point is opaque, see Figure~\ref{fig:exampleComposition-diagram}. The diagrams make apparent that this is just the opposite of the composition of corresponding inclusion maps of ordered bisets (the translucent word $t$ corresponds to the unique strictly increasing, left-right structure-preserving map $\restriction{\alpha_t}{[t]_0}\to \alpha_t$).

\begin{Proposition}\label{prop:compose-iota}
 In terms of the morphisms of finite ordered bisets defined in Remark~{\rm \ref{rk:translucent-words-and-morphisms-of-ordered-bisets}}, it holds that
 \[\iota_t\circ \iota_s=\iota_{s\circ t}\]
 for all composable translucent words $s$, $t$.
 In particular, the composition of translucent words is associative.
\end{Proposition}

\begin{proof}
 Both maps, $\iota_t\circ \iota_s$ and $\iota_{s\circ t}$, are strictly increasing maps from $|[s]_0|=|[s\circ t]_0|$ to $|t|$. A~strictly increasing map between finite totally ordered sets is uniquely determined by its image.
 Let $[t]_0=\{i_1<\cdots<i_n\}\subset \interval{|t|}$ and $[s]_0=\{j_1<\cdots<j_m\}\subset\interval{|s|}=\interval{n}$, so that $\iota_t(k)=i_k$ and $\iota_s(\ell)=j_\ell$. Then
 \begin{align*}
 \iota_{s\circ t}(\interval{m})=[s\circ t]_0=\{i_k\colon b_s(k)=0, k=1,\dots,m\}
 =\iota_t([s]_0)=\iota_t(\iota_s(\interval{m})
 \end{align*}
 because by definition of the composition $b_{s\circ t}(i)=0$ if and only if $b_t(i)=0$ and $b_s(k)=0$ for the unique $k$ with $i=i_k$.
\end{proof}

With this composition, translucent words form the morphisms of a category, also denoted~$\SWc$, whose objects are words in $\{\LL,\RR\}^\star$.
The empty word is a (the) terminal object, denoted $1$, with the unique morphism $(\alpha, \mathbf{1})\colon \alpha \to 1$.
In the category~$\SWc$ the identity morphism acting on a~word $\alpha$ corresponds to the translucent word $(\alpha, \mathbf{0})$, where $\mathbf{0}$ denotes the Boolean vector with the same length as~$\alpha$ and only zeros as entries.

\begin{Remark}
 Be aware that $t\in \SWc$ always means that $t$ is a translucent word, i.e., $t\in \operatorname{Mor}(\SWc)$. This is in contrast to the convention we use for most other categories; for another category $\mathcal C$, we write $c\in \mathcal C$ if $c\in \operatorname{Obj}(\mathcal C)$.
\end{Remark}

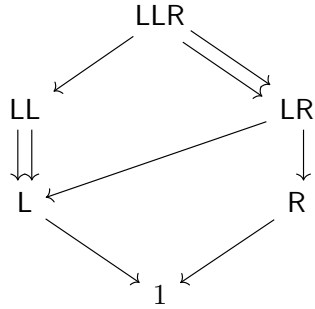
\begin{figure}[!ht]	\centering
	\begin{tikzcd}
		& \LL\LL\RR \arrow[rd, shift left] \arrow[rd, shift right] \arrow[ld] & \\
		\LL\LL \arrow[d, shift left] \arrow[d, shift right] & & \LL\RR \arrow[d, shift left] \arrow[lld] \\
		\LL \arrow[rd] & & \RR \arrow[ld] \\
		& 1 &
	\end{tikzcd}
	\caption{Objects and arrows in the category $\SWc$. The two inner-faces are commutative. While $\SWc$ can be drawn as a directed fat-graph, is it is not equivalent to the free category generated by this graph.}\label{fig:objectsarrows}
\end{figure}

Following the terminology introduced in \cite{content1980categories} (see also \cite{dur1986mobius,leinster2012notions}), the category $\SWc$ is a \emph{M\"obius category}, which means that any arrow $f$ of $\SWc$ is \emph{finitely decomposable} in the sense that
\[
	\big|\bigl\{(f',f'')\colon f'' \circ f' = f \bigr\}\big| < +\infty.
\]
This type of category is the proper algebraic framework for M\"obius inversion, see \cite{content1980categories}.

\subsection{Incidence category with coefficient in Set}\label{subsec:monoidalindexedmord-set}

\begin{Definition}
 The category $\Set^{\SWc}$, whose
\begin{itemize}\itemsep=0pt
 \item objects are the {\it $\SWc$-indexed collections}
	$\mathcal{X}=(\mathcal{X}(t))_{t \in \SWc}$, $\mathcal{X}(t) \in \Set$.
	\item morphisms from $\mathcal{X}$ to $\mathcal{Y}$
are the families of maps ($\Set$ morphisms)
\[f_{t}\colon \ \mathcal{X}(t) \to \mathcal{Y}(t),\qquad t \in \SWc,\]
\end{itemize}
is called {\it incidence category of $\SWc$ with coefficients in~$\Set$}.
\end{Definition}

\begin{Remark} The terminology incidence category we use in this paper is not to be confused with occurrences in the literature with a different meaning, see~\cite{szczesny2011incidence}.
\end{Remark}

Recall that the category $\Set $ is endowed with a monoidal structure given by the Cartesian product. As unit object we fix a one-point set $\{\Singleton\}$.
$\Set$ is also endowed with a categorical coproduct, the disjoint union of sets, denoted $\displaystyle\sqcup$ in what follows.

\begin{Definition}[monoidal composition of $\SWc$-indexed collections]
	For two $\SWc$-indexed collections of sets $\mathcal{X}$ and $\mathcal{Y}$,
	their {\it composition} is defined as the $\SWc$-indexed collection
	$\mathcal{X}\vertset \mathcal{Y}$, defined~as\looseness=-1
	\begin{equation*}
		(\mathcal{X}\vertset \mathcal{Y})(t)
		:= \bigsqcup_{\substack{r,s\in\SWc\\
		t = r \circ s}}
		\mathcal{X}(r) \times \mathcal{Y}(s), \qquad t \in \SWc.
	\end{equation*}
\end{Definition}
Define the $\SWc$-indexed collection $\mathcal{E}^\vertset$ by setting for any $t \in \SWc$,
\[\mathcal{E}^\vertset(t)=
\begin{cases}
 \varnothing &\text{if $b_t\neq \mathbf{0}$},\\
 \{\Singleton \} & \text{if $b_t= \mathbf{0}$}.
\end{cases}
\]
The collection $\mathcal{E}^\vertset$ acts as left and right unit for the composition $\vertset$, which means that for any $t \in \SWc$,
\begin{gather*}
	(\mathcal{X}\vertset \mathcal{E}^\vertset)(t) = \bigsqcup_{\substack{s\in \SWc \\ t=s\circ (\alpha_t,\mathbf{0}_{|t|})}} \mathcal{X}(s)\times \{\Singleton \} = \mathcal{X}(t) \times \{\Singleton \},
\\
	(\mathcal{E}^\vertset \vertset \mathcal{X})(t) = \bigsqcup_{\substack{s\in \SWc \\ t = (\restriction{\alpha_t}{[s]_0},\mathbf{0}_{|[s]_0|})\circ s}} \{\Singleton \} \times \mathcal{X}(s) = \{\Singleton \}\times \mathcal{X}(t).
\end{gather*}
Associativity for $\vertset$ follows from the associativity
of the monoidal product $\times$ on $\Set $, the coproduct $\sqcup$ on $\Set $
and the associativity of composition of morphisms in the category $\SWc$.
The composition $\vertset$ hence yields a monoidal structure on $\Set^{\SWc}$.

\subsection{Incidence category with coefficients in Vect}
\label{subsec:monoidalindexedmord-vect}

We now shift to the linear setting. The category of (complex) vector spaces is denoted $\Vect$ or $\Vect_\CC $ if we want to emphasize the field of scalars. A $\SWc$-indexed collection $\mathcal V=\mathcal{V}(t)_{t\in\SWc}$ of vector spaces $\mathcal{V}(t)\in \Vect$ is also called {\it linear $\SWc$-indexed collection}.

\begin{Definition}
 The category $\Vect^{\SWc}$, whose
\begin{itemize}\itemsep=0pt
 \item objects are the {\it $\SWc$-indexed collections of $\CC $-vector spaces}
	$(\mathcal{V}(t))_{t \in \SWc}$, $\mathcal{V}(t) \in \Vect$,
	\item morphisms from $\mathcal{V}$ to $\mathcal{W}$
are the families of linear maps ($\Vect$-morphisms)
\[f_{t}\colon \ \mathcal{V}(t) \to \mathcal{W}(t),\qquad t \in \SWc,\]
\end{itemize}
is called {\it incidence category of $\SWc$ with coefficients in $\Vect$}.
\end{Definition}

\begin{Remark}
 It is sometimes opportune to identify a collection $\mathcal V\in \Vect^\SWc$ with the $\SWc$-graded vector space $\bigoplus_{t\in\SWc} \mathcal V(t)$ and a morphism $f\colon \mathcal V\to \mathcal W$ with the grade-preserving map $f(v):=f_t(v)$ for $v\in\mathcal V(t)$.
\end{Remark}

All linear $\SWc$-collections in this work are spanned by explicitly given $\SWc$-indexed collection of sets. For a set $X$, we denote by \[\linearspan{X}=\bigoplus_{x\in X}\linearspan{x}\] the vector space formally spanned by $X$.
For a $\SWc$-collection $\mathcal{X}$ of sets, we denote by $\linearspan{\mathcal{X}}$ the linear $\SWc$-indexed collection spanned by $\mathcal{X}$,
\begin{equation*}
 (\linearspan{\mathcal{X}})(t)
 :=\linearspan{(\mathcal X(t))}
	=\bigoplus_{x\in\mathcal{X}(t)}\linearspan{x},\qquad t \in \SWc.
\end{equation*}
The {\it dual collection} of $\linearspan{\mathcal{X}}$ is the linear $\SWc$-indexed collection $\linearspan{\mathcal{X}}^{\ast}$ defined by
\[ \linearspan{\mathcal{X}}^{\ast}(t) := \left(\linearspan{\mathcal{X}}(t)\right)^{\ast},\qquad t \in \SWc .\]

There is a natural choice for a monoidal structure on the category of such linear collections.
In the following definition, we denote by $\otimes_{\CC}$ the tensor product on the category of complex vector spaces.
\begin{Definition}[``vertical'' composition of linear $\SWc$-indexed collections]
	For two linear $\SWc$-indexed collections $\mathcal{V}$ and $\mathcal{W}$ define the {\it vertical composition} $\mathcal{V}\overt \mathcal{W}$
	as the $\SWc$-indexed collection given by
	\begin{equation*}
		(\mathcal{V}\overt \mathcal{W})(t)
		:=
		\bigoplus_{t=r\circ s}
		\mathcal{V}(r)
		\otimes_{\CC}
		\mathcal{W}(s), \qquad t\in\SWc.
	\end{equation*}
\end{Definition}
Again, it follows from a straightforward verification that $\overt$ is monoidal, with unit $\uvert$ defined by
\begin{equation*}
	\uvert(t) =
	\begin{cases}
	 \linearspan{\{\Singleton\}} &\text{it $b_t=\mathbf{0}$},\\
	 \{0\} &\text{if $b_t \neq \mathbf{0}$}.
	\end{cases}
\end{equation*}
Note that ${\uvert}(1)=\linearspan{\{\Singleton\}}$ (1 the empty translucent word). Also, note that we use the same symbol for the unit in the linear and non-linear setting.

\begin{Notation} \label{notation:overt}
 We may write $x\overt y$ instead of $x\otimes y$ if $x\in\mathcal V(r)$ and $y\in\mathcal W(s)$ for composable translucent words $r$ and $s$.
\end{Notation}

\begin{Remark}\label{rk:span-and-vertical-composition}
 Clearly, we have for all $\mathcal{X},\mathcal{Y}\in \Set^\SWc$ a canonical isomorphism
 \[\linearspan {(\mathcal{X}\vertset \mathcal{Y})}\cong \linearspan{\mathcal{X}}\overt \linearspan{\mathcal {Y}},\]
 which will be used to identify the two.
\end{Remark}

\begin{Definition}
 A {\it $($linear$)$ $\SWc$-monoid} is a monoid in the monoidal category
 $\big(\Set^\SWc,\overt,\uvert\big)$ \big(or $\big(\mathbb \Vect^\SWc,\overt,\uvert\big)$\big), i.e., a triple $\big(\mathcal V,\vp,\eta^{\overt}\big)$ where $\mathcal V$ is a (linear) $\SWc$-indexed collection,
	$\vp\colon \mathcal V \overt \mathcal V \to \mathcal V$ and $\eta^{\overt}\colon\uvert \rightarrow \mathcal V$ are morphisms in $\Set^\SWc$ \big(or $\Vect^\SWc$\big) satisfying the unitality and associativity constraints:
	\begin{equation*}
		\vp \big(\eta^{\overt} \overt \mathrm{id}\big) = \vp \big(\mathrm{id}\overt \eta^{\overt} \big) = \mathrm{id},\qquad \vp \big(\vp \overt \mathrm{id}\big) = \vp\big(\mathrm{id}\overt \vp\big).
	\end{equation*}
\end{Definition}

\begin{Remark}
Let $\big(\mathcal{X},\mathsf{m}^\vertset_{\mathcal X},\eta^\vertset_{\mathcal X}\big)$ be a $\SWc$-monoid.
By Remark~\ref{rk:span-and-vertical-composition} it is clear that the $\Set$-monoidal product $\mathsf{m}_{\mathcal X}^\vertset\colon \mathcal{X}\vertset \mathcal{X} \to \mathcal{X}$ extends to a morphism of linear $\SWc$-indexed collections,
\begin{equation*}
	\vp_{\mathbb C\mathcal X}\colon \ \linearspan{\mathcal{X}}\overt \linearspan{\mathcal{X}} \to \linearspan{\mathcal{X}},
\end{equation*}
thus turning $\linearspan{\mathcal{X}}$ into a linear $\SWc$-monoid with unit $\eta^{\overt}_{\mathbb C\mathcal X}$ the linear extension of $\eta^\vertset_{\mathcal X}$.
\end{Remark}

\begin{Definition}
	A linear {\it $\SWc$-comonoid} is a triple $\big(F,\vcp,\varepsilon\big)$ where $F$ is a linear $\SWc$-indexed collection, $\vcp\colon F \to F \overt F$ and $\varepsilon\colon F \rightarrow \uvert$ satisfy the unitaly and coassociativity constraints:
	\begin{equation*}
		(\varepsilon \overt \mathrm{id}) \vcp = (\mathrm{id}\overt \varepsilon)\vcp = \mathrm{id},\qquad \big(\vcp \overt \mathrm{id}\big)\vcp = \big(\mathrm{id}\overt \vcp\big)\vcp.
	\end{equation*}
\end{Definition}
\begin{Remark}
	We only define \emph{linear} $\SWc$-comonoids,
	since we are interested in comonoids with an unshuffle splitting of the coproduct $\vcp$, which is a meaningless notion in $\Set$.
\end{Remark}
Of interest to us will be two linear monoids, which by the following recipe yield two linear $\SWc$-comonoids.
Assume for the rest of this section that $\mathcal X$ is a $\SWc$-monoid 
in which each set~$\mathcal X(t)$ is finite, $t \in \SWc$.
Thanks to this finiteness assumption, first,
\begin{equation*}
	(\linearspan{\mathcal{X}}\overt \linearspan{\mathcal{X}})^{\ast} \simeq \linearspan{\mathcal{X}}^{\ast}\overt \linearspan{\mathcal{X}}^{\ast}.
\end{equation*}
Second, the basis $\mathcal{X}$ yields an isomorphism between the collection $\linearspan{\mathcal{X}}$ and its dual $\linearspan{\mathcal{X}}^{\ast}$ by sending $\mathcal{X}$ to its dual basis. The algebraic dual of $\vp$,
\begin{equation*}
	\big(\vp\big)^{*}\colon \ \linearspan{\mathcal{X}}^{\ast} \to (\linearspan{\mathcal{X}}\overt \linearspan{\mathcal{X}})^{\ast} \simeq \linearspan{\mathcal{X}}^{\ast} \overt \linearspan{\mathcal{X}}^{\ast}
\end{equation*}
yields the structure of a comonoid on $\linearspan{\mathcal{X}}^{\ast} \simeq \linearspan{\mathcal{X}}$. More concretely, under the named isomorphism, $\big(\vp\big)^{*}$ and $\eta^\ast$ correspond to the structural morphisms $\vcp\colon \linearspan{\mathcal{X}}\to \linearspan{\mathcal{X}}\overt \linearspan{\mathcal{X}}$ and $\varepsilon\colon \linearspan{\mathcal{X}}\to \uvert$ defined by
\begin{equation*}
	\vcp_t(x) = \sum_{\substack{y,z\in \mathcal X(t) \\x = \vp(y\overt z)}} y \overt z,\qquad\varepsilon_t(x) =\begin{cases}
	 \Singleton & \text{if $b_t=\mathbf{0}$},\\
	 0 & \text{if $b_t\neq\mathbf{0}$}
	\end{cases}
\end{equation*}
for any $x\in \mathcal{X}(t)$.

\subsection{Horizontal semigroupal product}

We introduce in this subsection a second (non-unital) tensor product on the incidence category $\icsw$, which we name \emph{horizontal tensor product}. The horizontal tensor product turns out to be compatible with the previously introduced vertical tensor product $\overt$ in a way reminiscent of \emph{$2$-monoidal categories} or \emph{duoidal categories}. The reader is directed to the monograph \cite{aguiar2010monoidal} for a~detailed account of this notion from higher category theory.

In our case, however, the horizontal tensor product \emph{lacks existence of a unit} (we use the terminology semigroupal in the sequel), we, therefore, refrain from defining $2$-monoidal categories here.\looseness=-1

Nevertheless, we will show the existence of a functor, called \emph{four-points exchange relation} or \emph{interchange relations} satisfying the natural associativity constraints. This implies a rather rich structure on the category of horizontal semigroups (objects equipped with a binary associative multiplication morphism) that we leverage in the forthcoming sections.

The horizontal semigroupal product is built on splitting a translucent word $t$ at a translucent point.
More precisely, for $i\in [t]_0$ we denote the restrictions of $t$ to the $\stdl$-intervals of positions before or after $i$ with respect to the necklace order by
\begin{align*}
 t^{\stdl i}:=\restriction{t}{\{ j\colon j \stdl i \}},\qquad
 t^{\stdg i}:=\restriction{t}{\{ j\colon j \stdg i \}}.
\end{align*}
Recall that the unique endomorphism of the empty word $\varnothing\in\{\LL,\RR\}^{\star}$ is the empty translucent word, denoted $1$.

\begin{Example}
 For $t = ({\sf L R L L}, 0101)$, the necklace order is $1\stdl 3\stdl 4\stdl 2$, therefore
	\begin{align*}
		 t^{\stdl 1} = 1, \qquad t^{\stdg 1} = ({\sf RLL}, 101),
		 \qquad
		 t^{\stdl 3} = ({\sf L}, 0),\qquad t^{\stdg 3} = ({\sf RL}, 11).
	\end{align*}
\end{Example}
\begin{Definition}[``horizontal'' composition of linear $\SWc$-indexed collections]
	\label{def:ominus}
	For two $\SWc$-indexed collections $\mathcal{V}$ and $\mathcal{W}$, we
	define the {\it horizontal composition} $\mathcal{V} \ominus \mathcal{W}$
	as the $\SWc$-indexed collection\looseness=-1
	\begin{equation}\label{eqn:defominus}
		(\mathcal{V} \ominus \mathcal{W})(t) = \bigoplus_{i \in [t]_0} \mathcal{V}\big(t^{\stdl i}\big)\otimes \mathcal{W}\big(t^{\stdg i}\big).
 \end{equation}

Analogously, for $\SWc$-indexed collections of sets,
 \[(\mathcal{X} \ominus \mathcal{Y})(t) = \bigsqcup_{i \in [t]_0} \mathcal{X}\big(t^{\stdl i}\big)\times \mathcal{Y}\big(t^{\stdg i}\big),\]
 so that
 $\mathbb C (\mathcal X\ominus \mathcal Y)=\mathbb C \mathcal X \ominus \mathbb C\mathcal Y$.
\end{Definition}
It is straightforward to verify the associativity constraint for $\ominus$ for the (trivial) associator
\[
	\alpha^{\ominus}_{\mathcal{U},\mathcal{V},\mathcal W}\colon \ (\mathcal{U}\ominus \mathcal{V})\ominus \mathcal W \to \mathcal{U} \ominus (\mathcal{V}\ominus \mathcal W),
	\qquad (u\otimes v) \otimes w \mapsto u \otimes (v\otimes w),
\]
where $u \in \mathcal{U}(t^{\stdl i})$, $v \in \mathcal{V}(t^{\stdg i,\stdl j})$ and $w\in \mathcal{W}(t^{\stdg j})$ with $t^{\stdg i,\stdl j}$ the translucent word $t$ restricted to the interval $\rrbracket i,j\llbracket_{\stdl}$.
Here we use in a crucial way that this operation of cutting a translucent word $t$ at a certain position $i$ and forgetting about $t(i)$ is associative.

In fact, one has (omitting the associator $\alpha$)
\begin{align}\label{eq:horizontal-unbracketed}
	\mathcal{U}\ominus\mathcal{V}\ominus\mathcal{W}(t) = \bigoplus_{\substack{i,j\in [t]_0 \\i \stdl j }} \mathcal{U}\big(t^{\stdl i}\big)\ominus \mathcal{V}\big(t^{\stdg i,\stdl j}\big)\ominus \mathcal W\big(t^{\stdg j}\big).
\end{align}

The bifunctor $\ominus\colon \icsw{}^{\times 2}\to \icsw$ does not yield a monoidal product, since it cannot be endowed
with a unit, as we see below.%
\footnote{Being monoidal is a property, not a structure,
	compare \cite[Remark 2.2.9]{etingof2016tensor}.}
Therefore, we call it a {\it semigroupal product}
(a ``monoidal product without unit'').
Correspondingly, the equivalent of a monoid in such a category
will be called a {\it semigroup}.
We denote by $\sgominus$ the category of $\ominus$-semigroups
(with obvious morphisms).
\begin{Remark}
	If $[t]_0=\varnothing$, then the sum on the right side of equation~\eqref{eqn:defominus} is empty, thus equal to $\{0\}$. In particular, $(\mathcal{V}\ominus \mathcal{W})(\alpha,\mathbf{1}_{|\alpha|}) = \{0\}$ for any word $\alpha\in\{\LL,\RR\}^{\star}$. This also yields that $\ominus$ has no unit.
\end{Remark}
\begin{Remark}
	It is possible to design a horizontal tensor product $\ominus^{\prime}$ \emph{with unit} by
	\[
		(\mathcal{V} \ominus^{\prime} \mathcal{W})(t) = \sum_{i\in [b]_0} \mathcal{V}\big(t^{\stdl i}\big) \otimes \mathcal{W}\big(t^{\stdgeqs i}\big),
	\]
	where $t^{\stdgeqs i}$ is the restriction of $t$ to $\{x\colon x\stdgeq i\}$. This tensor product, while unital, seems less natural than $\ominus$ because of the choice on which side to include~$i$.
\end{Remark}

\begin{figure}[!ht]	\centering
	\includegraphics[scale=0.75]{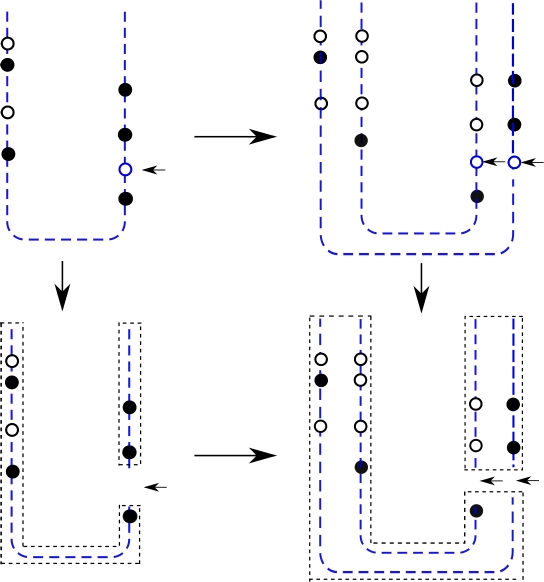}
	\caption{Four point exchange relation. The top left translucent word $t$
	is deconcatenated at the position~$i$ marked by the arrow, yielding two diagrams $t^{\stdl i}$ and $t^{\stdg i}$ shown in the bottom left. In the top right, we see two translucent words $r$ and $s$ with $t=r\circ s$ (of course, those are not uniquely determined by~$t$, the diagram only shows one possible decomposition of~$t$).
	The bottom right, which shows $r_-$, $r_+$, $s_-$, $s_+$, can be either obtained by decomposing the bottom left or splitting the top right diagrams, exemplifying that splitting and decomposition commute as long as the position $j$ at which $r$ is to be split corresponds to the position $i$ at which $s$ and $t$ are split. On the other hand, these diagrams also show how~$r$,~$s$ and~$j$ can be reconstructed from $r_-$, $r_+$, $s_-$, $s_+$ and $i$.}\label{fig:deuxmonoids}
\end{figure}
Our next aim is to introduce a natural transformation
\[R \colon \ \ominus \circ (\overt \times \overt) \to \overt \circ (\ominus \times \ominus) \circ \tau_{2,3},\]
where $\tau_{2,3}\colon \big(\icsw\big)^{\times 4}\to \big(\icsw\big)^{\times 4}$, $(\mathcal{V}_1,\mathcal{V}_2,\mathcal{V}_3,\mathcal{V}_4)\mapsto (\mathcal{V}_1,\mathcal{V}_3,\mathcal{V}_2,\mathcal{V}_4)$.

Given an element
\begin{align}
 (R_-\overt S_-)\ominus (R_+\overt S_+)&\in (\mathcal{V}_1 \overt \mathcal{V}_2) \ominus (\mathcal{V}_3 \overt \mathcal{V}_4)(t)\notag\\
 &=\bigoplus_{\substack{i\in[t]_0\\r_-\circ s_-=t^{\stdl i}\\r_+\circ s_+=t^{\stdl i}}}(\mathcal{V}_1(r_-) \otimes \mathcal{V}_2(s_-)) \otimes (\mathcal{V}_3(r_+) \otimes \mathcal{V}_4(s_+))\label{eq:exchange-domain}
\end{align}
we want to obtain an element{\samepage
\begin{align}
 (R_-\ominus R_+)\overt (S_-\ominus S_+)&\in (\mathcal{V}_1 \ominus \mathcal{V}_3) \overt (\mathcal{V}_2 \ominus \mathcal{V}_4)(t)\notag\\
 &=\bigoplus_{\substack{r\circ s=t\\ j\in [r]_0\\ i\in [s]_0}}(\mathcal{V}_1(r^{\stdl j}) \otimes \mathcal{V}_3(r^{\stdg j})) \otimes (\mathcal{V}_2(s^{\stdl i}) \otimes \mathcal{V}_4(s^{\stdg i})).\label{eq:exchange-image}
\end{align}}

\noindent
To this end, we will construct a map which maps each tuple $(i, r_-,r_+,s_-,s_+)$ with $i\in [t]_0$, $r_-\circ s_-=t^{\stdl i}$, $r_+\circ s_+=t^{\stdg i}$ to a tuple $(r,s, i',i) $ with
\begin{gather}
r\circ s=t,\quad\! i'\in[r]_0,\quad\! i\in [s]_0,\quad\! \text{and} \quad\!
 r_-=r^{\stdl i'},\quad\! s_-=s^{\stdl i},\quad\! r_+=r^{\stdg i'},\quad\! s_+=s^{\stdg i}\label{eq:conditions-for-construction}
\end{gather}
(we will also see that the map is uniquely determined by those conditions). This is the content of Construction~\ref{constrctone} below (see Figure~\ref{fig:deuxmonoids} for a graphical example).
 With such a map in hand, we can then simply map an element $(R_-\overt S_-)\ominus (R_+\overt S_+) \in(\mathcal{V}_1(r_-) \otimes \mathcal{V}_2(s_-)) \otimes (\mathcal{V}_3(r_+) \otimes \mathcal{V}_4(s_+))\subset (\mathcal{V}_1 \overt \mathcal{V}_2) \ominus (\mathcal{V}_3 \overt \mathcal{V}_4)(t) $ which lives in one of the direct summands of the direct sum in \eqref{eq:exchange-domain} to the element $(R_-\ominus R_+)\overt (S_-\ominus S_+)$ in the corresponding summand $(\mathcal{V}_1(r^{\stdl i'}) \otimes \mathcal{V}_3(r^{\stdg i'})) \otimes (\mathcal{V}_2(s^{\stdl i}) \otimes \mathcal{V}_4(s^{\stdg i}))\subset (\mathcal{V}_1 \ominus \mathcal{V}_3) \overt (\mathcal{V}_2 \ominus \mathcal{V}_4)(t) $ in \eqref{eq:exchange-image}.

Although we will not use those facts explicitly, one can easily observe that the map obtained from Construction~\ref{constrctone} is injective but not surjective and, hence, the natural transformation $R$ is composed of monomorphisms but not of isomorphisms.

\begin{Construction}\label{constrctone}
Pick a translucent word $t \in \SWc$. Consider an integer $i\in [t]_0$ and two pairs of composable translucent words
$r_{-}$, $s_{-}$ and $r_{+}$, $s_{+}$
such that{\samepage
\begin{equation*}
	r_{-}\circ s_{-} = t^{\stdl i},\qquad r_{+}\circ s_{+} = t^{\stdg i}.
\end{equation*}
Note that $r_-$, $r_+$ are determined by the choice of $s_-$, $s_+$.}

We define the translucent word $s \in \SWc$ as the unique translucent word of length $|t|$ fulfilling
\[{s}^{\stdl i}=s_-,\qquad {s}^{\stdg i}=s_+ \qquad \text{and}\qquad s(i)=t(i),\]
and the translucent word $r\in \SWc$ as $r:=\restriction{t}{[s]_0}$.
Clearly, $t=r\circ s$, and there is no other choice for $r$ and $s$ if \eqref{eq:conditions-for-construction} is to be fulfilled. Note that we used the original $i$. It remains to find a $i'\in[r]_0$ with $r_-=r^{\stdl i'}$, $r_+=r^{\stdg i'}$. Recall that composability implies $|r|=|[s]_0|$ and let $[s]_0=\{k_1<\dots<k_{|r|}\}$. The (unique) $i'$ with $i=k_{i'}$ has the desired properties. By definition of the composition, $r(i')=r\circ s(k_{i'})=t(i)$. Let us take for granted for the moment that $=(\restriction{t}{[s]_0})^{\stdl i'}
 =\restriction{t^{\stdl i}}{[s^{\stdl i}]_0}$, which is quite clear from the diagrammatic representation and will be proved below.
Then we can conclude from $s^{\stdl i}=s_-$ and $r_-\circ s_-=t^{\stdl i}$ that
\[
 r^{\stdl i'}
 =(\restriction{t}{[s]_0})^{\stdl i'}
 =\restriction{t^{\stdl i}}{[s^{\stdl i}]_0}=\restriction{t^{\stdl i}}{[s_-]_0}=r_-,\] and analogously we obtain \(r^{\stdg i'}=r_+\). In summary, the constructed tuple $(r,s,i',i)$ fulfills~\eqref{eq:conditions-for-construction}, as needed.
\end{Construction}

To prove the claim that $(\restriction{t}{[s]_0})^{\stdl i'}
=\restriction{t^{\stdl i}}{[s^{\stdl i}]_0}$, we need to some more notation. For any ordered set $A$, let $\varphi_A\colon A\to\interval{|A|}$ denote the unique monotone bijection. Then it is easy to check that for arbitrary subsets $A,B\subset \interval{|t|}$, one has
\[\restriction{(\restriction{t}{A})}{\varphi_A(B\cap A)}=\restriction{t}{A\cap B}=\restriction{(\restriction{t}{B})}{\varphi_B(A\cap B)}.\]
In the special case, recall the notation $[s]_0=\{k_1<\dots <k_{|r|}\}$ and note that $\varphi_{[s]_0}(i)=i'$. Then $\ell\stdl j$ is equivalent to $k_\ell\stdl k_j$ and we find
\begin{gather*}
(\restriction{t}{[s]_0})^{\stdl i'}\!=\restriction{(\restriction{t}{[s]_0})}{\{\ell:k_\ell\stdl i\}}=\restriction{(\restriction{t}{[s]_0})}{\varphi_{[s]_0}([s]_0\cap \{k\stdl i\})}=\restriction{(\restriction{t}{\{k\stdl i\}})}{\varphi_{\{k\stdl i\}}([s]_0\cap \{k\stdl i\})}=\restriction{t^{\stdl i}}{[s^{\stdl i}]_0}.
\end{gather*}

We are now in the position to define the natural transformation $R$ as follows.

Pick elements
\[R_-\in \mathcal{V}_1({r_-}),\qquad S_- \in \mathcal{V}_2(s_-),\qquad R_+ \in \mathcal{V}_3(r_+),\qquad S_+ \in \mathcal{V}_4 (s_+),\] so that
\[
	(R_-\overt S_-) \in (\mathcal{V}_1 \overt \mathcal{V}_2)\big(t^{\stdl i}\big),\qquad (R_+\overt S_+) \in (\mathcal{V}_3\overt \mathcal{V}_4)\big(t^{\stdg i}\big)
\]
and set, with $(r,s,i',i)$ obtained from Construction~\ref{constrctone},
\begin{gather*}
	R_{\mathcal{V}_1,\mathcal{V}_2,\mathcal{V}_3,\mathcal{V}_4}((R_-\overt S_-) \ominus (R_+ \overt S_+)):=
	 (R_- \ominus R_+) \overt (S_- \ominus S_+)\\
\qquad{} \in \big(\mathcal{V}_1\big(r^{\stdl i'}\big) \otimes \mathcal{V}_3\big(r^{\stdg i'}\big)\big) \otimes \big(\mathcal{V}_2\big(s^{\stdl i}\big) \otimes \mathcal{V}_4\big(s^{\stdg i}\big)\big)\subset (\mathcal{V}_1 \ominus \mathcal{V}_3) \overt (\mathcal{V}_2 \ominus \mathcal{V}_4)(t).
\end{gather*}

We now prove the main proposition of this section.
\begin{Proposition}\label{prop:exchange}
	The morphisms $R_{\mathcal{V}_1,\mathcal{V}_2,\mathcal{V}_3,\mathcal{V}_4}$, $\mathcal{V}_i \in \icsw$ yield a natural transformation
	\[
		R\colon \ \ominus \circ (\overt \times \overt) \to \overt \circ (\ominus \times \ominus) \circ \tau_{2,3}
	\]
	satisfying the associativity constraints that the diagrams in Figures~{\rm \ref{fig:assI}} and~{\rm \ref{fig:assII}} commute.

	\begin{figure}[!ht]\centering
		\begin{tikzcd}
			((A\overt B) \ominus (C\overt D))\ominus(E\overt F) \arrow[r, "\alpha"] \arrow[d, "R\ominus \mathrm{id}"] & (A\overt B)\ominus ((C\overt D)\ominus (E\overt F)) \arrow[d, "\mathrm{id}\overt R"] \\
			((A\ominus C) \overt (B\ominus D))\ominus(E\overt F) \arrow[d, "R"] & (A\overt B)\ominus ((C\ominus E)\overt (D\ominus F)) \arrow[d, "R"] \\
			((A\ominus C) \ominus E)\overt ((B\ominus D)\ominus F) \arrow[r, "\alpha^{\ominus}\overt\alpha^{\ominus}"] & (A\ominus (C\ominus E))\overt (B\ominus (D\ominus F))
		\end{tikzcd}
		\caption{First associativity constraint.}\label{fig:assI}
	\end{figure}

	\begin{figure}[!ht]\centering
		\begin{tikzcd}
			((A\overt B)\overt C)\ominus ((D\overt E)\overt F) \arrow[d, "R"] \arrow[r, "\alpha^{\ominus}\ominus\alpha^{\ominus}"] & (A\overt (B\overt C))\ominus(D\overt (E\overt F)) \arrow[d, "R"] \\
			((A\overt B)\ominus (D\overt E))\overt(C\ominus F) \arrow[d, "R\overt \mathrm{id}"] & (A\ominus D)\overt ((B\overt C)\ominus(E\overt F)) \arrow[d, "\mathrm{id}\overt R"] \\
			((A\ominus D)\overt (B\ominus E))\overt(C\ominus F) \arrow[r, "\alpha^{\ominus}"] & (A\ominus D)\overt ((B\ominus E)\overt(C\ominus F))
		\end{tikzcd}
		\caption{Second associativity constraint.}\label{fig:assII}
	\end{figure}
\end{Proposition}

\begin{figure}[!htb]\centering
	\includegraphics[scale=0.55]{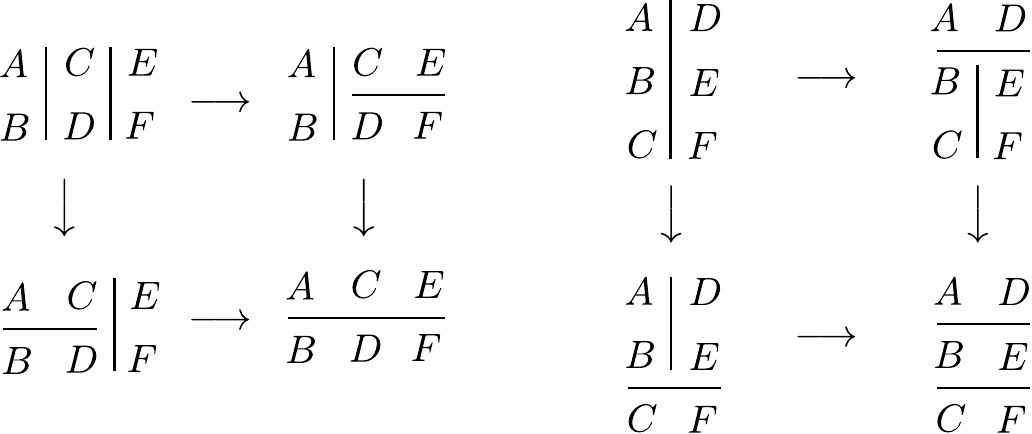}
	\caption{An artistic view on the associativity constraints.}\label{fig:associativity}
\end{figure}

One can depict the commutative diagram in Figures~\ref{fig:assI} and~\ref{fig:assII} using horizontal and vertical bars to indicate which of the horizontal product or vertical tensor product takes precedence. Omitting the (necessary) bracketing, one obtains Figure~\ref{fig:associativity}.

\begin{Remark} We insist upon the fact that the natural transformation $R$ is not an isomorphism. In fact, with the notations introduced in Construction \ref{constrctone}, the two integers $i$ and $i^{\prime}$ satisfy to the constraint $\varphi_{[s]_0}(i)=i^{\prime}$. In less technical terms, $r$ and $s$ are split at the same positions, when the letters of $r$ are indexed by the positions they have in the translucent letters of $t$ (placeholders). This indicates that $R$ is in fact not surjective, but it is easy to check that $R$ is injective.
\end{Remark}
\begin{proof}
 Let $t$ be a translucent word. Pick two integers $i \stdl j$. The notation $t^{\stdg i,\stdl j}=\restriction{t}{\rrbracket i,j\llbracket_\stdl}$ has been introduced at the beginning of the section. We pick three pairs of composable translucent words
 \[
 (r_-,s_-),\qquad (r_0,s_0),\qquad (r_+,s_+)
 \]
 such that
 \[
 r_-\circ s_-=t^{\stdl i},\qquad r_0\circ s_0 = t^{\stdg i,\stdl j}, \qquad \text{and} \qquad r_+\circ s_+ = t^{\stdg j}.
 \]
 First, note that instead of iterating, we can modify Construction~\ref{constrctone} to work for the ``unbracketet'' version of the threefold horizontal product \eqref{eq:horizontal-unbracketed}. As in the two-fold case, the only ambiguous point is where the splitting of $r$ into its given parts (here: $r_-$, $r_0$, $r_+$) takes place, and the modified construction clearly gives $i':=\varphi_{[s]_0}(i)$ and $j':=\varphi_{[s]_0}(j)$ for those positions. To prove the first associativity constraint (Figure~\ref{fig:assI}), we show that both iterations agree with the natural transformation $R_3\colon (A\overt B)\ominus (C\overt D)\ominus (E\overt F)\to (A\ominus C\ominus E)\overt (B\ominus D\ominus F)$ induced by the modified threefold construction.
 When we iterate Construction~\ref{constrctone} to find $r$ and $s$, the only problem is to determine and compare the positions where the words are to be split, and these splitting points determine the obtained natural transformation.
 We illustrate the splitting positions in a table:
 \begin{center}
 \begin{tabular}{|c|c|c|}\hline
 full word& factors & split position\\\hline\hline
 \multirow{2}{*}{$ t^{\stdl j} $}&$r_{-,0}$&$\varphi_{[s_{-,0}]_0}(\varphi_{\{k\stdl j\}}(i))$\\\cline{2-3}
 &$s_{-,0}$&$\varphi_{\{k\stdl j\}}(i)$\\\hline\hline
 \multirow{2}{*}{$t$}&$r$&$\varphi_{[s]_0}(j)$\\\cline{2-3}
 &$s$&$j$\\\hline

 \end{tabular}
 \end{center}
 If we want to compare with the threefold splitting, we have to map the positions on the first two rows to the corresponding positions in $\interval{|[s]_0|}$ and $\interval{|t|}$ via the maps $\varphi_{\{k\stdl j'\}}^{-1}$ and $\varphi_{\{k\stdl j\}}^{-1}$, respectively (this is hidden in the identification of $(\mathcal U\ominus \mathcal V)\ominus \mathcal W$ with $\mathcal U\ominus \mathcal V\ominus \mathcal W$ as in \eqref{eq:horizontal-unbracketed}). For the first splitting position of $s$, $\varphi_{\{k\stdl j\}}^{-1}(\varphi_{\{k\stdl j\}}(i))=i$ holds trivially. It is not hard to see that $\varphi_{[s_{-,0}]_0}\circ\varphi_{\{k\stdl j\}}=\varphi_{[s]_0\cap\{k\stdl j\}}=\varphi_{\{k\stdl j'\}}\circ\varphi_{[s]_0}$, so for the first splitting position of $r$ we have $\varphi_{\{k\stdl j'\}}^{-1}(\varphi_{[s_{-,0}]_0}(\varphi_{\{k\stdl j\}}(i)))=\varphi_{[s]_0}(i)$, as needed. This shows that (under the mentioned identification) we have $(R\overt \mathrm{id})\circ R=R_3$. We can argue analogously to show $(\mathrm{id}\overt R)\circ R=R_3$, which concludes the proof of the first associativity constraint Figure~\ref{fig:assI}.

 Pick now two composable triples
 \[
 t^{\stdl i}=q_- \circ r_- \circ s_- \qquad \text{and}\qquad t^{\stdg i}=q_+ \circ r_+ \circ s_+.
 \]
 Set also
 \begin{align*}
 (qr)_{-} := q_- \circ r_-,\qquad (qr)_{+} = q_+ \circ r_+,\qquad
 	(rs)_{-} = r_{-}\circ s_{-},\qquad (rs)_{+} = r_+ \circ s_+,
 \end{align*}
 so that
 \[
 t^{\stdl i} = (qr)_{-} \circ s_-,\qquad t^{\stdg i} = (qr)_{+} \circ s_+,\qquad
 t^{\stdl i} = q_-\circ (rs)_-, \qquad t^{\stdg i} = q_+ \circ (rs)_+.
 \]
 When we iterate Construction~\ref{constrctone} to find $q$, $r$ and $s$, the only ambiguity lies in the positions where the words are split (or merged, coming from their given splittings). Recall that for an ordered set $A$, we denote $\varphi_A\colon A\to \interval{|A|}$ the unique monotone bijection. If we start by merging $(qr)_{\pm}$ and $s_{\pm}$ to obtain $qr$ and $s$ with $qr\circ s=t$, this happens at $i'=\varphi_{[s]_0}(i)$ and $i$, respectively. If we then merge $q_{\pm}$ and $r_{\pm}$ by applying Construction~\ref{constrctone} with respect to $qr$, this happens at $\varphi_{[r]_0}(\varphi_{[s]_0}(i))$ and $\varphi_{[s]_0}(i)$, respectively. The following tabular summarizes the result:
 \begin{center}
 \begin{tabular}{|c|c|c|}\hline
 full word& factors & split position\\\hline\hline
 \multirow{2}{*}{$t$}&$qr$&$\varphi_{[s]_0}(i)$\\\cline{2-3}
 &$s$&$i$\\\hline\hline
 \multirow{2}{*}{$qr$}&$q$&$\varphi_{[r]_0}(\varphi_{[s]_0}(i))$\\\cline{2-3}
 &$r$&$\varphi_{[s]_0}(i)$\\\hline

 \end{tabular}
 \end{center}
 If we argue analogously for the other iteration, we get:
 \begin{center}
 \begin{tabular}{|c|c|c|}\hline
 full word& factors & split position\\\hline\hline
 \multirow{2}{*}{$t$}&$q$&$\varphi_{[rs]_0}(i)$\\\cline{2-3}
 &$rs$&$i$\\\hline\hline
 \multirow{2}{*}{$rs$}&$r$&$\varphi_{[s]_0}(i)$\\\cline{2-3}
 &$s$&$i$\\\hline

 \end{tabular}
 \end{center}
 Therefore, the associativity constraint in Figure~\ref{fig:assII} is equivalent to $\varphi_{[rs]_0}(i)=\varphi_{[r]_0}(\varphi_{[s]_0}(i))$. This follows easily from Proposition~\ref{prop:compose-iota} because, as a map between ordered sets, $\iota_s$ is nothing but $\varphi_{[s]_0}^{-1}$.
\end{proof}

 Existence of such a four-points exchange relation has major consequences for the categories of horizontal semigroups and vertical comonoids, as stated in the next proposition. For the proof, we refer the reader to \cite[Chapter 4]{aguiar2010monoidal}.

\begin{Proposition}	\label{prop:duoid}
	Pick $\big(\mathcal{V},\hp_\mathcal{V}\big)$ and $\big(\mathcal{W},\hp_\mathcal{W}\big)$ two $\ominus$-semigroups, then $\mathcal{V}\overt \mathcal{W}$ is a $\ominus$-semigroup if equipped with the product $m_{\mathcal{V}\overt \mathcal{W}}^{\ominus}$ defined by
	\begin{equation}\label{eqn:exchange-m^ominus}
		\hp_{\mathcal{V}\overt \mathcal{W}} := \big(\hp_{\mathcal{V}} \overt \hp_{\mathcal{W}}\big) \circ R_{\mathcal{V},\mathcal{W},\mathcal{V},\mathcal{W}}.
	\end{equation}

	Pick two $\SWc$-comonoids $\big(\mathcal{V},\vcp_{\mathcal{V}}\big)$ and $\big(\mathcal{W},\vcp_\mathcal{W}\big)$ then $\mathcal{V}\ominus \mathcal{W}$ is a $\SWc$-comonoid if equipped with the vertical coproduct
	\begin{equation*}
		\vcp_{\mathcal{V}\ominus \mathcal{W}} = R_{\mathcal{V},\mathcal{W},\mathcal{V},\mathcal{W}} \circ \big(\vcp_\mathcal{V} \ominus \vcp_\mathcal{W}\big).
	\end{equation*}
\end{Proposition}
This proposition implies that the category $\sgominus$,
is monoidal for the vertical tensor product, in fact $\uvert$ is a unit, since it can be equipped with a horizontal product
 \[
 \hp_{\uvert}\colon \ \uvert \ominus \uvert \to \uvert
 \]
 defined by, for any $w^{-}\ominus w^{+} \in \uvert \ominus \uvert (\alpha, 1_{\alpha}) $, $w=z\cdot \Singleton$, $w^{\prime} = z^{\prime}\cdot \Singleton$
 \[
 \hp_{\uvert}(w^{-}\ominus w^{+}) = zz^{\prime} \cdot \Singleton.
 \]

Likewise, the category of $\overt$-comonoids is a semigroup category if equipped with the horizontal semigroupal product $\ominus$.

As mentioned previously, all linear $\SWc$-indexed collections considered in this work come equipped with preferred \emph{finite bases}; they are all of the form $\linearspan{\mathcal{X}}$ with $\mathcal{X}$ an $\SWc$-indexed collections of finite sets. Additionally, they are (vertically) {\it augmented}, i.e., they come with a~morphism of linear $\SWc$-indexed collections
\[
	\eta\colon \ \uvert \to \linearspan{\mathcal{X}},
\]
called the \emph{augmentation morphism}, which induces an isomorphism for each translucent word $t$ with $b_t=\mathbf{0}_{|t|}$; more concretely one has
\begin{equation*}
	\linearspan{\mathcal X}{(\alpha,\mathbf{0})} \simeq_{\eta_{(\alpha,\mathbf{0})}^{-1}} \linearspan{\{\Singleton\}}
\end{equation*}
for all $\alpha\in\{\LL,\RR\}^\star$.

In the sequel, we use the notation ${\mathcal{X}}^+$ for the $\SWc$-indexed collection with
\begin{equation}
 \label{eq:Xplus}
	{\mathcal X^+}(t):=\begin{cases}
 {\mathcal X}(t) &\text{if $b_t\neq \mathbf{0}$},\\
	 \varnothing &\text{if $b_t= \mathbf{0}$}.
	\end{cases}
\end{equation}
The resulting linear collection is then given by
\begin{equation*}
	\linearspan{\mathcal X^+}(t):=\begin{cases}
 \linearspan{\mathcal X}(t) &\text{if $b_t\neq \mathbf{0}$},\\
	 \{0\} &\text{if $b_t= \mathbf{0}$}.
	\end{cases}
\end{equation*}
Furthermore, when considering $\SWc$-comonoid structures on collections of the form $\linearspan{\mathcal X}$, we will always assume \emph{nilpotency} for the coproduct $\vcp$, which means that for every $c \in \mathcal{X}(t)$, $t \in \SWc$, $b_t\neq {\bf 0}$,
\begin{align*}
	 & \vcp_t(c) = c\otimes_\CC \eta_{(\alpha_t,\mathbf{0})}(\Singleton) + \eta_{\restriction{t}{[t]_0}}(\Singleton) \otimes_{\CC} c + \bvcp(c), \\
	 & \vcp_{(\alpha,{\bf 0})}(\eta_{(\alpha,{\bf 0})}(\Singleton)) = \eta_{(\alpha,{\bf 0})}(\Singleton) \otimes \eta_{(\alpha,{\bf 0})}(\Singleton),
\end{align*}
with ${\bvcp}$ the reduced coproduct
\begin{equation*}
	\bvcp\colon \ \linearspan{\mathcal{X}}^+ \to \linearspan{\mathcal{X}}^+ \otimes \linearspan{\mathcal{X}}^+
\end{equation*}
and $\bvcp{}^n(c) = 0$ for a certain integer $n$ depending on $c$ where we define recursively $\bvcp{}^n = ({\bvcp} \otimes \mathrm{id}) \circ \bvcp{}^{n-1}$, ${\bvcp}{}^1=\bvcp$.

\subsection{Freely generated horizontal semigroups}
\label{subsec:freely-generated-ominus-semigroups}

\begin{Definition}[opaque collection]
 For each linear $\SWc$-indexed collection $\mathcal V$ we define its {\it opaque part} $\mathcal V_{\mathbf 1}$ as
 \[\mathcal V_{\mathbf 1}(t)=
 \begin{cases}
 \mathcal V(t)&\text{if $b_t=\mathbf 1$,}\\
 \{0\} &\text{else.}
 \end{cases}
 \]
	We call a linear $\SWc$-indexed collection $\mathcal{V}$ {\it opaque} if
	$\mathcal V=\mathcal V_{\mathbf 1}$.
	The full sub-category of $\icsw$ of all opaque $\SWc$-indexed collections is denoted $\icsw_{\mathbf{1}}$.

 Analogously, a $\SWc$-indexed collections of sets $\mathcal X$ is opaque if it coincides with $\mathcal X_{\mathbf 1}$ defined as
 \[\mathcal X_{\mathbf 1}(t)=
 \begin{cases}
 \mathcal X(t)&\text{if $b_t=\mathbf 1$,}\\
 \varnothing &\text{else},
 \end{cases}
 \]
 and $\Set^\SWc_{\mathbf 1}$ is the full sub-category of $\Set^\SWc$ consisting of all opaque collections.
\end{Definition}
\begin{Remark}
 Opaque collections $\mathcal{V}$ correspond $1:1$ to collections of vector spaces indexed by words in $\{\LL,\RR\}^{\star}$.
\end{Remark}
We consider a functor $U\colon\sgominus\to \icsw_{\mathbf{1}}$, from the category of $\ominus$-semigroups to the category of all linear opaque $\SWc$-indexed collections defined by
\[
U(\mathcal{V}, \hp_{\mathcal{V}})(t)= \begin{cases} \mathcal{V} (t), & t=(\alpha,{\bf 1}), \\ \{0\}, & \text{otherwise}. \end{cases}
\]
\begin{Remark}
 \label{rmq:graded-subspace}
 The functor $U$ is not a forgetful functor in the sense that it simply ``forgets'' about the horizontal product $\hp_{\mathcal{V}}$: it does also set certain components of $\mathcal{V}$ to $\{0\}$.

 Note that $U(\mathcal V)(t)$ is a subspace of $\mathcal V(t)$ for all $t\in \SWc$; either the two are the same or $U(\mathcal V)(t)$ is trivial. This means that $U(\mathcal V)$, seen as a graded vector space, is a graded subspace of $\mathcal V$. This will play an important role in the next proof.
\end{Remark}
A left-adjoint for $U$ is a functor $\mathcal{F}\colon \icsw_{\mathbf{1}} \to \sgominus$ satisfying
\begin{gather*}
	\mathrm{Hom}_{\icsw_{\mathbf{1}}}(\mathcal{V},U(\mathcal{W})) \simeq \mathrm{Hom}_{\sgominus}(\mathcal{F}(\mathcal{V}),\mathcal{W}),\\ \mathcal{V}\in \icsw_{\mathbf{1}},\qquad \mathcal{W} \in \sgominus,
\end{gather*}
where we have used the notation
\[
	\mathrm{Hom}_{\sgominus}(\mathcal{V},\mathcal{W}) := \{ f\in \mathrm{Hom}_{\icsw}(\mathcal{V},\mathcal{W})\colon \hp_\mathcal{W} \circ (f \ominus f) = f\circ\hp_{\mathcal{V}}\}.
\]
\begin{Proposition}\label{prop:freefunctor}
Define the functor $\mathcal{F}\colon \icsw_{\mathbf{1}} \to \sgominus$, for any opaque collection $\mathcal{V} \in \icsw_{\mathbf{1}} $, as
\begin{equation*}
	\mathcal{F}(\mathcal{V}) = \bigoplus_{n\geq 1} {\mathcal{V}}^{\ominus n} = \mathcal{V} \oplus \mathcal{V}^{\ominus 2} \oplus \mathcal{V}^{\ominus 3} \oplus \cdots,\qquad \mathcal{F}(f) = \bigoplus_{p\geq 1} f^{\ominus p}
\end{equation*}
equipped with its canonical horizontal semigroupal product $\hp_{\mathcal{F}(\mathcal{V})}$. Then $\mathcal{F}$ is a left-adjoint for~$U$.
\end{Proposition}
\begin{proof}
For the sake of the proof, we extend the functor $\mathcal{F}$ by the same formula to the incidence category $\icsw$.
 Let $\mathcal V\in \Vect_{\mathbf{1}}^\SWc$ be an opaque collection and $\mathcal W\in \sgominus$ a $\ominus$-semigroup.
 Note that the product $\hp_{\mathcal{W}}$ yields (thanks to associativity) a morphism of $\SWc$-collections
 \[
 \mathsf{m}_{\mathcal{W}}^{\ominus p}\colon \ \mathcal{W}^{\ominus p} \to \mathcal{W}.
 \]
 Altogether, the morphisms $\mathsf{m}_{\mathcal{W}}^{\ominus p}$, $p \geq 1$ yield a morphism of $\SWc$-collections
 \[\underline{\hp_{\mathcal{W}}}\colon \ \mathcal{F}(\mathcal{W}) \to \mathcal{W},\]
 which restricts to a morphism of $\SWc$-collections, also denoted $\underline{\hp_{\mathcal{W}}}$, from $\mathcal{F}(U(\mathcal{W}))$ to $\mathcal{W}$ (we interpret $\mathcal{F}(U(\mathcal{W}))$ as a graded subspace of $\mathcal F(W)$ as in Remark~\ref{rmq:graded-subspace}).
 Pick a morphism of opaque collections $f \colon \mathcal{V}\to U(\mathcal{W})$ and define
 \[f^{\ominus}\colon \ \mathcal{F}(\mathcal{V})\to \mathcal{W} ,\qquad f^{\ominus} := \underline{\hp_{\mathcal{W}}}\circ \mathcal{F}(f).\]
 Then $f^{\ominus}$ is a morphism of $\ominus$-semigroups. This stems from the fact that
 \[
 \mathcal{F}(f) \circ \hp_{\mathcal{F}(\mathcal{V})} = \hp_{\mathcal{F}(\mathcal{W})} \circ (\mathcal{F}(f) \ominus \mathcal{F}(f))
 \]
 and (from associativity of $\hp_{\mathcal{W}}$)
 \[
 \underline{\hp_{\mathcal{W}}} \circ \hp_{\mathcal{F}(W)} = \hp_{\mathcal{W}} \circ \underline{\hp_{\mathcal{W}}} \ominus \underline{\hp_{\mathcal{W}}}.
 \]
 Pick a morphism of horizontal semigroups $g\colon \mathcal{F}(\mathcal{V}) \to \mathcal{W}$, then $g$ is a morphism of $\SWc$-collections, hence
 \[
 g(\mathcal{V}(\alpha,\mathbf{1})) \subset \mathcal{W}(\alpha,\mathbf{1})
 \]
 (when considering $\mathcal{V} \subset \mathcal{F}(\mathcal{V})$).
 This implies that $g$ yields a morphism between the opaque collections $\mathcal{V}$ and $U(\mathcal{W})$). This ends the proof.
 \end{proof}

 Of course, the analogous statement for $\Set_{\mathbf 1}^{\SWc}$ with
 \[\mathcal F(\mathcal X):=\bigsqcup_{p\geq 1} \mathcal X^{\ominus p},\qquad U(\mathcal{X}, \hp_{\mathcal{X}})(t)= \begin{cases} \mathcal{X} (t), & t=(\alpha,{\bf 1}), \\ \varnothing, & \text{otherwise}, \end{cases}
 \]
 holds likewise.

\begin{Remark}\label{rk:freelygenerated}
 In the following, we say that a horizontal semigroup $\mathcal{V}$ is {\it freely generated by its opaque part} if it is isomorphic as a horizontal semigroup to $\mathcal{F}(U(\mathcal{V}))$. This means that to every element $v \in \mathcal{V}(t)$ corresponds a unique sequence $v_1,\dots,v_p$ where $p=|[t]_0|$ and $v_i \in U(\mathcal{V})$, $i=1,\dots,p$, such that
 \[
 v = \hp_{\mathcal{V}}(v_1\ominus \cdots \ominus v_p).
 \]
\end{Remark}

\subsection{Monoids of bipartitions and noncrossing bipartitions}\label{subsec:incompletebinoncrossing}

We will exploit the results proved in this section in the proofs of our main theorems in Section~\ref{sec:momentscumulantsrelations}.
Let $X$ be a set. An {\it incomplete partition of $X$} is a partition of a subset of $[X]_1\subset X$. An incomplete partition $\pi$ of $X$ becomes a partition $\widetilde \pi$ of $X$ by adding the block $[X]_0:=X\setminus [X]_1$, $\widetilde \pi=\pi\cup \{[X]_0\}$. We refer to the blocks of $\pi$ as opaque blocks and to the block $[X]_0$ of $\widetilde \pi$ as the translucent block of $\pi$ (although it is not formally an element of $\pi$ as a set). An incomplete partition defines a translucent-opaque structure on $X$, points contained in opaque blocks are opaque and points contained in the translucent block are translucent. If $X$ is a translucent set from the beginning, we modify the concept of partition so that it fits with the above observations.

\begin{Definition}
 A {\it partition of a translucent set} $X$ is a partition $\pi$ of the set of opaque points~$[X]_1$. We refer to the set of translucent points as the {\it translucent block} of $\pi$ although it is not formally an element of $\pi$ regarded as the set of its blocks.
\end{Definition}

\begin{figure}[!ht]\centering
	\includegraphics[scale=0.48]{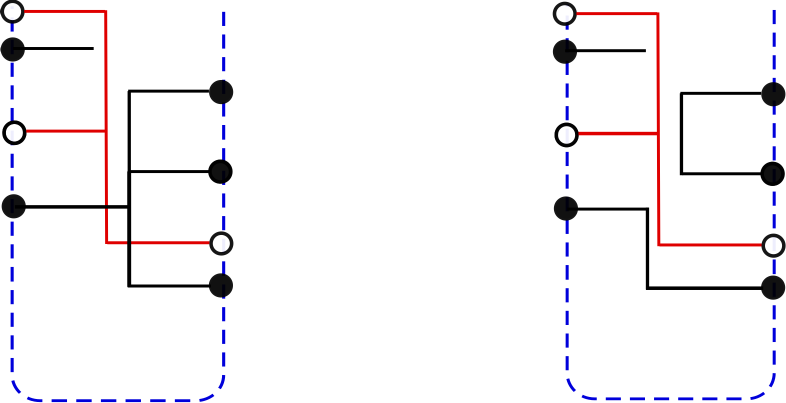}
	\caption{On the left an incomplete bipartition and on the right an incomplete noncrossing bipartition.
		The underlying translucent word is $({\sf LLRLRLRR},01101101)$ in both examples.}	\label{fig:incompletebipartition}
\end{figure}

\begin{Remark}
 In diagrams, the translucent block is coloured red, see for example Figure~\ref{fig:incompletebipartition}.
\end{Remark}

A partition of a translucent set is an incomplete partition of the underlying set which is compatible with the translucent-opaque structure.
The same applies to finite translucent ordered bisets and therefore to translucent words. Because of the importance of the concept, we spell it out concretely.

\begin{Definition} A {\it bipartition of a translucent word} $t$ is a partition of $[t]_1$, i.e., a bipartition of the associated finite translucent ordered biset $X_t=(\interval{|t|},\alpha_t,b_t)$. The set of all bipartitions of a translucent word $t$ is denoted
$\BPc(t)$. The collection
$\bpinc = (\bpinc(t))_{t \in \SWc}$ is called the
{\it $\SWc$-indexed collection of incomplete bipartitions}.
\end{Definition}

\begin{Remark}In the remainder of this article, if we speak of an {\it incomplete bipartition}, we always mean a bipartition of some translucent word in the sense of the previous definition.
\end{Remark}

We will need to extend the notion of restriction to partitions. Let $\pi$ be a~bipartition of a~translucent word $t$ and $I=\{i_1<\dots<i_p\}$ a subset of $\interval{|t|}$. Then $\restriction{\pi}{I}$ is the unique bipartition of~$\restriction{t}{I}$ such that $i,j\in\interval{p}$ are in the same block of $\restriction{\pi}{I}$ if and only if $i_k,i_\ell$ are in the same block of~$\pi$. It is straightforward to check that the translucent block of $\restriction{\pi}{I}$ is $\{k\colon i_k\in [t]_0\}=[\restriction{t}{I}]_0$, so that $\restriction{\pi}{I}$ is indeed a partition of the translucent word~$\restriction{t}{I}$.

With these notions at hand, we are now ready to define the composition of incomplete bipartitions.

\begin{figure}[!htb]\centering
	\includegraphics[scale=0.45]{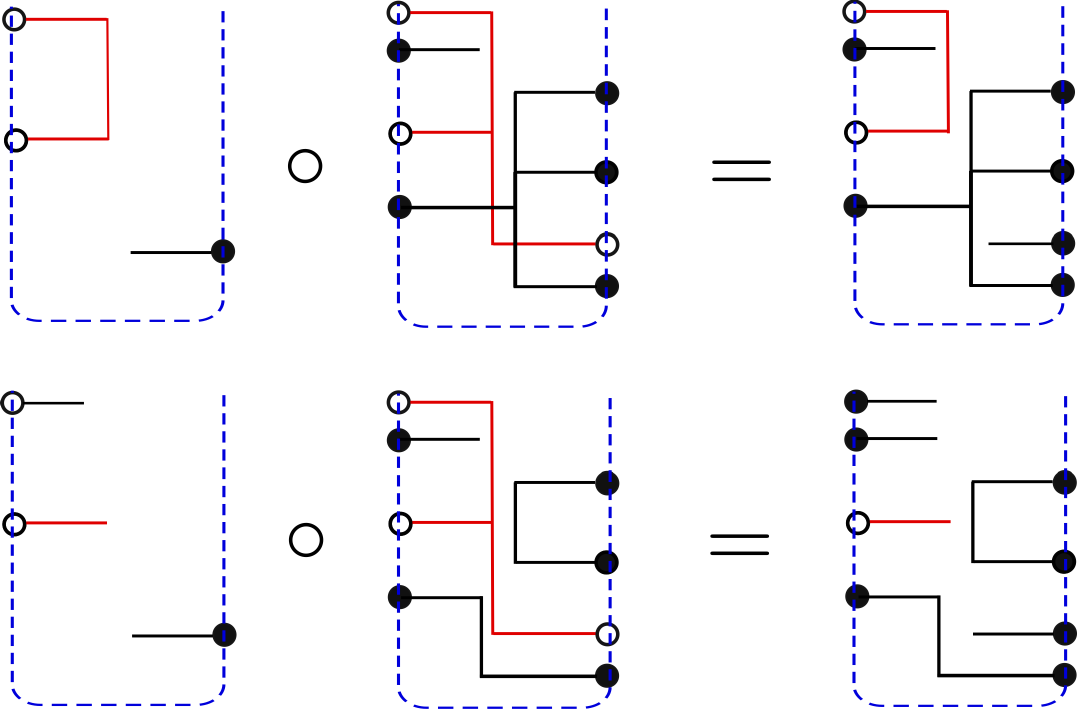}
 \caption{Two examples of compositions between bipartitions and noncrossing bipartitions.} \label{fig:incompletebipartitioncomposition}
\end{figure}

\begin{Definition}\label{def:incompletebipartitioncomposition}
 Let $r$ and $s$ be composable translucent words, $\rho\in\bpinc(r)$, $\sigma\in \bpinc(s)$. The composition ${\rho\circ\sigma}$ is the unique bipartition $\pi$ of $t:=r\circ s$ such with $\restriction{\pi}{[s]_0}=\rho$ and $\restriction{\pi}{[s]_1}=\restriction{\sigma}{[s]_1}$. Roughly speaking, the translucent block of $\sigma$ is substituted by the bipartition $\rho$. See Figure~\ref{fig:incompletebipartitioncomposition} for an example.
\end{Definition}

Henceforth, we adopt the same notation $\overt$ for the vertical monoidal product $\circ$ on $\SWc$-collections of sets and on $\SWc$-collections of vector spaces.

\begin{Observation}The collection $\bpinc$ is a $\SWc$-monoid with the product
\[\vpBP\colon \ \bpinc \overt \bpinc \to \bpinc,\qquad \vpBP(\rho,\sigma):=\rho\circ\sigma. \]
\end{Observation}
The unit for the vertical product $\vpBP$ is the morphism
\[
\eta_{\BPc} \colon \ \uvert \to \bpinc
\]
sending $\Singleton \in \BPc(\alpha,\mathbf{0})$ to the unique incomplete bipartition in $\BPc(\alpha,\mathbf{0})$, the empty partition.

For applications in noncommutative probability, noncrossing partitions are most important.

\begin{Definition}
 An incomplete bipartition $\pi\in \bpinc(t)$ is called {\it noncrossing} if the extended bipartition $\pi\cup\{[t]_0\}$ is a noncrossing bipartition of $\alpha_t\in\{\LL,\RR\}^\star$. The set of all noncrossing bipartitions of $t$ is denoted $\bncinc{(t)}$. The collection of all incomplete noncrossing bipartitions is denoted $\bncinc:=(\bncinc(t))_{(t\in\SWc)}$.
\end{Definition}

It is obvious that the composition of noncrossing incomplete bipartitions is again an incomplete noncrossing bipartition. Therefore, we document for future reference:

\begin{Observation}
 The collection $\bncinc$ $($and consequently its linearization $\mathbb C\bncinc)$ is a $\SWc$-monoid with the product
 \[\vpBNC\colon \ \bncinc \overt \bncinc \to \bncinc,\qquad \vpBNC(\rho\overt\sigma):=\rho\circ\sigma \]
 $($or the linearization thereof$)$.
\end{Observation}

\begin{figure}[!htb]
 \centering
 \includegraphics[scale=0.5]{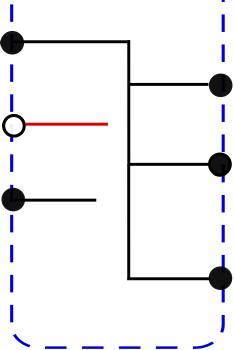}
 \caption{Example of an incomplete bipartition. The collection $\BNCc$ is not in the image of the functor~$\mathcal{F}$.} \label{fig:contreexemple}
\end{figure}

Neither the collection $\mathbb C\BPc$ nor $\mathbb C\BNCc$ lies in the image of the functor $\mathcal{F}$ of Proposition \ref{prop:freefunctor}, i.e., those collections are not freely generated horizontal semigroups. Indeed, $A=\mathcal F(B)$ implies $A_{\mathbf 1}= \mathcal F(B)_{\mathbf 1}=B_{\mathbf 1}$. Now consider the noncrossing bipartition \[\pi\in \BNCc(\LL\RR\LL\RR\RR\RR,110111)\subset \BPc(\LL\RR\LL\RR\RR\RR,110111)\] from Figure~\ref{fig:contreexemple}. Then
\[
\pi\notin\mathcal F(\BPc_{\mathbf 1})(\LL\RR\LL\RR\LL\RR,110111)=\BPc_{\mathbf 1}^{\ominus 2}(\LL\RR\LL\RR\LL\RR,110111)\supset \BNCc_{\mathbf 1}^{\ominus 2}(\LL\RR\LL\RR\LL\RR,110111),
\] where by definition the blocks have to be contained in the $\stdl$-intervals $\rrbracket{-\infty_\stdl},3\llbracket_\stdl$ and $\rrbracket 3,\infty_\stdl\llbracket_\stdl$.
Below, we define a subcollection of incomplete bipartitions (following~\cite{CNS15:random-variables} we call them \emph{shaded bipartitions}) freely generated by the opaque collection of noncrossing bipartitions.

Recall that given a partition $\pi$ of $\interval{n}$, $n\geq 1$, we write $\sim_{\pi}$ for the equivalence relation induced by $\pi$ on $\interval{n}$.

\def\u{{\color{red}\scalebox{0.9}{$\pmb{|}$}}}

\begin{Definition}[shaded noncrossing bipartition]	\label{def:spadeincompletenonbipartition}
	Pick a translucent word $t\in \SWc$ and an incomplete noncrossing bipartition $\pi \in \BNCc(t)$ (recall that $\pi$ is a partition of $[t]_1$). We say that $\pi$ is a {\it \sincp} of type $t$ if
		 for all integers $1\leq k< \min [t]_0$ with $\alpha(k)=\alpha(\min [t]_0)$ we have
		 \[ j \sim_{\pi} k\implies j < \min [t]_0 \qquad \text{and} \qquad \alpha(j)=\alpha(\min[t]_0).\]
\end{Definition}

\begin{Remark}An alternative description of a partition $\pi \in \shBNCc(t)$ is as follows. The spines-and-ribs diagram of $\pi \cup [t]_0$ to which a \emph{straight vertical chord emerging from the top} of the diagram and connected to the translucent block $[t]_0$ is added has no crossings. In other words, the translucent block is an outer block of the associated noncrossing partition of $\interval{|t|}$ with respect to necklace order. As such, \sincps~are \emph{shaded diagrams} in the sense of~\cite{CNS15:random-variables} with one top chord. In the following figures, we choose to colour the translucent block red.\looseness=-1
\end{Remark}

\begin{figure}[!ht]\centering
	\includegraphics[scale=0.7]{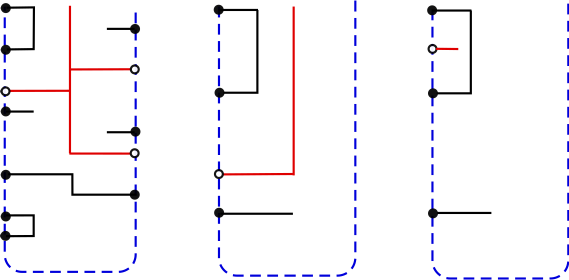}
	\caption{Two partitions on the left are \sincps, the one on the right is not \emph{not} \sincp.}\label{fig:incompletebinoncrossingpartition}
 \end{figure}

See Figure~\ref{fig:incompletebinoncrossingpartition} for examples and counterexamples. Of course, any {\sincp} is an incomplete noncrossing bipartition.

We denote by $\shBNCc \subset \bncinc$ the $\SWc$-collection comprising all {\sincp}s.

Before we observe that $\shBNCc$ is indeed freely generated by noncrossing bipartitions, we introduce some notation.

\begin{Notation}\label{notation:INT}
For a translucent word $t\in\SWc$,
we let $\INT{t}$ denote the set of maximal $\stdl$-intervals inside $[t]_1 \subset \llbracket |t| \rrbracket$,
\begin{equation*}
	\INT{t} = \big\{I_1^{\stdl}(t) \stdl \cdots \stdl I^{\stdl}_{k(t)}(t)\big\}
\end{equation*}
and call $I_j^{\stdl}(t)$ the $j$th $\stdl$-interval.
\end{Notation}

\begin{Example}
	Consider $t = ({\sf LRLLLR}, 011101)$, then
	\begin{align*}
		I^{\stdl}_1(t) = \{3,4\},\qquad I^{\stdl}_2(t) = \{2,6\},
	\end{align*}
\end{Example}
\begin{Remark}
	This sequence of intervals $I_j^{\stdl}(t)$, $1 \leq j\leq k(t)$ has been introduced in the proof of Proposition~\ref{prop:freealgebra} in a slightly different way to include empty intervals separating consecutive integers in $[t]_0$.
\end{Remark}

\begin{Observation}\label{rk:fondrk}
	A fundamental observation is the following:
	\sincps~can be reconstructed from their restrictions to the maximal $\stdl$-intervals.
	More precisely, given a translucent word $t$, 
	each block of a \sincp~$\pi$ of type $t$ is contained in one of the intervals $I^{\stdl}_j(t)$, $1 \leq j\leq k(t)$ which make up $\INT{t}$. Thus,{\samepage
	\begin{equation*}
		\shBNCc(t)\ni \pi \mapsto \restriction{\pi}{I_1^{\stdl}}, \dots, \restriction{\pi}{I_{k(t)}^{\stdl}} \in \shBNCc(\restriction{t}{I_1^{\stdl}})\times\cdots\times \shBNCc(\restriction{t}{I_{k(t)}^{\stdl}})
	\end{equation*}
	is well-defined and a bijection.}

 As a consequence, the collection of shaded noncrossing bipartitions is identified as the free $\ominus$-semigroup over the collection of noncrossing partitions, i.e.,
 \[\shBNCc=\mathcal F(\BNCc_{\mathbf 1}).\]
 Indeed,
 \[\mathcal F(\BNCc_{\mathbf 1})(t)=\BNCc_{\mathbf 1}^{\ominus k(t)+1}(t)=\BNCc_{\mathbf 1}(\restriction{t}{I_1^{\stdl}})\times\cdots\times \BNCc_{\mathbf 1}(\restriction{t}{I_{k(t)}^{\stdl}})\cong \shBNCc(t)\]
 because the opaque parts of $\shBNCc$ and $\BNCc$ agree.
\end{Observation}

\begin{Remark} We compare the restriction of the product $\mathsf m_{\mathcal{NC}}^{\overt}$ to (one-sided) noncrossing partitions and the gap-insertion operad introduced in~\cite{ebrahimi2020operads}.
	First, the monoidal product $\mathsf m_{\mathcal{NC}}^{\overt}$ restricts to the sub-$\SWc$-collection of $\shBNCc$ comprising all shaded noncrossing partitions; the {\sincps} with type $({\LL}^n,b)$, $n\geq 1$ and $b$ a Boolean word with length $n$.

	Pick a shaded noncrossing bipartition $\pi$ with type $t=({\LL}^n,b)$ and denote by $I_0,\dots,I_p$ the connected components (i.e., maximal $\stdl$-intervals) of $[t]_1$. The shaded noncrossing partition~$\pi$ yields noncrossing partitions $\pi_0,\dots,\pi_p$, obtained as the restrictions of~$\pi$ to the intervals $I_0,\dots,I_p$.
	Pick a noncrossing partition $\pi^{\prime}$ with type $({\LL}^p, \mathbf{1})$. Denote by $\rho$ the gap-insertion operadic product introduced in~\cite{ebrahimi2020operads}. Then
	\begin{equation*}
		\rho(\pi^{\prime},\pi_0,\dots,\pi_p) = \mathsf m_{\mathcal{NC}}^{\overt}(\pi^{\prime} \overt \pi).
	\end{equation*}
\end{Remark}

Finally, we can of course linearize our monoids to obtain monoids in $\Vect^\SWc$.

\begin{Observation}The linear extensions of $\vpBP$ and $\vpBNC$ turn the linear collections $\linearspan{\bpinc}$ and $\linearspan{\bncinc}$ into linear $\SWc$-monoids.
\end{Observation}

Later we will also need incomplete versions of interval and monotone partitions, which we define now (which do not give rise to $\overt$-monoids).

\begin{Definition}[incomplete interval bipartitions]
	An {\it interval bipartition} of type $t$ is an incomplete bipartition $\pi\in \BPc(t)$ such that each block of $\pi$ is an interval for the restriction of the order $\stdl$ to $[t]_1$; note that $[t]_0$ is not considered a block of $\pi$ and therefore is not necessarily a~$\stdl$-interval.
\end{Definition}
Given a translucent word $t$, we denote by $\BBc(t)$ the set of incomplete interval bipartitions.
\begin{figure}[!ht]\centering
	\includegraphics[scale=0.7]{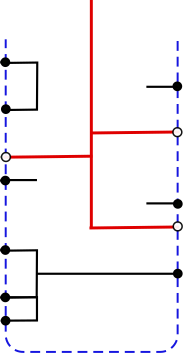}
	\caption{Example of an incomplete interval bipartition.}\label{fig:incompletebi-Boolean}
\end{figure}

\begin{Definition}[incomplete monotone bipartitions]
An {\it monotone bipartition of type $t$} is a~noncrossing bipartition $\pi \in \BNCc(t)$ equipped with a total order (represented by an injective labeling $\lambda\colon \pi \to \interval{|\pi|}$) such that, when the order is extended to $\{[t]_{0}\}\cup \pi$ by declaring $[t]_0$ the minimum block (represented by $\lambda([t]_0)=0$), the resulting partition is a monotone bipartition of~$\alpha_t$ (i.e., a monotone partition of $\interval{|t|}$ with respect to necklace order). In the sequel, we denote monotone partitions as pairs $(\pi,\lambda)$ to make the labelling explicit. We use the term \emph{incomplete monotone partition} to refer to any monotone partition of type $t$ for some translucent word $t$.

Given a translucent word $t \in \SWc$, we denote by $\MBNCc(t)$ the set of all monotone bipartitions of type $t$.
\end{Definition}

\begin{Remark}One could define the set $\shMBNCc$ of shaded monotone bipartitions of type $t$ similar to the noncrossing case, but it amounts to the same thing; indeed, since the block $[t]_0$ is the minimal block, it must be an outer block and can always be extended to the top. In diagrams of incomplete monotone partitions, we still draw the red extension in order to make the comparison with the noncrossing case easier.

Similarly, an incomplete interval bipartition is automatically shaded in the obvious sense.
\end{Remark}

\section{Words on a two-faced family of random variables}\label{sec:words}

Recall that the letters $\LL$ and $\RR$ shall act as placeholders, meant to be substituted by left or right random variables.
In this section, we formalize this through the introduction of a vertical monoid in the category $\icsw$ (depending on a choice of a finite set of left and right variables). This monoid should be thought of as a two-faced analogon of
\emph{double bar construction} of an algebra from which the shuffle point-of-view on moment-cumulant relations is derived \cite{ebrahimi2015cumulants}.
It is supported by a $\SWc$-collection of \emph{incomplete words} on random variables. Some of the entries of these {incomplete words} are drawn from the set of random variables we fixed and the others are drawn from the placeholders $\LL,\RR$.

Finally, we define \emph{concatenation} of two incomplete words as a horizontal product, yielding the structure of a $\ominus$-semigroup. Besides, this semigroup is freely generated by words on random variables (the opaque part of the collection of incomplete words, see Remark \ref{rk:freelygenerated}.)

With the aforementioned vertical composition, this second product on incomplete words yields (almost) a \emph{dimonoid} on incomplete words (to be precise, it is not quite a dimonoid because one of the operations only forms a semigroup, not a monoid).

\subsection{Collection of incomplete words on random variables}

Let $(\mathcal A,\varphi)$ be an algebraic probability space and $\big(\mathcal A^\LL, \mathcal A^\RR\big)$ a pair of faces, i.e., a pair of subalgebras of $\mathcal A$.
The unital linear functional $\varphi$ will not play a role and can be forgotten about while reading this section, but of course, it will be important later on when we discuss moment-cumulant relations.

\newcommand\bulletAccent[1]{\accentset{\bullet}{#1}}

Pick two \emph{finite} subsets $\bulletAccent{S}^{\LL} \subset \mathcal{A}^{\LL}$, $\bulletAccent{S}^{\RR}\subset\mathcal{A}^{\RR}$ of random variables in this probability space and define
\[
	\bulletAccent{S} = \bulletAccent{S}^{\LL}\sqcup \bulletAccent{S}^{\RR}
\]
as their disjoint union. We introduce more notation.
\begin{enumerate}\itemsep=0pt
	\item
	 Put
	 $
	 S^{\LL}:=\bulletAccent{S}^{\LL}\cup\{\sf L \}$, $S^{\RR}:=\bulletAccent{S}^{\RR}\cup\{\sf R \}
	 $
	 and
	 $S := S^{\LL} \sqcup S^{\RR} = \bulletAccent{S} \sqcup \{ \sf L, R\}$.

	\item As for any alphabet, we denote $S^\star$, $\bulletAccent{S}^\star$, etc.\ the sets of words with letters (or entries) from the set, including the empty word $\varnothing$.
	\item Words in $\bulletAccent S^\star$ are referred to as {\it words on random variables} or {\it complete words} for short. Words in $S^\star$ are referred to as {\it incomplete words on random variables} or {\it incomplete words}, for short.
	\item We put
	\begin{alignat*}{3}
	 & S^{(\LL,0)}:=\{\LL\},\qquad && S^{(\RR,0)}:=\{\RR\},&\\
	 & S^{(\LL,1)}:=\bulletAccent S^{\LL},\qquad &&
	 S^{(\RR,1)}:=\bulletAccent S^\RR&
	\end{alignat*} and, for a translucent word $t$ of length $n$,
 \begin{align*}
 \IWc(t):=S^{t(1)}\times S^{t(2)}\times\dots \times S^{t(n)}\subset S^\star;
 \end{align*}
 recall the convention $t(k)=(\alpha_t(k),b_t(k))$. The $\SWc$-indexed collection $\IWc=(\IWc(t))_{t\in \SWc}$ is called {\it collection of incomplete words on random variables}.
	\item We use the set decomposition \[S^\star=\bigsqcup_{t\in\SWc} \IWc(t)\]
 to define the {\it type} of a word in $\mathcal{W}$ as the translucent word $\Type(w)$, or shortly $t_w$, if $w\in \IWc(t_w)$.
 For example, if $a^\LL\in \bulletAccent S^{\LL}$ and $a^{\RR}\in \bulletAccent S^{\RR}$, then $\Type\big(a^\LL a^\RR \LL a^\RR\RR\big)=({\LL\RR\LL\RR\RR, 11010})$
	\item For an incomplete word $w$, we write $[w]_1:=[\Type(w)]_1$ for the set of \emph{opaque positions} and $[w]_0:=[\Type(w)]_0$ for the set of \emph{translucent} positions.
	\item For a word in $w\in\IWc(t)$, we define $\source{w}:=\source{t}=\alpha_t$ and $\target{w}:=\target{t}=\restriction{\alpha_t}{[t]_0}$.
\end{enumerate}

\begin{Example}
 Let $\bulletAccent{S}^{\LL} = \big\{a^\LL, b^\LL \big\}$, $\bulletAccent{S}^{\RR} = \big\{a^\RR\big\}$.
 Then
 \begin{align*}
 S = \big\{\LL,\RR,a^\LL,b^\LL,a^\RR\big\}.
 \end{align*}

 And then, for example,
 \begin{align*}
 \IWc( \LL\RR\LL,011 )
 =
 \big\{ \LL a^\RR a^\LL, \LL a^\RR b^\LL \big\}.
 \end{align*}
\end{Example}

\begin{Remark}
Again, we will mostly work with the linear $\SWc$-indexed collection $\linearspan{\IWc}$ instead.
\end{Remark}

\subsection{Monoid and comonoid of incomplete words}\label{ss:DmonoidDcomonoid}

Given two words $w\in\IWc{(s)},w'\in\IWc{(t)}$ with $\source s=\target t$, we define their \textit{composition} as the unique word $w\circ w'\in \IWc(s\circ t)$ with \[\restriction{(w\circ w')}{[t]_1}=\restriction{w'}{[t]_1} \qquad \text{and} \qquad \restriction{(w\circ w')}{[t]_0}=w.\] To put it another way, the word $w$ overwrites the placeholders in $w'$.

\begin{Example}
Continuing the example from Figures~\ref{fig:exampleComposition} and \ref{fig:exampleComposition-diagram}, let
	\begin{align*}
		w = x^\LL \LL x^\RR y^\RR \in \IWc(\LL\LL\RR\RR, 1011),\qquad
 w^{\prime}=\LL a^\LL a^\RR \LL\RR b^\LL \RR b^\LL
 \in \IWc(\LL\LL\RR\LL\RR\LL\RR\RR, 01100101).
	\end{align*}
	Then
	\begin{align*}
		w \circ w' = x^\LL a^\LL a^\RR \LL x^\RR b^\LL y^\RR b^\LL \in \IWc(\LL\LL\RR\LL\RR\LL\RR\RR, 11101111).
	\end{align*}
\end{Example}

\begin{Observation}
 It is not difficult to see that the composition of words defined above is associative. With
 \[\mathsf{m}^{\overt}\colon \ \IWc\overt \IWc\to \IWc, \qquad (w,w')\mapsto w\circ w'\]
 and the unit $\eta^{\overt}\colon \mathcal E^{\overt} \to \IWc$ defined by
	\begin{equation*}
		\eta^{\overt}_{(\alpha,\mathbf{0})} (\Singleton) = \alpha,
	\end{equation*}
	$\IWc$ becomes a $\overt$-monoid.
 By linearly extending $\mathsf{m}^{\overt}$ and $\eta^{\overt}$, we obtain a linear $\overt$-monoid structure on $\CC\IWc$ with
 \begin{align*}
 \vpS\colon \ \CC\IWc\overt \CC\IWc\to \CC\IWc,\qquad \vpS(w\overt w')=w\circ w'.
 \end{align*}
 $($Recall Notation~{\rm \ref{notation:overt}}, i.e., we write $w\overt w'$ instead of $w\otimes w'$ if the types of $w$ and $w'$, and therefore $w$ and $w'$ themselves, are composable.$)$
\end{Observation}

We consider now the coproduct $\vcpS$ dual to the composition product $\vpS$:
\begin{equation*}
	\vcpS\colon \ \linearspan{\IWc} \to \linearspan{\IWc}\overt \linearspan{\IWc},\qquad
	\vcpS(w) =
	\sum_{w=w'\circ w''}
	w' \overt w''.
\end{equation*}
On the right-hand side of the formula, we again used Notation~\ref{notation:overt}.

We will now give a more concrete formula for $\vcpS$ by using the operations of restriction and translucidation, extended to incomplete words.

Let $w$ be an incomplete word and an subset $I=\{i_1<\dots<i_p\}\subset \interval{|w|}$. The restriction of~$w$ to~$I$ is already defined for words over arbitrary alphabets, in this case, we get
\[\restriction{w}{I}:=w(i_1)\dots w(i_p).\]
Translucidation has already been defined for translucent words, see equation~\eqref{eq:tranlucidation}.
For the incomplete word $w$, we define the translucidation at $I$ as the word $\translucidation{w}{I}$ with letters
\[\translucidation{w}{I}(k)=\begin{cases}
 \LL &\text{if $k\in I$ and $w(k)\in S^\LL$},\\
 \RR &\text{if $k\in I$ and $w(k)\in S^\RR$},\\
 w(k)&\text{if $k\notin I$}.
\end{cases}\]

\begin{example*}
$ w = a^{\ell}b^{r}\LL a^{r} \RR$, $\Type(w)=\LL\RR\LL\RR\RR$, $I=\{2,5\}$, $\translucidation{w}{I}=a^{\ell}\RR\LL a^r\RR$, $\restriction{w}{I} = b^r\RR$.
\end{example*}

Note that $\Type(\restriction{w}{I})=\restriction{\Type(w)}{I}$ and $\Type(\translucidation{w}{I})=\translucidation{\Type(w)}{I}$. If $I\supset [w]_0$, then the translucent words $\Type(\restriction{w}{I})$ and $\Type(\translucidation{w}{I})$ are composable and
\begin{equation}
\label{eqn:factorization}
w=\restriction{w}{I}\circ\translucidation{w}{I}.
\end{equation}
Furthermore, any factorization of $w$ is of the form \eqref{eqn:factorization} for a unique subset $[w]_0\subset I\subset\interval{|w|}$. This yields the following concrete description of the coproduct $\vcpS$.

\begin{Observation} \label{obs:coproduct}
 For any incomplete word $w$,
 \[\vcpS(w)=\sum_{\substack{I \subset \interval{|w|} \\ I\supset{[w]_0}}}{\restriction{w}{I}}\overt{\translucidation{w}{I}}.\]

An {\it admissible cut} of a word $w$ is a pair $(\ell,u)$ of incomplete words, with
\[
\ell = \restriction{w}{I},\qquad u=\translucidation{w}{I}
\]
and $ [w]_{0} \subset I \subset \llbracket |w| \rrbracket$. We denote by ${\rm Adm}(w)$
the set of all admissible cuts of a word $w$. With these notations,
\[
\vcpS(w)=\sum_{(\ell,u) \in {\rm Adm}(w)} \ell \overt u.
\]
\end{Observation}

\begin{Example}
	\begin{gather*}
		\vcpS\big(a_1^{\LL}\LL\RR a_2^{\RR} a_3^{\RR}\big)
		 = \LL\RR \otimes a_1^{\LL} \LL\RR a_2^{\RR} a_3^{\RR}
		+a_1^{\LL}\LL\RR a_2^{\RR} a_3^{\RR}\otimes \LL\LL\RR\RR\RR\\
\hphantom{\vcpS\big(a_1^{\LL}\LL\RR a_2^{\RR} a_3^{\RR}\big)=}{}		 + a_1^{\LL}\LL\RR \otimes \LL\LL\RR a_2^{\RR} a_3^{\RR}
		+ \LL\RR a_2^{\RR}\otimes a_1^{\LL}\LL\RR\RR a_3^{\RR}
		+ \LL\RR a_3^{\RR}\otimes a^{\LL}\LL\RR a_2^{\RR}\RR\\
\hphantom{\vcpS\big(a_1^{\LL}\LL\RR a_2^{\RR} a_3^{\RR}\big)=}{} + a_1^{\LL}\LL\RR a_2^{\RR} \otimes \LL\LL\RR\RR a_3^{\RR}
		+ a_1^{\LL}\LL\RR a_3^{\RR}\otimes \LL\LL\RR a_2^{\RR}\RR
		+ \LL\RR a_2^{\RR} a_3^{\RR}\otimes a_1^{\LL}\LL\RR\RR\RR.
	\end{gather*}
\end{Example}

\subsection{Double monoid of words on random variables}\label{sec:bimonoid}
As for translucent words, we also write $w^{\stdl i}:=\restriction{w}{\{x:x\stdl i\}}$ and $w^{\stdg i}:=\restriction{w}{\{x:x\stdg i\}}$ for an incomplete word $w$ and an integer $i\in\interval{|w|}$.
\begin{Definition}[horizontal semigroup]	\label{ominusproduct}
	We define the structure of a $\ominus$-semigroup on $\IWc$,
	\begin{equation*}
		\mS\colon \ \CC\IWc\ominus\CC\IWc \to \CC\IWc,
	\end{equation*}
	as follows. Pick a translucent word $t \in \SWc$ with $[t]_{0}\neq \varnothing$, an integer $i \in [t]_0$ and \[
	w^{-}\ominus w^{+} \in (\CC\IWc\ominus\CC\IWc) (t)
	\]
	such that
	$w^{-}\in\IWc(t^{\stdl i})$, $w^{+}\in\IWc(t^{\stdg i})$. Then, $\mS(w^{-}\ominus w^{+})$ is the unique incomplete word in~$\IWc(t)$ satisfying
	\[\mS(w^{-} \ominus w^{+})^{\stdl i}=w^{-},\qquad \mS(w^{-}\ominus w^{+})(i)=\alpha_t(i),\qquad \mS(w^{-}\ominus w^{+})^{\stdg +}=w^{+}. \]
\end{Definition}
Further below, we prove the associativity of $\mS$.
The above definition looks cumbersome, but it is not difficult to grasp, see Example~\ref{ex:horizontalcompoun}.

\begin{Remark}\label{rk:cmphorieontalmult}
For a tensor $w^{-}\ominus w^{+} \in \IWc\ominus \IWc$, we replace the symbol $\otimes$ by $\ominus$ to emphasize the belonging to $\IWc\ominus \IWc$.

Be aware that given two incomplete words $w^{-}$ and $w^{+}$ computing the horizontal product $\mS(w^{-}\ominus w^{+})$ between those two words is, in general, meaningless.
In fact, a pair of incomplete words $w^-$, $w^+$ does not (in general) yield an element of the tensor product $\CC\IWc\ominus\CC\IWc$. For this to be true, the types $t_{w^-}$ and $t_{w^+}$ of $w^-$ and $w^{+}$ respectively should be {\it compatible}, in the sense that if $\alpha_{t_{w^-}}$ contains a letter ${\RR}$ then $\alpha_{t_{w^{+}}}(k)=\RR$ for all $k\in \interval{|w^+|}$, i.e., $\alpha_{t_{w^+}}$ is a word on just the letter $\RR$ (note that this is equivalent to the symmetric condition $\alpha_{t_{w^+}}(j)=\LL\implies \alpha_{t_{w^-}}\in\{\LL\}^\star$).

Even in the case where the words \emph{are} in fact compatible,
one needs to specify the translucent word $t$ such that $w^{-}\ominus w^{+} \in (\IWc\ominus \IWc) (t)$.
For example, take $w^{-}={\LL}a^{\LL}$ and $w^{+}=b^{\RR}$. The two words $w^{-}$ and $w^{+}$ are compatible. Then $w^{-}\ominus w^{+}$ can be either interpreted as an element in $(\IWc\ominus \IWc)(\LL\LL\LL\RR)$ or in $(\IWc\ominus \IWc)(\LL\LL\RR\RR)$. In the following, we will always be specific about which component of the collection $\IWc\ominus \IWc$ an element $w^{-}\ominus w^{+}$ should be considered as part of.
\end{Remark}
\begin{Remark}
	In Definition~\ref{ominusproduct}, suppose, for example, that $i=1$, which implies $t^{\stdl i}=\varnothing$ and~$w^{-}$ is the empty word. Then, $\mS(w^{-}\ominus w^{+})=\alpha_t(1)w^{+}$.
\end{Remark}

\newcommand\us[2]{\underset{#1}{#2}}

\begin{Example}
	\label{ex:horizontalcompoun}
	Let $w^- = \LL\RR a^{\LL}a^{\RR}$ and $w^{+}={b^{\RR}{\sf RR}}$. These words are compatible. We consider
	\[
	w^{-}\ominus w^{+} \in (\CC\IWc\ominus\CC\IWc)(t),\qquad
	t := ({\sf \RR\LL\RR\RR\RR\RR\LL\RR}, 10000011).
	\]
	Denote $\dot k$ the $k$th element in $\interval{|t|}=\{1,\dots, 8\}$ with respect to necklace order $\stdleq$, i.e.,
	\[\dot 1=2, \qquad \dot2=7,\qquad \dot 3=8,\qquad \dot 4=6,\qquad \dot 5= 5,\qquad \dot 6=4,\qquad \dot 7 =3, \qquad \dot8=1. \]

	The integer $i \in [t]_0$ of Definition~\ref{ominusproduct} is then given by $5$; indeed, the word $w^-\ominus w^ +$ is composed~-- in necklace order~-- of the four letters of $w^ -$ in positions $\dot 1$, $\dot 2$, $\dot 3$, $\dot4$, then a placeholder in position $i=\dot5=5$, and finally $w^ +$ on positions $\dot6$, $\dot7$, $\dot8$. Note how one can always deduce the position $i$ just from knowing the inputs, $w^-$, $w^+$ and $t$.

	To find the word $\mS(w^ -\ominus w^ +)$, we also have to remember that the letters of $w^ -$ and $w^ +$ are given in natural order, not in necklace order. To keep track of the different labellings of the positions according to natural and necklace order, it is best to draw a table with the following rows:
	\begin{itemize}\itemsep=0pt
	 \item $j\in \interval{|t|}$ in natural order,
	 \item the letters of $\alpha_t$ in natural order,
	 \item the corresponding necklace order positions $\dot k$ with $\dot k=j$,
	 \item the letters of $w^-\ominus w^+$: \begin{itemize}\itemsep=0pt
	 \item letters of $w^-$ in natural order (i.e., from the left in the table) to the first necklace order positions (here: $\dot1,\dots, \dot 4$ in natural order $\dot1=2<\dot4=6<\dot2=7<\dot3=8$),
	 \item $\alpha_t(i)$ in position $i$ (here: $\dot5$),
	 \item letters of $w^+$ in the natural order to the final positions in necklace order (i.e., the remaining open positions; here: $\dot8=1<\dot7=3<\dot6=4$).
	 \end{itemize}
	\end{itemize}
	In the concrete example, we obtain, using colors to visually identify $w^-$ and $w^+$ inside $\mS(w^{-}\ominus w^{+})$:
 \begin{align}\label{eq:array-for-m-ominus}\hspace*{32mm}
	\begin{array}{|c|c|c|c|c|c|c|c|c|}
			\hline
			\vrule height12pt width0pt depth7pt
			\text{nat.\ pos.}&1&2&3&4&5&6&7&8\\
			\hline
			\vrule height12pt width0pt depth7pt
			\alpha_{t}&{\RR} & {\LL} & {\RR} & {\RR} & {\RR} & {\RR} & {\LL} & {\RR} \\
			\hline
			\vrule height12pt width0pt depth7pt
			\text{neckl.\ pos.} & \dot 8& \dot 1&\dot 7&\dot6&\dot5&\dot4&\dot2&\dot3 \\
			\hline
			\vrule height12pt width0pt depth7pt
			\mS(\textcolor{blue}{w^{-}}\ominus \textcolor{red}{w^{+}})&\textcolor{red}{b^\RR} & \textcolor{blue}{\LL} & \textcolor{red}{\RR} & \textcolor{red}{\RR} & \RR & \textcolor{blue}{\RR} & \textcolor{blue}{a^\LL} & \textcolor{blue}{a^\RR} \\
			\hline
		\end{array}
 \end{align}
\end{Example}
Owing to Propositions~\ref{prop:exchange} and \ref{prop:duoid}, the vertical composition $\IWc\overt\IWc$ is endowed with a~horizontal semigroupal product:
\begin{equation*}
	\mSS\colon \ (\CC\IWc\overt \CC\IWc) \ominus (\CC\IWc\overt \CC\IWc) \to \CC\IWc\overt \CC\IWc.
\end{equation*}
\begin{Example}
	\label{ex:horizontacompodeux}
	Consider the words $w_1^- = \LL \RR a^{\LL}$, $w_2^- = \LL \RR \LL a^{\RR}$. The two words $w^{-}_1$ and $w^{-}_2$ are composable and yield $w^-:=w^-_1 \circ w^{-}_2 = \LL\RR a^{\LL} a^{\RR}$. Likewise for the following words:
	$w^{+}_1 = b^{\RR} {\RR}$, $w^{+}_2=\RR \RR b^{\RR}$, $w^+:=w^{+}_1\circ w^{+}_2 = b^{\RR}\RR b^{\RR}$. Since $w^-\LL\RR a^{\LL} a^{\RR}$ and $w^+=b^{\RR}\RR b^{\RR}$ are compatible words, the above four words yield an element in{\samepage
	\[\bigl((\CC\IWc\overt \CC\IWc) \ominus (\CC\IWc\overt \CC\IWc)\bigr)(t)\]
	with $t =({\sf RLRRRRLR}, 10000011)$.}

Recall that, by definition,
 \[\mSS((w^{-}_1\overt w^{-}_2)\ominus (w^{+}_1 \overt w^{+}_2))=\mS(w^{-}_1\ominus w^{+}_1)\overt \mS(w^{-}_2\ominus w^{+}_2)\in (\CC\IWc\overt\CC\IWc)(t).\] We draw a table similar to \eqref{eq:array-for-m-ominus}, but the composable pairs $w_1^-\overt w_2^-$ and $w_1^+\overt w_2^+$ are written in two separate rows. We use colours again so that $w^-$ and $w^+$ can be easily recognized and also write vertical and horizontal tensor products simply as superposition and juxtaposition.
$$\begin{array}{|c|c|c|c|c|c|c|c|c|}
			\hline
			\vrule height12pt width0pt depth7pt
			\text{nat.\ pos.}&1&2&3&4&5&6&7&8\\
			\hline
			\vrule height12pt width0pt depth7pt
			\alpha_t&{\RR} & {\LL} & {\RR} & {\RR} & {\RR} & {\RR} & {\LL} & {\RR} \\
			\hline
			\vrule height12pt width0pt depth7pt
			\text{neckl.\ pos.} & \dot 8& \dot 1&\dot 7&\dot6&\dot5&\dot4&\dot2&\dot3 \\
			\hline
			\vrule height12pt width0pt depth7pt
			\multirow{2}{*}{$\mSS\left(
 \begin{matrix}
 \textcolor{blue}{w^{-}_1}\\\textcolor{blue}{w^{-}_2}
 \end{matrix}
 \middle\vert
 \begin{matrix}
 \textcolor{red}{w^{+}_1} \\ \textcolor{red}{w^{+}_2}
 \end{matrix}
 \right)$}	&\textcolor{red}{b^\RR} & \textcolor{blue}{\LL} & \textcolor{red}{\RR} & & \RR & \textcolor{blue}{\RR} & \textcolor{blue}{a^\LL} & \\	&\textcolor{red}{\RR} & \textcolor{blue}{\LL} & \textcolor{red}\RR & \textcolor{red}{b^\RR} & \RR & \textcolor{blue}{\RR} & \textcolor{blue}{\LL} & \textcolor{blue}{a^\RR} \\
\hline
 \end{array}$$
 The last row contains $\mS(w_2^-\ominus w_2^+)\in (\CC\IWc \ominus \CC\IWc)(t)$ and is found just as before in Example~\ref{ex:horizontalcompoun}. In the row showing $\mS(w_1^-\ominus w_1^+)\in (\CC\IWc \ominus \CC\IWc)([\mS (w_2^-\ominus w_2^+)]_0) $, the entries are obtained by the same rules, but restricted to the positions in which $\mS(w_2^-\ominus w_2^+)$ is translucent, i.e., those positions in which $\mS(w_2^-\ominus w_2^+)$ has a placeholder. In our usual notation, the result reads \begin{equation*}
		\mSS((w^{-}_1\overt w^{-}_2)\ominus (w^{+}_1 \overt w^{+}_2)) = b^{\RR} \LL \RR \RR \RR a^{\LL} \overt \RR \LL \RR b^{\RR} \RR \RR \LL a^{\RR}.
 \end{equation*}
 From the table, one immediately sees the composability of the obtained words, and could also easily read off their composition.
\end{Example}

\begin{Proposition}\label{prop:freealgebra} The horizontal semigroup of incomplete words is freely generated by its opaque part, that is
	\begin{equation*}
		\IWc \simeq \mathcal{F}(U(\IWc))
	\end{equation*}
with the adjoint functors $\mathcal F$ and $U$ as defined in Section~{\rm \ref{subsec:freely-generated-ominus-semigroups}}.
\end{Proposition}
\begin{proof}Pick a translucent word $t \in \SWc$ and a non-empty incomplete word on random variables $w \in \IWc(t)$. Let $[t]_0 = \{i_1 \stdl \cdots \stdl i_p\}$.
	Denote by $w_j$ the restriction of $w$ to $\rrbracket i_{j-1},i_j\llbracket_\stdl$, for $1\leq j \leq p+1$, with the convention that $i_0=-\infty_\stdl$, $i_{p+1}=+\infty_\stdl$ (see Notation~\ref{notation:infinity}) and set
	\[\Phi(w) := w_0\ominus \cdots \ominus w_{p+1} \in \mathcal{F}(U(\IWc))(t).\]

	Then $\Phi$ yields in injective morphism $\Phi\colon \IWc\to\mathcal{F}(U(\IWc))$ of linear $\SWc$-collections. The inverse morphism $\Phi^{-1}$ satisfies
	\begin{equation*}
		\Phi^{-1}(w_1\ominus \cdots \ominus w_p) = \mS(w_1\ominus \cdots \ominus w_p),
	\end{equation*}
	where $w_1 \ominus \cdots \ominus w_p \in \mathcal{F}(U(\IWc))(t)$.
\end{proof}

In the following, we use the shorter notation
\[
\mathcal{R} := \big(\mS \overt \mS\big) \circ R_{\IWc,\IWc,\IWc,\IWc}.
\]
\begin{Observation}
 Pick $w^-\ominus w^+\in (\IWc\ominus \IWc)(t)$ and set $w = \mS(w^{-}\ominus w^{+})$.
 The map $\mathcal{R}$ yields a \emph{bijection} between the set of admissible cuts of $w$ and pairs of an admissible cut of $w^-$ and an admissible cut of $w^+$.
\end{Observation}

\begin{Example}
	Take $t=({\sf RLRR{R}RL},1001001)$, $i=5$, $w=b^{\RR} \LL \RR b^{\RR} {\RR}\RR a^{\LL}$, $w^- = \LL\RR a^{\LL}$, $w^{+}={b^{\RR}{\sf Rb^{\RR}}}$. In the following table, we listed cuts of $w$ and pairs of a cut of $w^-$ and $w^+$ that correspond through $\mathcal{R}$.
$$\renewcommand{\arraystretch}{1.3}
		\begin{array}{|l|c|c|c|}
			\hline
			 {\rm Adm}(w) & \LL \RR b^{\RR} \RR \RR a^{\LL}, \ b^{\RR} \LL\RR\RR\RR\RR\LL & b^{\RR} \LL\RR\RR\RR a^{\LL}, \ \RR \LL\RR b^{\RR}\RR\RR\LL & b^{\RR}\LL\RR b^{\RR}\RR\RR, \ \RR\LL\RR\RR\RR\RR a^{\LL} \\
			\hline
			 {\rm Adm}(w^{+}) & \RR b^{\RR}, \ b^{\RR} \RR \RR & b^{\RR} \RR, \ \RR \RR b^{\RR} & b^{\RR} \RR b^{\RR}, \ \RR\RR\RR \\
			\hline
			 {\rm Adm}(w^-) & \LL\RR a^{\LL}, \ \LL\RR\LL & \LL\RR a^{\LL}, \ \LL\RR\LL & \LL\RR, \ \LL\RR a^{\LL} \\
			\hline
		\end{array}
$$
\end{Example}

\begin{Proposition}	\label{cor:vcpominusmorphism}
	\begin{equation*}
		\vcpS \circ \mS = \mSS\circ\big(\vcpS \ominus \vcpS\big).
	\end{equation*}
\end{Proposition}
\begin{proof} Pick compatible words with $w^-\ominus w^+\in (\IWc\ominus\IWc)(t)$ and set $w=\mS(w^-\ominus w^{+})$. Then one has
	\begin{align*}
		\vcpS(w)=\sum_{(\ell,u) \in {\rm Adm}(w)} \ell \otimes u &= \sum_{\substack{(\ell^{-},u^{-}) \in {\rm Adm}(w^-) \\ (\ell^{+},u^{+}) \in {\rm Adm}(w^{+})}} \mS((\ell^-\otimes u^-) \ominus (\ell^+ \otimes u^+)) \\
		&= \mSS\big(\vcpS(w^{-})\ominus \vcpS(w^{+})\big)
	\end{align*}
	as needed.
\end{proof}

\subsection{Unshuffle structure}
We now come to the core of the present article: we show that $\vcpS$ introduced in Section~\ref{ss:DmonoidDcomonoid} splits into two parts $\vcp_{\prec}$ and $\vcp_{\succ}$ satisfying the co-shuffle relations.

We start with a short account of shuffle algebras and related structures. The abstract notion of a shuffle product and its decomposition into two nonassociative products goes back to the work of Eilenberg, MacLane, and Sch\"utzenberger in the 1950s. For introductory materials on shuffle algebras, we refer the reader to \cite[Section~5]{loday1763dialgebras}. To avoid confusion, we use $\otimes_\CC$ for the tensor product over the category of vector spaces.

\begin{Definition}\label{def:unshuffle}
	A {\it counital unshuffle coalgebra} is a coaugmented coassociative coalgebra $C = \CC\oplus \bar{C}$ with coproduct
	\begin{equation*}
		\Delta(c):= 1 \otimes_{\CC} c + c \otimes_{\CC} 1 + {\bar{\Delta}}(c)
	\end{equation*}
	such that the {\it reduced coproduct}
	on $\bar{C}$ splits, $\bar{\Delta} = \Delta_{\prec} + \Delta_{\succ}$ with
	\begin{align*}
		&(\Delta_{\prec} \otimes_{\CC} \mathrm{id}) \circ \Delta_\prec = (\mathrm{id} \otimes_{\CC} \bar{\Delta})\circ \Delta_{\prec}, \\
		&(\Delta_{\succ} \otimes_{\CC} \mathrm{id}) \circ \Delta_{\prec} = (\mathrm{id} \otimes_{\CC} \Delta_{\prec})\circ \Delta_{\succ}, \\
		&(\bar{\Delta} \otimes_{\CC} \mathrm{id}) \circ \Delta_{\succ} = (\mathrm{id} \otimes_{\CC} \Delta_{\succ})\circ \Delta_{\succ}.
	\end{align*}
\end{Definition}
The maps $\Delta_{\prec}$ and $\Delta_{\succ}$ are called respectively the left and right half-unshuffle coproducts.

\begin{Definition}A {\it unital shuffle algebra} $A$ is an augmented algebra, which means \[{A} = \CC\oplus \bar{A},\]
	where $\bar{A}$ is a two-sideal ideal of $A$ equipped with an algebra product
	$\centerdot$, which splits into two (nonassociative) bilinear products $\prec$ and $\succ$ on $\bar{A}$,
	\begin{equation*}
		a \centerdot b = a \prec b + a \succ b,\qquad a,b \in \bar{A},
	\end{equation*}
	satisfying the {\it shuffle relations}, with $a,b,c \in \bar{A}$,
	\begin{align*}
		 (a \prec b) \prec c = a \prec (b\centerdot c),\qquad
		a \succ (b \succ c) = (a \centerdot b) \succ c , \qquad
		 (a \succ b) \prec c = a \succ ( b \prec c).
	\end{align*}
\end{Definition}

In the following paragraphs, Definition~\ref{def:unshuffle} undergoes a very light ``internalization''\footnote{See, for example, the corresponding article on nLab~\cite{nlab:internalization}.} procedure yielding the definition of a coaugmented conilpotent unshuffle coalgebra but in the category of linear $\SWc$-collections. To keep the presentation contained, and because the study of unshuffle algebras in the category of $\SWc$-indexed collections in full generality is not the core of the present work, we restrain from giving a general definition of such an object (but we give an example of it). The reader will, without much effort, be able to extract it from Theorem~\ref{thm:unshuffle}.

We now turn to the definition of the \emph{unshuffle structure on incomplete words on random variables.}

Recall, \eqref{eq:Xplus},
that we denote by ${\IWc^+}$ the $\SWc$-indexed collection with ${\IWc^+}(t)={\IWc(t)}$ if $b_t \neq \mathbf{0}$ and $\IWc^+(t) = \varnothing$ if $b_t = \mathbf{0}$.

For $w\in \IWc^+$
and a set
$A \subset \interval{|w|}$, we denote by $\bigwedgedot A$ the minimum in $A$ with respect to the necklace order $\stdl$.
In particular, $\bigwedgedot [w]_1$ is the position of the first letter of $w$ in necklace order.
\begin{Definition}
	We define two $\SWc$-collection morphisms, called left and right {\it half-unshuffle coproduct}, respectively,
	\[\vcpSpr, \vcpSsu\colon \ \linearspan{\IWc^{+}} \to \linearspan{\IWc^+} \overt \linearspan{\IWc^+}\] by, for $w\in \IWc^+$,
	\begin{equation*}
		\vcpSpr(w)
		= \sum_{\substack{I \subsetneq \interval{|w|} \\ [w]_0 \subsetneq I \\ \bigwedgedots [w]_1 \in I}} \restriction{w}{I} \otimes \translucidation{w}{I},\qquad
 \vcpSsu(w) = \sum_{\substack{I \subsetneq \interval{|w|} \\ [w]_0 \subsetneq I \\ \bigwedgedots [w]_1 \not\in I}} \restriction{w}{I} \otimes \translucidation{w}{I}.
	\end{equation*}
	We also put $\barvcpS=\vcpSpr + \vcpSsu$ and call it {\it reduced unshuffle coproduct}.
\end{Definition}

\begin{Example}
	$\bulletAccent S = \big\{a^{\LL},b^{\LL}\big\}\sqcup \big\{a^{\RR},b^{\RR}\big\}$,
	then, for example,
	\begin{gather*}
		\vcpSpr\big( \cyan{a^{\LL}}{\LL}b^{\LL}b^{\RR} \big) 
		 =
		\cyan{a^{\LL}}{\LL} \otimes {\LL} {\LL} b^{\LL} b^{\RR}
		+
		\cyan{a^{\LL}}{\LL} b^{\LL} \otimes {\LL} {\LL} L b^{\RR}
		+
		\cyan{a^{\LL}}{\LL} b^{\RR} \otimes {\LL} {\LL} b^{\LL} ,\\
		\vcpSsu\big( \cyan{a^{\LL}}{\LL} b^{\LL} b^{\RR} \big)
		 =
		\LL b^{\LL} b^{\RR} \otimes \cyan{a^{\LL}}\LL \LL \RR
		+
		\LL b^{\LL} \otimes \cyan{a^{\LL}}\LL \LL b^{\RR}
		+
		\LL b^{\RR} \otimes \cyan{a^{\LL}}\LL b^{\LL} R,
	\end{gather*}
	where, for legibility, we highlight the letter at position $\bigwedgedot[w]_1$.
 \end{Example}

\begin{Theorem}\label{thm:unshuffle}
	The linear $\SWc$-indexed collection $\linearspan{\IWc}$ equipped with the coproduct $\vcpS$ is an \emph{unshuffle conilpotent coaugmented comonoid}, by which we mean the following.

	First, there exist two $\SWc$-indexed collection morphisms $\varepsilon\colon \IWc\to \uvert$ and $\eta\colon \uvert \to \IWc$, named co-unit and co-augmentation map, defined by
\begin{gather*}
\varepsilon_{(\alpha,\mathbf{0})}(\alpha) = \Singleton \qquad \text{and} \qquad \varepsilon_{t}= 0 \text{ for } t\neq (\alpha,\mathbf{0}), \\
\eta_{(\alpha,\mathbf{0})}(\raisebox{-1pt}{\scalebox{0.7}{\PHcat}}) = \alpha.
	\end{gather*}
	and satisfying
	\begin{align*}
 (\varepsilon \overt \mathrm{id})\circ\vcp = \mathrm{id},\qquad (\mathrm{id}\overt \varepsilon)\circ\vcp = \mathrm{id},\qquad \varepsilon \circ \eta = \mathrm{id}.
	\end{align*}
	Second, the coproduct $\vcpS$ splits for $w\in \mathcal W^+$ as
	\begin{equation*}
		\vcpS(w) = w \overt \eta(\Singleton) + \eta(\Singleton)\overt w + \barvcpS(w).
	\end{equation*}
	Third, the unshuffle relations are satisfied,
	\begin{gather*}
		 \big(\vcpSpr \overt \mathrm{id}\big)\circ \vcpSpr = \big(\mathrm{id} \overt \barvcpS\big) \circ\vcpSpr,\qquad
		\big(\mathrm{id} \overt \vcpSsu\big)\circ\vcpSsu = \big(\barvcpS \overt \mathrm{id}\big)\circ\vcpSsu, \\
		 \big(\vcpSsu \overt \mathrm{id}\big)\circ\vcpSpr = \big(\mathrm{id} \overt {\vcpSpr}\big)\circ\vcpSsu.
	\end{gather*}
	Finally, $\barvcpS$ is nilpotent in the sense that for every $w\in \mathcal W^+$, there is an $n$ such that ${\barvcpS}^n=0$ $\big($defined recursively, ${\barvcpS}^1=\barvcpS$, ${\barvcpS}^{k+1}:=\big(\barvcpS\overt \mathrm{id}\big){\barvcpS}^{k}\big)$.
\end{Theorem}
\begin{proof}The only statements for which proofs are needed are the unshuffle relations and the nilpotency of $\bar{\Delta}^{\overt}_{\mathcal{W}}$. Nilpotency follows from the fact that if $w\in \mathcal{W}^+$, then $\bar{\Delta}^{\overt}(w)$ is a sum of tensor $w^{\prime}\overt w^{\prime}$ where $w$ and $w^{\prime}$ are translucent words with strictly fewer of opaque letters than~$w$ has. Iterating, we eventually end up applying $\bar{\Delta}^{\overt}$ on some $\eta(\Singleton)$, which gives~$0$.

We now turn to proving the unshuffle relations.
Pick $w$ an incomplete word. Pick a pair of subsets $[w]_0 \subset X_1 \subset X_2 \subset \interval{|w|}$. We denote by $k\mapsto \overline k\colon X_1 \to \interval{|X_1|}$ the unique strictly increasing bijection (i.e., the inverse of the map $k\mapsto i_k$ for $X_1=\{i_1\stdl\cdots \stdl i_p\}$). Then,
\[
\restriction{(\translucidation{w}{X_2})}{X_1} = \translucidation{(\restriction{w}{X_1})}{\overline{X}_2}.
\]
Coassociativity of $\vcpS$ entails
\begin{gather*}
\sum_{\substack{[w]_0 \subset X_2 \subset X_1}}\!\!\!\!\restriction{w}{X_2} \otimes \restriction{(\translucidation{w}{X_2})}{X_1} \otimes \translucidation{w}{X_1} = \sum_{\substack{[w]_0 \subset X_1 \\ \overline{[w]}_0\subset\overline{X}_2 \subset \interval{|X_1|}}} \!\!\!\!\restriction{(\restriction{w}{X_1})}{\overline{X}_2}\otimes \translucidation{(\restriction{w}{{X}_1})}{\overline{X}_2}\otimes\translucidation{w}{X_1}.
\end{gather*}
Now, $\bigwedgedot [w]_1 \in X_2$ if and only if (with obvious notation) $\overline{\bigwedgedot {[w]_1}} \in \overline{X}_2$ and moreover, $\bigwedgedot [\restriction{w}{X_1}]_1 = \overline{\bigwedgedot {[w]_1}}$.

This implies the following chain of equalities,
\begin{align*}
		 \big(\vcpSpr \overt \mathrm{id}\big) \circ \vcpSpr(w) &=
		 \sum_{\substack{[w]_0 \subsetneq X_1 \\ \overline{[w]}_0\subsetneq\overline{X}_2 \subsetneq \interval{|X_1|} \\ \bigwedgedots [\restriction{w}{X_1}]_1 \in X_2 }} \restriction{(\restriction{w}{X_1})}{\overline{X}_2}\otimes \translucidation{(\restriction{w}{{X}_1})}{\overline{X}_2}\otimes\translucidation{w}{X_1}\\
		 &=\sum_{\substack{[w]_0 \subsetneq X_2 \subsetneq X_1\\\bigwedgedots [w]_1 \in X_2}} \restriction{w}{X_2} \otimes \restriction{(\translucidation{w}{X_2})}{X_1} \otimes \translucidation{w}{X_1} \\
		 &= \big(\mathrm{id} \overt \barvcpS\big) \circ \vcpSpr(w).
\end{align*}
We prove likewise the remaining two identities. If $\bigwedgedot [w]_1 \in X_2^c$, then $\bigwedgedot [\translucidation{w}{X_2}]=\bigwedgedot [w]_1$. This alone implies
\begin{align*}
		 \big(\mathrm{id} \overt \vcpSsu\big)\circ\vcpSsu(w)
		 &=\sum_{\substack{[w]_0 \subsetneq X_2 \subsetneq X_1 \\ \bigwedgedot [w]_1 \in X_1^c}} \restriction{w}{X_2} \otimes \restriction{(\translucidation{w}{X_2})}{X_1} \otimes \translucidation{w}{X_1} \\
		 &= \sum_{\substack{[w]_0 \subsetneq X_1 \\ \overline{[w]}_0\subsetneq\overline{X}_2 \subsetneq \interval{|X_1|} \\ \bigwedgedots [w]_1 \in X_1^c}} \restriction{(\restriction{w}{X_1})}{\overline{X}_2}\otimes \translucidation{(\restriction{w}{{X}_1})}{\overline{X}_2}\otimes\translucidation{w}{X_1} \\
		 &= \big(\barvcpS \overt \mathrm{id}\big)\circ\vcpSsu(w).
\end{align*}
Finally, if $\bigwedgedot[w]_1\in X_1$, since $\bigwedgedot[\restriction{w}{{X}_1}]_1 = \overline{\bigwedgedot {[w]_1}}\in \llbracket |X_1| \rrbracket \setminus \overline{X}_2$ is equivalent to $\bigwedgedot [w]_1 \in X_1\setminus X_2$, one gets
\begin{align*}
		 \big(\vcpSsu \otimes \mathrm{id}\big)\circ\vcpSpr(w) &=
		 \sum_{\substack{[w]_0 \subsetneq X_1 \\ \overline{[w]}_0\subsetneq\overline{X}_2 \subsetneq \interval{|X_1|} \\ \bigwedgedots [\restriction{w}{X_1}]_1 \in \llbracket |X_1| \rrbracket \setminus \overline{X}_2 }} \restriction{(\restriction{w}{X_1})}{\overline{X}_2}\otimes \translucidation{(\restriction{w}{{X}_1})}{\overline{X}_2}\otimes\translucidation{w}{X_1}\\
		 &=\sum_{\substack{[w]_0 \subsetneq X_2 \subsetneq X_1 \\ \bigwedgedots [w]_1 \in X_1 \setminus X_2}} \restriction{w}{X_2} \otimes \restriction{(\translucidation{w}{X_2})}{X_1} \otimes \translucidation{w}{X_1}\\
		 &= \big(\mathrm{id} \otimes {\vcpSpr}\big)\circ \vcpSsu(w).
	\end{align*}
	This ends the proof.
\end{proof}

\begin{Proposition} One has
	\label{prop:distributivite}
	\begin{equation*}
		\vcp_{\prec,\succ}\circ \mS = \mSS\circ\big(\vcp_{\prec,\succ} \ominus \vcpS\big).
	\end{equation*}
\end{Proposition}

\begin{proof}The proof follows the same line of arguments as in Corollary~\ref{cor:vcpominusmorphism}
\end{proof}

We now turn our attention to the dual structures.
Given two linear functionals
\[\ell_{1}, \ell_{2}\colon \ \bigoplus_{t\in\SWc} \CIWc(t) \to \CC \]
we define their (full shuffle) convolution product $\ell_1\star\ell_2$ by
\begin{equation*}
	\ell_1 \star \ell_2 = m_{\CC}\circ (\ell_1 \overt \ell_2)\circ \vcpS,
\end{equation*}
where $m_\CC$ is the product of complex numbers.

Since the reduced coproduct $\bvcp$ splits into two parts, the bilinear product $\star$ splits as well into two nonassociative products.

Note that we identify a linear functional in $(\linearspan{\IWc})^{\ast}$ with $\ell\circ\eta = 0$ with the functional in~$(\linearspan{\IWc^{+}})^{\ast}$ obtained by restriction.

\begin{Definition} For two linear functionals $\ell_1, \ell_2 \in (\linearspan{\IWc^{+}})^{\ast}$, we define the left and right {\it half-shuffle convolution products}
	\begin{equation*}
		\ell_1 \prec \ell_2 := m_{\CC} \circ (\ell_1 \overt \ell_2) \circ \vcpSpr,\qquad \ell_1 \succ \ell_2 := m_{\CC} \circ (\ell_1 \overt \ell_2) \circ \vcpSsu.
	\end{equation*}
\end{Definition}

\begin{Corollary}The space of linear functionals which are constant on the coaugmentation part of $\CIWc$,
	\[\mathbb C\varepsilon \oplus (\linearspan{\IWc^{+}})^{\ast}=\{\ell (\linearspan{\IWc^{+}})^{\ast}\mid \ell\circ\eta\in\mathbb C\varepsilon\},\] is a unital shuffle algebra with respect to the splitting
	\[\ell_1\star\ell_2=\ell_1\prec \ell_2 + \ell_1\succ\ell_2\]
	of the shuffle convolution into left and right half-shuffle convolution.
Furthermore, for every $w\in (\linearspan{\IWc^{+}})$, there exists an $n_0\in\mathbb N$ such that $(\ell_1\star\cdots\star \ell_n)(w)=0$ for all $n>n_0$ and all $\ell_1,\dots, \ell_n\in (\linearspan{\IWc^{+}})^\ast$.
\end{Corollary}
\begin{proof}
 From $\vcpS(\eta(\Singleton))=\eta(\Singleton)\overt\eta(\Singleton)$, we easily conclude that $(\linearspan{\IWc^{+}})^*$ is a two-sided ideal. The other claims follow immediately from Theorem~\ref{thm:unshuffle}.
\end{proof}

The two products $\prec$ and $\succ$ are partially extended to $\mathbb C\varepsilon \oplus (\linearspan{\IWc^{+}})^{\ast}$.
This amounts to ascribing meaning to $\ell \prec \varepsilon$ and $\varepsilon \succ \ell$ with $\ell \in {(\linearspan{\IWc}^{+})}^{\ast}$ by setting
\begin{equation*}
	\ell \prec \varepsilon = \ell,\qquad \varepsilon \succ \ell = \ell.
\end{equation*}
and $\varepsilon \prec \ell = 0 = \ell \succ \varepsilon$.
Notice that the expressions $\varepsilon \prec \varepsilon$ and $\varepsilon \succ \varepsilon$ are not defined.

The reader is directed to \cite{manchon2011short} for the definition of and a short introduction to preLie algebras.
\begin{Proposition}\label{prop:group-G-and-preLie-g}
	Denote by $G$ the set of linear functionals on $\linearspan{\IWc}$ equal to one on the co-augmentation part of $\linearspan{\IWc}$,
	\begin{equation*}
		G = \{\ell \colon \linearspan{\IWc}\to\CC,\,\ell \circ \eta=\varepsilon\}.
	\end{equation*}
	Then, $G$ is a group if equipped with the convolution product $\star$.
	 Besides, the set
	\begin{equation*}
		\mathfrak{g} = \{ \ell \colon \linearspan{\IWc} \to \CC,\, \ell \circ \eta = 0\}
	\end{equation*}
	is a preLie algebra for the preLie product
	\begin{equation}\label{eq:preLie-from-shuffle}
		\ell_1 \vartriangleleft \ell_2 := \ell_1 \prec \ell_2 - \ell_2 \succ \ell_1,\qquad \ell_1,\ell_2\in\mathfrak{g}.
	\end{equation}
\end{Proposition}

\begin{proof}Invertibility of $\ell\in G$ follows from the ``argumentwise nilpotency'': $\ell^{-1}(w)=\sum_{k=0}^\infty (\varepsilon-\ell)^{\star k}(w)$ where the right hand side is a finite sum.
 The second statement is a special case of the well-known fact that, for any shuffle algebra, \eqref{eq:preLie-from-shuffle} defines a preLie product, see, e.g., \cite[Section~3.1]{ebrahimi2015cumulants}.
\end{proof}

In the next two and final sections, we exhibit the two sub preLie algebras of $\mathfrak{g}$ relevant to our work. We end this section with the definition of the two ``time-ordered'' exponentials or half-shuffle exponentials
\[
	\exp_{\prec}(k) := \varepsilon + \sum_{n=1}^{\infty} k^{\prec n},\qquad \exp_{\succ}(k) := \varepsilon + \sum_{n=1}^{\infty} k^{n \succ},
\]
where $k^{\prec 1} = k$, $k^{\prec (n+1)} = k \prec k^{\prec n}$ and
$k^{1 \succ } = k$, $k^{(n+1)\succ} = k^{n \succ}\succ k$.

\section{Lie theory of two-faced independences}\label{sec:momentscumulantsrelations}
In this section, we exhibit two sub-preLie algebras of the preLie algebra $\mathfrak{g}$ introduced in the previous section in Proposition~\ref{prop:group-G-and-preLie-g}. We need a few additional items of notation.
Pick an incomplete word $w$ with length $n$. Given a subset $I \subset [w]_1$, we write $\hat{I} := I \cup [w]_0 \subset [n]$. Pick two linear forms $f,g \in \mathfrak{g}$, then
\begin{align*}
 (f\vartriangleleft g)(w) = (f\prec g - g\succ f)(w) = \sum_{\substack{I \subset [w]_1 \\ {\wedgedot}[w]_1\in I}}f(\restriction{w}{\hat{I}})g(\translucidation{w}{\hat{I}}) - \sum_{\substack{I\subset [w]_1 \\ \bigwedgedots [w]_1 \not\in I}}g(\restriction{w}{\hat{I}})f(\translucidation{w}{\hat{I}}).
\end{align*}

\begin{figure}[!ht]\centering
	\includegraphics[scale=0.45]{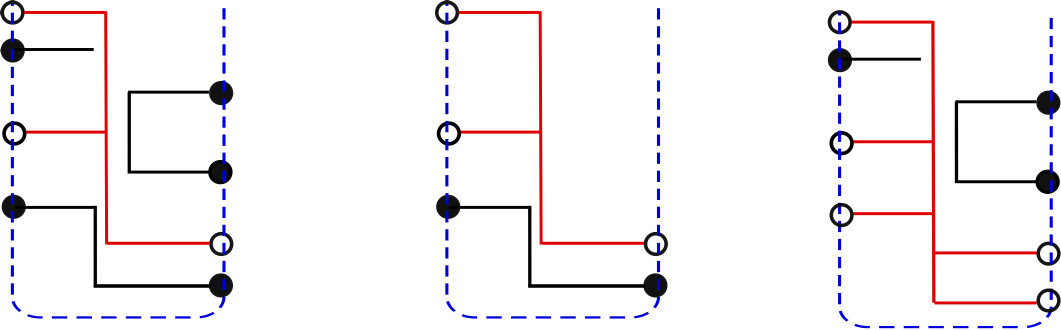}
 \caption{An incomplete bipartitition $\pi$, its restriction $\restriction{\pi}P$ and its translucidation $\translucidation{\pi}P$.
 Here, $P$~contains $[t]_0$ and the two-element block on the bottom.}
		\label{fig:exampleRestriction-diagram}
\end{figure}

Besides translucidation and restrictions of translucent words and incomplete words on random variables, it will be convenient to consider the same operations on incomplete bipartitions.
Let~$t$ be a translucent word
and $\pi \in \BPc( t )$ an incomplete bipartition.
For a subset of blocks of $P \subset \pi \cup \{[t]_0\}$, $\restriction{\pi}{P}$ is the restriction of $\pi$ to $P$ (this includes re-indexing) and $\translucidation{\pi}{P}$ is the incomplete bipartition obtained from $\pi$ by turning translucent all blocks in $\pi\setminus P$. One has
\[
 \restriction{\pi}{P} \in \BPc(\restriction{t}{{\rm supp} (P)}),\qquad \translucidation{\pi}{P} \in \BPc(\translucidation{t}{{\rm supp}(P)}),
\]
where ${\rm supp}(P):=\bigcup_{W\in P}W$.
See Figure~\ref{fig:exampleRestriction-diagram} for an example.

\subsection{preLie algebra I}\label{sec:preliealgebraI}
Recall Notation~\ref{notation:INT},
\(
	\INT{t} = \{I_1^{\stdl}(t) \stdl \cdots \stdl I^{\stdl}_{k(t)}(t)\}.
\)
For $w \in \IWc$ we write $I_j^\stdl(w) = I_j^\stdl(t_w)$ and $\INT{w} = \INT{t_w}$.

\begin{Definition}\label{def:morphismU}
	We introduce the following subsets of linear functionals (the index $\mathsf I$ stands for interval):
	\begin{itemize}\itemsep=0pt
		\item Define $\shgrI$ to be the set of all linear forms on $\linearspan{\IWc}$ that are multiplicative in the following sense,
\begin{equation}\label{eqn:definitionmprecs}
			 f \in \shgrI \Leftrightarrow f(w) =
			 \prod_{J \in \INT{w}} f(\restriction{w}{J}),\qquad w\in \IWc.
		 \end{equation}

		\item
		 Define $\shalI$ to be the set of all linear forms $f$ on $\linearspan{\IWc}$ such that
					\smallskip
		 \begin{itemize}\itemsep=0pt
			 \item
 The support of $f$ is contained in the set of words whose $[t]_1$ has only one connected component, i.e.,
 \begin{align}\label{eq:m(S):noninterval}
 \forall w \in \IWc\colon \ ( |\INT{w}| \not= 1 \Rightarrow f(w) = 0 ).
 \end{align}
			 \item
 For words $w\in\IWc(t)$ with $|\INT{w}|=1$, the value $f(w)$ only depends upon the restriction $\restriction{w}{[t]_1}$ of $w$ to the places tagged $1$:
 \begin{equation}\label{eq:m(S):interval}
 f(w) = f(\restriction{w}{[t]_1}).
 \end{equation}
		 \end{itemize}
	\end{itemize}
\end{Definition}

\begin{Remark}
	\label{rk:monoidc}
	Notice that for a linear form $f\in \shgrI$ the above formula implies $f(w)=1$ if $w$ is a word in $\{\sf L,R \}^{\star}$ since in that case, the product on the right-hand side of \eqref{eqn:definitionmprecs} is empty (thus conventionally equal to $1$). Hence, Definition~\ref{def:morphismU} is equivalent to
	\begin{equation*}
		\shgrI = \big\{\alpha \in \mathrm{Hom}_{\sgominus}\big(\big(\IWc,\hpIWc\big), \CC^{\ominus}\big)\colon \alpha(u) = 1, u \in \{\LL,\RR\}^{\star}\big\},
	\end{equation*}
	where $\CC^{\ominus}$ is the linear $\SWc$-indexed collection with $\CC^{\ominus}(t)=\CC$, $t\in\SWc$.
\end{Remark}

We denote by $\SWint \subset \SWc$ the set of all translucent words $t$ with $|\INT{t}|=1$, i.e., $[t]_1$~a~$\stdl$-interval. We also use the notation $\IWint$ for the set of all incomplete words with $|\INT{w}|=1$, i.e., those words whose translucent type lies in $\SWint$. At times we call words in $\IWint$ \emph{interval words}.

\begin{Proposition}	\label{prop:prelieproduct}
	The preLie product $\vartriangleleft$ restricts to linear space $\shalI$. Furthermore, for a word $w\in\IWc(t)$, $|\INT{w}|=1$ the following formula holds for all $f,g \in \shalI$,
	\begin{equation}\label{eqn:formulaprelie}
		(f \vartriangleleft g)(w)=
 \sum_{\substack{I_1, I_2, J \subset [t]_1\\ I_1,I_2\not= \varnothing \\ I_1 \stdl J \stdl I_2 ~ \stdl {\rm-interval}, \\ I_1\sqcup J \sqcup I_2 =[t]_1 }}
 f(\restriction{w}{{I}_1\cup {I}_2})g(\restriction{w}{J}).
	\end{equation}
\end{Proposition}

\begin{figure}[!ht]\centering
	\includegraphics[scale=0.45]{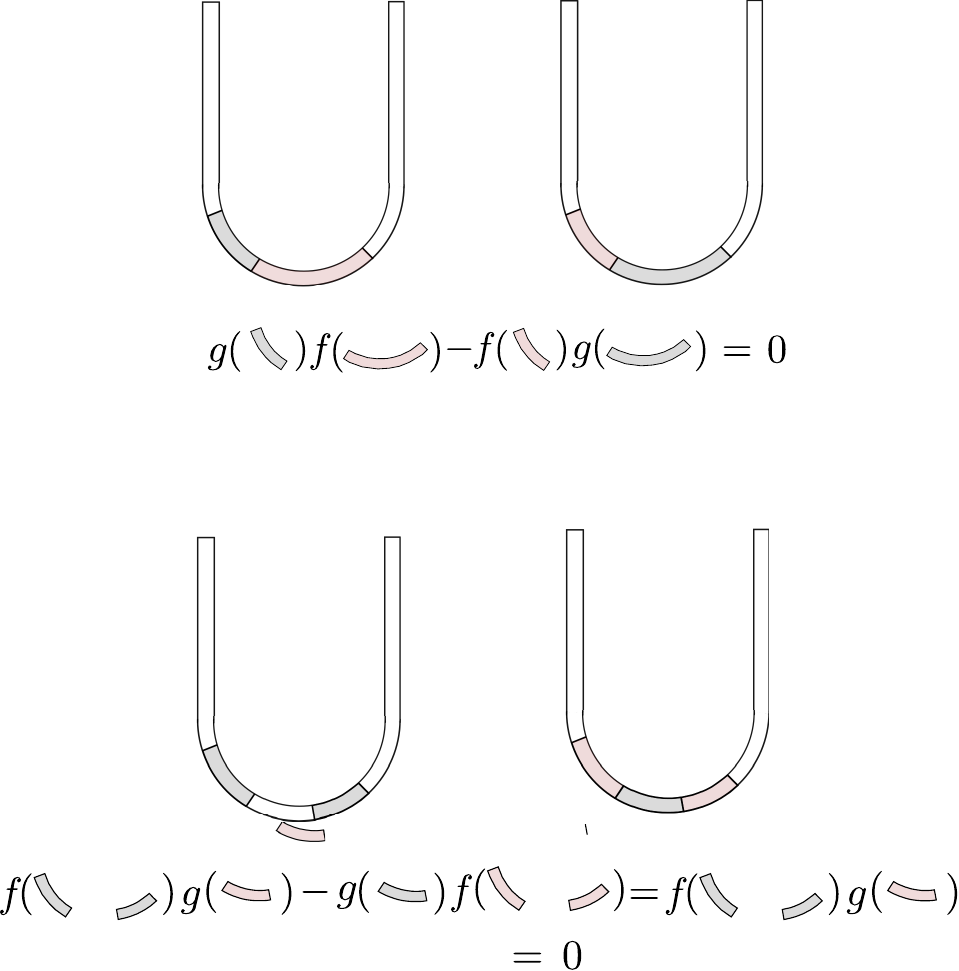}
	\caption{Some example sets appearing in the calculation of $f\vartriangleleft g$.} \label{fig:prop55}
\end{figure}

\begin{Example}
	Take $t = (\LL\RR\RR\LL\LL, 10111)$.
	Then $[t]_1$ is a $\stdl$-interval.
	Then, for example, $w := a^{\LL}\RR a^{\RR} b^{\LL} c^{\LL} \in \IWc(t)$
	and for $f,g \in \shalI$, one has
	\begin{equation*}
		(f \vartriangleleft g)(w)
		=
		f\big(a^{\LL} a^{\RR}\big)g\big(b^{\LL} c^{\LL}\big)
		+ f\big(a^{\LL}a^{\RR}c^{\LL}\big)g\big(b^{\LL}\big)
		+ f\big(a^{\LL}a^{\RR}b^{\LL}\big)g\big(c^{\LL}\big).
	\end{equation*}
\end{Example}
\begin{proof}[Proof of Proposition~\ref{prop:prelieproduct}]
	Let $f,g \in \shalI$ be linear functionals.

 We will show first that $f\vartriangleleft g(w)=0$ whenever $|\INT{w}|\neq 1$.
 \begin{itemize}\itemsep=0pt
 \item
 $|\INT{w}|=0$ means that $w\in\{\LL,\RR\}^{\star}$ is a word containing only placeholders. Then $\vcpSpr(w)$, $\vcpSsu(w)$ are zero and, therefore, $(f\vartriangleleft g)(w)=0$.
 \item
 If $|\INT{w}|>2$, then it is immediate to see that, in every term $w^{\prime}\overt w^{\prime\prime}$ occurring in the sum expressing $\vcpSpr(w)$, $\vcpSsu(w)$, one has $|\INT{w'}|\geq2$ or $|\INT{w''}|\geq2$ (or both), so that $f(w')=0$ or $g(w'')=0$. Consequently, $f\vartriangleleft g(w)=0$.
 \item
 Similarly, if $\INT{w}=\{I_1\stdl I_2\}$, there are exactly two ways to write $w=w'\circ w''$ with interval words $w'$, $w''$, namely
 \[w=\restriction{w}{I_1\cup [w]_0}\circ \translucidation{w}{I_1}=\restriction{w}{I_2\cup [w]_0}\circ \translucidation{w}{I_2},\]
 and, clearly,
 $\bigwedgedot\interval{|w|}\in I_1$.
 By \eqref{eq:m(S):noninterval}, all other summands in the sums expressing $\vcpSpr(w)$, $\vcpSsu(w)$ do not contribute when calculating $f\vartriangleleft g(w)$.
 By \eqref{eq:m(S):interval}, $f(\restriction{w}{I_1\cup [w]_0})=f(\restriction{w}{I_1})$ and $g(\translucidation{w}{I_1})=g(\restriction{w}{I_2})$,
 therefore $f\prec g(w)=f(\restriction{w}{I_1})g(\restriction{w}{I_2})$.
 Analogously, $g\succ f(w)=g(\restriction{w}{I_2})f(\restriction{w}{I_1})$. Therefore, $f\vartriangleleft g(w)=0$.
 \end{itemize}
 This proves that \eqref{eq:m(S):noninterval} holds for $f\vartriangleleft g$.

 It is easy to see that \eqref{eqn:formulaprelie} implies that \eqref{eq:m(S):interval} holds for $f\vartriangleleft g$. Therefore, showing \eqref{eqn:formulaprelie} will conclude the proof of the proposition.

 We start the second half of the proof with some general observations. If $w$ is an interval word and $\varnothing\neq I\subset [w]_1$, then
 \begin{itemize}\itemsep=0pt
 \item $\restriction{w}{I}$ is always an interval word;
 \item $\translucidation{w}{I}$ is an interval word if and only if one of the following situations applies: either $I$ is a $\stdl$-interval containing exactly one of the $\stdl$-endpoints of the $\stdl$-interval $\bigwedgedot [w]_1$, $\bigveedot [w]_1$, or $I=I_1\cup I_2$ is the union of two $\stdl$-intervals each containing one $\stdl$-endpoint, $\bigwedgedot [w]_1\in I_1$, $\bigveedot [w]_1\in I_2$.
 \end{itemize}

 Let $w\in \SWc$, $|\INT{w}|=1$. Using \eqref{eq:m(S):interval}, we can write
 \[(g \succ f)(w)
 =
 \sum_{\substack{[w]_0 \subset X\\ \bigwedgedots [w]_1 \notin X}}
 g(\restriction{w}{X})f(\translucidation{w}{X})
 =
 \sum_{\substack{J\subset[w]_1\\ \bigwedgedots [w]_1 \notin J}}
 f(\translucidation{w}{J})g(\restriction{w}{J}).
 \]
 The summands on the right-hand side vanish except when $\restriction{w}{J}$ and $\translucidation{w}{J}$ are both interval words. By the previous observations, $\restriction{w}{J}$ is automatically an interval word. On the other hand, $\translucidation{w}{J}$ is an interval word if and only if $J$ is a $\stdl$-interval containing $\bigveedot [w]_1$. In that case, $I:=[w]_1\setminus J$ is a $\stdl$-interval containing $\bigwedgedot [w]_1$ and $f(\translucidation{w}{J})=f(\restriction{w}{I})$ by \eqref{eq:m(S):interval}. This amounts to
 \[
 (g\succ f)(w)
 =
 \sum_{\substack{I\stdl J\text{ $\stdl$-intervals}\\I\cup J=[w]_1}}
 f(\restriction{w}{I})g(\restriction{w}{J}).
 \]
 For the second half product, we find
 \[(f \prec g)(w)
 =
 \sum_{\substack{[w]_0 \subset X\\ \bigwedgedots [w]_1 \in X}}
 f(\restriction{w}{X})g(\translucidation{w}{X})
 =
 \sum_{\substack{I\subset[w]_1\\ \bigwedgedots [w]_1 \in I}}
 f(\restriction{w}{I})g(\translucidation{w}{I}).
 \]
 If $I$ is a $\stdl$-interval, we get exactly the same terms as in the calculation of $(f\succ g)(w)$. However, for $J:=[w]_1\setminus I$ to be a $\stdl$-interval, it is enough that $I=I_1\cup I_2$ is the union of two $\stdl$-intervals $I_1\stdl I_2$ with $\bigwedgedot[w]_1\in I_1$, $\bigveedot [w]_1\in I_2$, giving additional summands corresponding to both~$I_1$ and~$I_2$ nonempty.

 By taking the difference between $(f\prec g)(w)$ and $(g\succ f)(w)$, we finally find~\eqref{eqn:formulaprelie}.
 \end{proof}

\begin{Proposition}\label{prop:monoid}
	The set $\shgrI$ is a monoid for the convolution product $\star$.
\end{Proposition}
\begin{proof}
The proof exploits the compatibility between the monoidal product $\overt$ and the semi-groupal product $\ominus$
exposed in Section~\ref{sec:bimonoid}. Following Definition~\ref{def:morphismU}, we note that $\shgrI$ is a set of morphisms of $\ominus$-semigroups,
	\begin{gather}
		f \in \shgrI
		\iff
		f(\hpIWc(w\ominus w^{\prime}))=f(w)f(w^{\prime}) \qquad \text{and}\nonumber\\
 f(u)=1 \quad \text{for all } w\ominus w^{\prime} \in \IWc\ominus \IWc, \quad u\in \{\LL,\RR\}^{\star}.\label{eq:f-in-M_I}
 	\end{gather}
	From this viewpoint, $f\in \shgrI(S)$ is considered valued in the linear $\SWc$-indexed collection $\CC^{\ominus}$ introduced in Remark~\ref{rk:monoidc}.
	{The $\IWc$-indexed collection} $\CC^{\ominus}$ supports a vertical product \[\vp_{\CC}\colon \ \CC^{\ominus} \overt {\CC}^{\ominus}\to \CC^{\ominus}\] (this product already appears, in disguise, in the definition of the convolution product $\star$):
	\begin{equation*}
		\vp_{\CC}(z \overt z^{\prime}) = zz^{\prime},\qquad z\overt z^{\prime} \in \mathbb{C}^{\ominus}\overt \mathbb{C}^{\ominus}
	\end{equation*}
	and another product $\hp_{\CC}\colon \CC^{\ominus} \ominus \CC^{\ominus} \to \CC^{\ominus}$,
	\begin{equation*}
		\hp_{\CC}(z \ominus z^{\prime}) = zz^{\prime},\qquad z \ominus z^{\prime} \in \mathbb{C}^{\ominus}\ominus \mathbb{C}^{\ominus}.
	\end{equation*}
 Both products come from the interpretation of the multiplication of complex numbers either as a vertical product or as a horizontal product.
	Of course, both products are compatible in the following sense:
 \[({\sf m}^{\overt}_{\CC}\circ (\hp_{\CC}\otimes \hp_{\CC})) \circ R_{\CC^{\ominus},\CC^{\ominus},\CC^{\ominus},\CC^{\ominus}} = {\sf m}^{\overt}_{\CC} \ominus {\sf m}^{\overt}_{\CC},\] which is equivalent to saying that $\vp_{\CC}$ is a morphism of $\ominus$-semigroups. Then \eqref{eq:f-in-M_I} is equivalent to
	\begin{equation*}
		f\in \shgrI \iff f\circ \hpIWc = \hp_{\CC} \circ (f \ominus f),\qquad f(u)=1, u \in \{\LL, \RR\}^{\star}.
	\end{equation*}
	Pick two linear functionals $f,g\in \shgrI$ and $w^{-}\ominus w^{+} \in \mathbb{C}\mathcal{W}\ominus \mathbb{C}\mathcal{W}$. Then, using Corollary~\ref{cor:vcpominusmorphism} and \eqref{eqn:exchange-m^ominus},
	\begin{align*}
		(f\star g)\big(\hpIWc(w^{-}\ominus w^{+})\big) & = \big(\vp_{\CC}\circ(f\overt g) \circ \vcp \circ \hpIWc\big) (w^{-}\ominus w^{+}) \\
		 & = \big(\vp_{\CC} \circ (f\overt g)\circ \hp_{\IWc\overt\IWc}\big)\big(\vcp(w^{-})\ominus\vcp(w^{+})\big) \\
		 & = \big(\vp_{\CC} \circ (f \overt g) \big(\hp_{\IWc}\overt \hp_{\IWc}\big)\circ R_{\IWc,\IWc,\IWc,\IWc}\big) \big(\vcp(w^{-})\ominus \vcp(w^{+})\big) \\
		 & = \big(\vp_{\CC}\circ \big(\hp_{\CC} \overt \hp_{\CC}\big)\circ ((f\ominus f) \overt(g \ominus g))\circ R_{\IWc,\IWc,\IWc,\IWc}\big)\\
&\quad {} \times\big(\vcp(w^{-}) \ominus \vcp(w^{+})\big).
	\end{align*}
	Owing to naturality of $R$, and using that
	\begin{gather*}
\vp_{\CC}\circ \big(\hp_{\CC}\overt \hp_{\CC}\big) \circ R_{\CC^{\ominus},\CC^{\ominus},\CC^{\ominus},\CC^{\ominus}}((z_1\overt z_2 ) \ominus(z_3 \overt z_4))\\
\qquad{} =z_1z_2z_3z_4=\hp_{\CC}\circ (\vp_{\CC}\ominus \vp_{\CC}) ((z_1\overt z_2) \ominus(z_3 \overt z_4)),
\end{gather*}
	we continue the previous calculation and get
	\begin{align*}
		(f\star g)\big(\hpIWc(w^{-}\ominus w^{+})\big) & =
\big(\vp_{\CC}\circ (\hp_{\CC}\overt \hp_{\CC}) \circ R_{\CC^{\ominus},\CC^{\ominus},\CC^{\ominus},\CC^{\ominus}} \circ ((f\overt g) \ominus(f \overt g))\big) \\
&\quad {}\times \big(\vcp(w^-) \ominus \vcp(w^+)\big) \\
		 & =\big(\hp_{\CC}\circ (\vp_{\CC}\ominus \vp_{\CC}) \circ ((f\overt g) \ominus(f \overt g))\big) \big(\vcp(w^-) \ominus \vcp(w^{+})\big) \\
		 & =\hp_{\CC} \big((f\star g)(w^-) \ominus (f\star g)(w^{+})\big),
	\end{align*}
	which shows that $f\star g \in \shgrI$. At this point, we have shown that $\shgrI$ is a semigroup. Every element of $\shgrI$ is invertible, what is left to prove is that this inverse is also in $\shgrI$. We postpone the proof of this point after Theorem~\ref{thm:shuffleexpo}.
\end{proof}

We will prove later on (Corollary~\ref{cor:shgrIisagroup} of Theorem~\ref{thm:shuffleexpo}) that $\shgrI$ is not only a monoid but a~group. The convolution product on $\shgrI$ is defined as the dual of the \emph{vertical} coproduct $\vcp_{\IWc}$, while elements in $\shgrI$ are compatible with respect to a horizontal product.
In that respect, $\shgrI$~does not fit in the classical theory of \emph{pro-algebraic groups.} To emphasize the differences, we briefly recall what a pro-algebraic group is.
Such a group $G_B(A)$ can be represented as the convolution group of algebra morphisms of a (plain, classical) graded connected augmented co-nilpotent bialgebra (in the category Vect) $\big(\bar{B} = \CC \cdot 1_B \oplus B, \Delta, m, \varepsilon\big)$ with values in a commutative unital algebra $A$,
			\[
				G_B(A) = \mathrm{Hom}_{\mathrm{Alg}}(B,A).
			\]
			\emph{An infinitesimal morphism} $\kappa\colon B \to A$ is a map such that $\kappa(1_B)=0$ and $\kappa(ab)=\varepsilon(a)\kappa(b) + \kappa(a)\varepsilon(b)$. The space $\mathfrak{g}_B(A)$ of infinitesimal morphisms is a Lie algebra for the bracket $[f,g]=f\star_{B}g-g\star_B f$. The two sets $G_{B}(A)$ and $\mathfrak{g}_B(A)$ are in bijection through the exponentials and its inverse the logarithm map:
			\[
				\exp_{\star_B}(\alpha) := \varepsilon + \sum_{n\geq 1} \frac{1}{n!}\alpha^{\star_B n },\qquad \ln_{\star_B}(\varepsilon + \alpha) := \sum_{n\geq 1} \frac{(-1)^{n-1}}{n}\alpha^{\star_b n}
			\]
			with $\alpha$ an infinitesimal morphism. We state first an equivalent result in our setting and continue the discussion then.
The following theorem follows from Theorem~\ref{thm:shuffleexpo} and conilpotency of~$\big(\mathcal{W},\vpS\big)$.
\begin{Theorem}
	\label{prop:liegroupun}
	The Lie algebra $\shalI$ and the group $\shgrI$ are in bijection through
	\[
		\exp_{\star}(m) : = \varepsilon + \sum_{n\geq 1} m^{\star_B n},\qquad \ln (\varepsilon+m) = \sum_{n\geq 1}\frac{(-1)^{n-1}}{n} m^{\star_B n}.
	\]
\end{Theorem}

\begin{Remark}
A \emph{double monoid} in the category $\SWc$ is a tuple $\big(O, \hp_O, \vp_O, \eta^{\overt}_O\big)$ with $O$ a~$\SWc$-collection and
	\begin{equation*}
		\hp_{O}\colon\ O \ominus O \to O,\qquad \vp_O\colon \ O \overt O \to O,\qquad \eta^{\overt}_O\colon \ \uvert \to O
	\end{equation*}
	such that $\big(O,\vp_O,\eta^{\overt}_O\big)$ is a monoid in $\big(\SWc,\overt,\uvert\big)$ and $\big(O,\hp_O\big)$ is a horizontal semi-group. Additionally, we require that $m^{\overt}_O$ is a morphism of $\ominus$-semigroups ($O \overt O$ is equipped with the product $m^{\ominus}_{O\overt O}= \big(m^{\ominus}_O \overt m^{\ominus}_O\big) \circ R_{O,O,O,O}$).
	Then, the proof of Proposition \ref{prop:monoid} generalizes to this context and the convolution product
	\begin{equation*}
		f \star g := \hp_{O}\circ (f \overt g) \circ \vcp, \qquad f,g \in \mathrm{Hom}_{\icsw}(\CC\IWc,O)
	\end{equation*}
restricts to a group product on the set $\mathrm{Hom}_{\mathrm{Semigroup}(\ominus)}(\IWc,O)$ of $\ominus$-semigroupal morphisms $f\colon \CC\IWc \to O$ with $f\circ\eta= \eta^{\overt}_O$. Besides, Proposition~\ref{prop:liegroupun} has a straightforward generalization in this setting.
\end{Remark}

In the remaining part of this section, we compute the two shuffle exponentials $\exp_{\prec, \succ}$ and the (shuffle) exponential $\exp_{\star}$ by using the defining fixed-point equations:
\[
	\exp_{\prec}(k) = \varepsilon + k\prec \exp_{\prec}(k),\qquad \exp_{\succ}(b) = \varepsilon + \exp_{\succ}(b) \succ b.
\]

We are now ready to state our main theorem.

\begin{Theorem}[left half-shuffle fixed point equation]
	\label{thm:lefthalfshuffle}
	Pick a linear form $k$ in the Lie algebra~$\shalI$. The solution $M$ of the following left half-shuffle fixed point equation
	\begin{equation*}
		M = \varepsilon + k \prec M
	\end{equation*}
	is in $\shgrI$ and satisfies the formula
	\begin{equation}\label{eqn:toprove}
		M(w) = \sum_{\pi \in \shBNCc(t)} \prod_{V\in\pi}k(\restriction{w}{V}),\qquad w\in \IWc(t).
	\end{equation}
\end{Theorem}
\begin{proof}
	We prove first using induction on the number of on-letters $|[t]_1|$ that
	\begin{equation}\label{eqn:multiplicative}
 M(w) = \prod_{J\in \INT{w}} M(\restriction{w}{J}),\qquad w \in \IWc(t),\qquad t \in \SWc.
	\end{equation}
	First, if $w \in \IWc(t)$ and $|[t]_1|=0$, equation \eqref{eqn:multiplicative} is trivial.

	Pick $N\geq 1$ and assume that \eqref{eqn:multiplicative} holds for any word $w \in \IWc(t)$, with $t\in \SWc$ a translucent word with at most $N$ ``opaque'' positions; $|[t]_1|\leq N$. Since $k \in \shalI$, for every subset $I \subset [t]_1$ with
 $\bigwedgedot [t]_1 \in I {\rm~and~} I \cap I^{\stdl}_j(t) \neq \varnothing$ for some $j\geq 2$,
 one has $k(\restriction{w}{\hat{I}})=0$ (recall the notation $\hat I = I \cup [t]_0$).
 We thus obtain (and recalling that for any incomplete word $v$ we denote by $v_{i}^\stdl$ its restriction to the $i$th interval in $\INT{v}$),
	\begin{align*}
 M(w) & = (k \prec M)(w) = \sum_{\substack{I \subset I^{\stdl}_1 \\ \bigwedgedots [t]_1 \in I}} k(\restriction{w}{\hat{I}})M(\translucidation{w}{\hat{I}}) \\
 & = \sum_{\substack{I \subset I^{\stdl}_1 \\ \bigwedgedots [t]_1 \in I}} k(\restriction{w}{\hat{I}})M\big((\translucidation{w}{\hat{I}})^{\stdl}_1\big)M\big((\translucidation{w}{\hat{I}})^{\stdl}_2\big) \cdots M\big((\translucidation{w}{\hat{I}})^{\stdl}_{j}\big),
	\end{align*}
	where we have used the inductive hypothesis applied to the word $\translucidation{w}{\hat{I}}$ to obtain the last equality since $|[\translucidation{t}{\hat{I}}]_1|\leq |[t]_1|-1$. Since in the last equality the sum ranges over $I \subset I_1^{\stdl}(w)$,
	\begin{equation*}
		(\translucidation{w}{\hat{I}})^{\stdl}_1 = \restriction{w_1^{\stdl}}{I_1^{\stdl} \setminus I},\qquad (\translucidation{w}{\hat{I}})^{\stdl}_2 = w^{\stdl}_2,\qquad \dots, \qquad (\translucidation{w}{\hat{I}})^{\stdl}_{j} = w^{\stdl}_{j},
	\end{equation*}
	from which we deduce
	\begin{align*}
 M(w)
 &= \sum_{\substack{I \subset I^{\stdl}_1 \\ \bigwedgedots [t]_1 \in I}} k\big(\restriction{w_1^{\stdl}}{\hat{I}}\big)M\big(\restriction{w_1^{\stdl}}{I_1^{\stdl} \setminus I}\big)M\big(w^{\stdl}_2\big) \cdots M\big(w^{\stdl}_{j}\big)\\
					 &= (k\prec M)\big(w^{\stdl}_1\big)M\big(w^{\stdl}_2\big)\cdots M\big(w^{\stdl}_{j}\big) \\
					&=M\big(w^{\stdl}_1\big)M\big(w^{\stdl}_2\big)\cdots M\big(w^{\stdl}_{j}\big).
	\end{align*}

We prove next the formula \eqref{eqn:toprove}. Define $ G\colon \CC \IWc\to \CC$ by
	\begin{equation*}
		G(w) := \sum_{\pi\in{\shBNCc}(t)}\prod_{V\in \pi} k(\restriction{w}{V}),\qquad w\in \IWc(t),\qquad t \in \SWc.
	\end{equation*}

We prove now that $G$ is the solution to the fixed point equation, $G = \varepsilon + k \prec G$.

First, owing to Remark~\ref{rk:fondrk}, it is straightforward to check that $G$ is a multiplicative function in $\shgrI$. It will thus be sufficient to prove that for any word $w\in \IWc(w,\mathbf{1}_n)$, $w\in \{\LL,\RR\}^{\star}$, one has
	\begin{equation*}
		G(w) = (k\prec G)(w).
	\end{equation*}
	Pick such a word $w \in \IWc(t)$. Then it holds that
	\begin{equation*}
 (k\prec G)(w) = \sum_{\substack{I \subset [t]_1 \\ \bigwedgedots [t]_1 \in I}} k(\restriction{w}{\hat{I}})G(\translucidation{w}{\hat{I}}) = \sum_{\substack{I \subset [t]_1 \\ \bigwedgedots [t]_1 \in I}} k(\restriction{w}{\hat{I}}) \sum_{\pi \in \shBNCc(\translucidation{t}{\hat{I}})} \prod_{W\in \pi} k(\restriction{w}{W}).
	\end{equation*}

Pick a noncrossing bipartition $\pi$ in $\BNCc(t)$ and call $\pi^{\wedgedot}$ the block that contains $\bigwedgedot[t]_1$. First, the translucidation $\translucidation{\pi}{\pi^{\wedgedots}}$ of $\pi$ yields a shaded noncrossing bipartition in
	${\shBNCc}(\translucidation{t}{\pi^{\wedgedots}})$.

Secondly,
	\[
		\BNCc(t)\ni \pi \mapsto \big(\pi^{\wedgedot}, \translucidation{\pi}{\pi^{\wedgedots}}\big)
	\]
	is a bijection onto the set of pairs $(I, \pi)$ with $\bigwedgedot[t]_1\in I\subset \llbracket 1,|t|\rrbracket$ and $\pi\in {\shBNCc}(\translucidation{t}{I})$. One concludes the proof with the following equality
	\begin{equation*}
 \sum_{\pi \in \BNCc(t)}\prod_{W\in \pi} k(\restriction{w}{W})= \sum_{\substack{I \subset [t]_1 \\ \bigwedgedots [t]_1 \in I}} k(\restriction{w}{\hat{I}}) \sum_{\pi \in {\shBNCc}(\translucidation{t}{\hat{I}})} \prod_{W\in \pi} k(\restriction{w}{W}).\tag*{\qed}
	\end{equation*}\renewcommand{\qed}{}
\end{proof}

\begin{Example}
	Pick $w=a^{\RR}a^{\LL}$, then
	\begin{align*}
		M(w) = k\big(a^{\LL}\big)M\big(a^{\RR}\big) + k\big(a^{\RR} a^{\LL}\big) = k\big(a^{\LL}\big)k\big(a^{\RR}\big) + k\big(a^{\RR} a^{\LL}\big).
	\end{align*}
	Pick $w = a^{\RR} a^{\LL} a^{\RR} a^{\LL}$,
	\begin{gather*}
		M(w) = k\big(a^{\LL}\big) M\big(a^{\RR}a^{\RR} a^{\LL}\big) + k\big(a^{\RR} a^{\LL} \big)M\big(a^{\RR} a^{\LL}\big) + k\big(a^{\LL} a^{\RR}\big)M\big(a^{\RR}\big)M\big(a^{\LL}\big) \\
\hphantom{M(w) =}{}
+ k\big(a^{\LL} a^{\LL}\big)M\big(a^{\RR}a^{\RR}\big)+ k\big(a^{\RR}a^{\LL}a^{\RR}\big)M\big(a^{\LL}\big) + k\big(a^{\RR} a^{\LL} a^{\LL}\big)M\big(a^{\RR}\big) \\
\hphantom{M(w) =}{}
+ k\big(a^{\LL} a^{\RR} a^{\LL}\big)M\big(a^{\RR}\big)+ k\big(a^{\RR} a^{\LL} a^{\LL} a^{\LL}\big) \\
\hphantom{M(w)}{}
 = k\big(a^{\LL}\big)\big(k\big(a^{\LL}\big)M\big(a^{\RR}a^{\RR}\big) + 2k\big(a^\RR a^{\LL}\big)k\big(a^{\RR}\big)\big) + k\big(a^{\RR}a^\RR a^{\LL}\big)\big) \\
\hphantom{M(w) =}{}
+ k\big(a^{\RR} a^{\LL}\big)\big(k\big(a^{\LL}\big)k\big(a^{\RR}\big) + k\big(a^{\RR} a^{\LL}\big)\big) \\
\hphantom{M(w) =}{}
+ k\big(a^{\LL} a^{\RR}\big)k\big(a^{\RR}\big)k\big(a^{\LL}\big) + k\big(a^{\LL}a^{\LL}\big)\big(k\big(a^{\RR}\big)k\big(a^{\RR}\big)+k\big(a^\RR a^\RR\big)\big) \\
\hphantom{M(w) =}{}
		 + k\big(a^{\RR} a^{\LL} a^{\RR}\big)k\big(a^{\LL}\big) + k\big(a^{\RR} a^{\LL} a^{\LL}\big)k\big(a^{\RR}\big) + k\big(a^{\LL} a^{\RR} a^{\LL}\big)k\big(a^{\RR}\big) +k\big(a^{\RR} a^{\LL} a^{\RR} a^{\LL}\big) \\
\hphantom{M(w)}{}
 = k\big(a^{\LL}\big)\big(k\big(a^{\LL}\big)\big(k\big(a^{\RR}a^{\RR}\big)+k\big(a^\RR\big)k\big(a^\RR\big)\big) + 2k\big(a^\RR a^{\LL}\big)k\big(a^{\RR}\big)\big) + k\big(a^{\RR}a^\RR a^{\LL}\big)\big)\\
\hphantom{M(w) =}{}
 + k\big(a^{\RR} a^{\LL}\big)\big(k\big(a^{\LL}\big)k\big(a^{\RR}\big) + k\big(a^{\RR} a^{\LL}\big)\big)
		 + k\big(a^{\LL} a^{\RR}\big)k\big(a^{\RR}\big)k\big(a^{\LL}\big) \\
\hphantom{M(w) =}{}
+ k\big(a^{\LL}a^{\LL}\big)\big(k\big(a^{\RR}\big)k\big(a^{\RR}\big)+k\big(a^\RR a^\RR\big)\big) + k\big(a^{\RR} a^{\LL} a^{\RR}\big)k\big(a^{\LL}\big) \\
\hphantom{M(w) =}{}
		 + k\big(a^{\RR} a^{\LL} a^{\LL}\big)k\big(a^{\RR}\big) + k\big(a^{\LL} a^{\RR} a^{\LL}\big)k\big(a^{\RR}\big)
		 +k\big(a^{\RR} a^{\LL} a^{\RR} a^{\LL}\big).
	\end{gather*}
\end{Example}

\begin{figure}[!ht]\centering
	\includegraphics[scale=0.7]{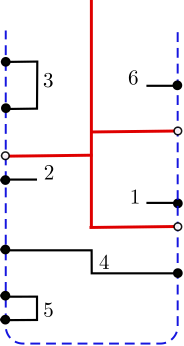}
	\caption{Example of an incomplete monotone bipartitions.}\label{fig:incompletebinoncrossingpartitionmonotone}
\end{figure}

\begin{Theorem}[shuffle exponential]
	\label{thm:shuffleexpo}
	Pick a linear function $m$ in the Lie algebra $ \shalI$. Then, for any word $w\in \IWc(t)$, one has
	\begin{equation*}
		\exp_{\star}(m)(w) = \sum_{(\pi,\lambda) \in \MBNCc(t)} \frac{1}{|\pi|!}\prod_{W\in\pi}m(\restriction{w}{W}).
	\end{equation*}
\end{Theorem}
\begin{proof}Pick $m\in \shalI$.
	We let $\MBNCc^{\ell}(t)$ denotes the set of monotone bipartitions of type $t$ with~$\ell$ blocks (the block~$[t]_{0}$ is not counted). We prove by induction on $\ell \geq 1$ that for any word~$w$ in~$\IWc(t)$,
	\begin{equation} \label{eqn:toshow}
		m^{\overt \ell}\circ \big(\vcpS\big)^{\ell} = \sum_{(\pi,\lambda) \in \MBNCc^{\ell}(t)} \prod_{W\in\pi}m(\restriction{w}{W}),
	\end{equation}
	where $\big(\vcp\big)^{\ell} = \big(\big(\vcp\big)^{\ell-1}\overt \mathrm{id}\big)\circ \vcp$, $\big(\vcp\big)^{0}=\mathrm{id}$.

	The formula \eqref{eqn:toshow} holds if $\ell=1$. Pick an integer $N\geq 1$ and suppose that equation \eqref{eqn:toshow} holds for every integer $\ell \leq N$.

	Using $\big({\vcp}\big)^{N+1} = \big(\big({\vcp}\big)^{N}\overt \mathrm{id}\big)\circ {\vcp}$ and the inductive hypothesis, one gets
	\begin{align*}
		\big(m^{\overt (N+1)}\circ \big(\vcp\big)^{N+1}\big)(w) & = \sum_{I \subset [t]_1} \big(m^{\overt N} \circ \big(\vcp\big)^{N}\big)(\restriction{w}{\hat{I}})m(\translucidation{w}{\hat{I}}) \\
		&= \sum_{\substack{I \subset [t]_1 \\ \translucidation{t}{\hat{I}} \in {\SWint}}} \big(m^{\overt N} \circ \big(\vcp\big)^{N}\big)(\restriction{w}{\hat{I}})m(\translucidation{w}{\hat{I}}) \\
		 & =\sum_{\substack{I \subset [t]_1 \\ \translucidation{t}{\hat{I}} \in \SWint}} \bigg[\sum_{{(\pi^{\prime},\lambda') \in \MBNCc^N(\restriction{t}{\hat{I}})}}\prod_{W\in\pi^{\prime}}m(\restriction{w}{W})\bigg] m(\translucidation{w}{\hat{I}}) \\
		 & = \sum_{\substack{I\subset [t]_1 \\ [t]_1 \setminus I\ \stdl \text{-interval}}} \bigg[\sum_{{(\pi^{\prime},\lambda') \in \MBNCc^N(\restriction{t}{\hat{I} })}}\prod_{W\in\pi^{\prime}}m(\restriction{w}{W})\bigg] m(\translucidation{w}{\hat{I}}).
	\end{align*}
	Pick a monotone partition $(\pi^{\prime}, \lambda')\in \MBNCc^N(\restriction{t}{\hat{I}})$ with $I \subset [t]_1$ and $J:=[t]_1\setminus I$ a $\stdl$-interval.
	There exists a unique incomplete bipartition $\hat{\pi} \in \BPc(t)$ such that $\restriction{\hat{\pi}}{\hat{I}} = \pi^{\prime}$ and $\translucidation{\hat{\pi}}{\hat{I}} = \{[t]_1 \setminus I\}=\{J\}$ an incomplete one block partition of $\translucidation{t}{I}$, namely $\hat\pi=\pi'\circ \{J\}$.
	We claim that $\hat{\pi}$ is an incomplete noncrossing bipartition: this comes from the fact that $J$ a $\stdl$-interval and $\pi^{\prime}$ is noncrossing.

	Besides, $\hat{\pi}$ is a shaded noncrossing bipartition.
	The translucent block of $\hat{\pi}$ is $[t]_0$.
	Pick an integer $1\leq k< \min [t]_0$ (consequently, $k\in [t]_1$) with $\alpha_t(k)=\alpha_t(\min [t]_0)$.

	If $k \in {I}$, then $j \sim_{\hat\pi} k$ implies $j \in I$ because $I$ is the complement of the block $J$ inside $[t]_1$. Thus, $j < \min [t]_0 \text{ and } \alpha(j)=\alpha(\min [t]_0)$
	since $\restriction{\hat{\pi}}{\hat{I}}$ is a shaded noncrossing bipartition.

	If $k$ is not in $I$, since $J$ is a $\stdl$-interval, one must have $\alpha_{t}(k) = \alpha_{t}(j)$ for all $j\in J$.
	For example, suppose $\alpha_{t}(\min [t]_0) = {\sf L}$ one has $\alpha_{t}(k)={\sf L}$ and $k\stdl \min [t]_0$. Again, since $J$ is a $\stdl$-interval, this last inequality implies that $j \stdl \min [t]_0$ and, thus, $\alpha_t(j)={\sf L}$ for any $j \in [t]_1 \setminus I$.

 The monotone order $\lambda^{\prime}$ on $\pi^{\prime}$ yields an order $\hat{\lambda}$ on $\hat{\pi}$ labelling the block $J=[t]_{1}\setminus I$ by $|\hat\pi|=|\pi'|+1$ highest.
	Because $J$ is a $\stdl$-interval, $\hat{\pi}$ endowed with the order induced by $\hat\lambda$ described above, is a monotone bipartition.

	Finally, to any monotone partition $(\pi, \lambda)$ in $\MBNCc(t)$ corresponds a unique pair $((\pi^{\prime},\lambda^{\prime}),I)$ with $I \subset [t]_1$, $[t]_1 \setminus I \text{ a }\stdl \text{-interval}$ and $(\pi^{\prime},\lambda^{\prime})\in\MBNCc(\restriction{t}{I})$
	such that $\pi=\hat{\pi}$ and $\lambda = \hat{\lambda}$; indeed, $J=[t]_1\setminus I$ has to be the highest block of $\pi$, which is always a $\stdl$-interval. This concludes the proof.
\end{proof}
\begin{Corollary}
	\label{cor:shgrIisagroup}
	The monoid $\shgrI$ is a group with Lie algebra $\shalI$, that is
	\[
	\exp_{\star}(\shalI) = \shgrI.
	\]
\end{Corollary}
\begin{proof}
	First, Theorem~\ref{thm:shuffleexpo} implies $\exp_{\star}(\shalI)\subset \shgrI$.
	Pick $M \in \shgrI$ and define $m \in\shalI$ inductively on the length of the word $\restriction{w}{[t]_0}$, for $w\in \IWc(t)$,
 with $|\INT{w}|=1$
	\begin{align*}
		 & m(w)
		 :=m(\restriction{w}{[t]_1})
		 :=M(\restriction{w}{[t]_1})
		 -\sum_{\pi\in \MBNCc(\restriction{t}{[t]_1})}\frac{1}{|\pi|!}
		 \prod_{V\in \pi}m(\restriction{w}{V})
	\end{align*}
	with initial conditions $m(s) = M(s)$, $s \in S$. Again, Theorem~\ref{thm:shuffleexpo} yields $\exp_{\star}(m)=M$. Hence $\exp_{\star}(\shalI) = \shgrI$.
\end{proof}

We end this section with the computation of the right half-shuffle exponential of an element of $\shalI$.

\begin{Theorem}
	Pick a function $b$ in the Lie algebra $\shalI$. The solution $M$ of the following right half-shuffle fixed point equation
	\begin{equation}
		\label{eqn:fixedpointright}
		M = \varepsilon + M \succ b
	\end{equation}
	is in $\shgrI$ and satisfies the formula
	\begin{equation*}
		M(w) = \sum_{\pi \in \BBc(t)} \prod_{V\in\pi}b(\restriction{w}{V}),\qquad w\in \IWc(t),\qquad t\in \SWc.
	\end{equation*}
\end{Theorem}
\begin{proof}
	To prove $M\in \shgrI$, the reasoning we used in proving the same property for the solution of the left half-shuffle fixed point equation in Theorem~\ref{thm:lefthalfshuffle} applies here. We omit the details for brevity and move on to the proof of the second statement, the formula for $M$ namely.

	We use the notation introduced in the proof of Theorem~\ref{thm:lefthalfshuffle}. We set
	\begin{equation*}
		G(w) := \sum_{\pi\in\BBc(t)}\prod_{V\in \pi}G(\restriction{w}{V}),\qquad w \in \IWc(t),\qquad t\in\SWc,
	\end{equation*}
	and show that $G$ satisfies to the fixed point equation \eqref{eqn:fixedpointright}. Pick an opaque word $t \in \SWc$, i.e., $b_t = \mathbf{1}$, and a word $w\in\IWc(t)$. Then one gets
	\begin{align}
		\label{eqn:recright}
		(M\succ b)(w) & = \sum_{\substack{I \subset \llbracket |t| \rrbracket,\\ \bigwedgedots\interval{|t|} \in \interval{|t|}\setminus I}} M(\restriction{w}{I})b(\translucidation{w}{I}) = \sum_{\substack{I \subset \interval{|t|},\\ \bigwedgedots\interval{|t|} \in \interval{|t|}\setminus I}} b(\translucidation{w}{I}) \sum_{\pi \in \BBc(\restriction{t}{I})} \prod_{V\in\pi}b(\restriction{w}{V}) \nonumber \\
		 & =\sum_{\substack{\bigwedgedots\interval{|t|}\in I\subset \interval{|t|} \\ I\text{ $\stdl$-interval}}} b(\restriction{w}{I})\sum_{\pi \in \BBc(\restriction{t}{\interval{|t|}\setminus I})} \prod_{V\in\pi}b(\restriction{w}{V}).
	\end{align}

	To any interval bipartition $\pi^{\prime}$ in $\BBc(t)$, one associates a unique pair $(I,\pi)$ of an interval $I\subset \interval{|t|}$ for the order $\stdl$ and an interval bipartition $\pi \in \BBc(\translucidation{t}{I})$: $I$ is the block of $\pi^{\prime}$ which contains $\bigwedgedot \llbracket |t| \rrbracket$
	and $\pi$ is the restriction of $\pi^{\prime}$ to $\llbracket |t| \rrbracket\setminus I$. The function
	$\BBc(t) \ni \pi' \mapsto (I, \pi) $
	is a bijection between the set $\BBc(t)$ and the set of pairs $(I,\pi)$ with $I$ an interval of $\llbracket |t| \rrbracket$ which contains $\bigwedgedot \llbracket |t| \rrbracket$ and $\pi \in \BBc(\translucidation{t}{I})$.

	The double sum on the right-hand side of equation \eqref{eqn:recright} can thus be rewritten as a sum over interval bipartitions in $\BBc(t)$ as follows:
	\begin{equation*}
		(M\succ b)(w) = \sum_{\substack{\bigwedgedots [t]_1\in I\subset [t]_1 \\
		\text{$I$ $\stdl$-interval}}} b(\restriction{w}{I})\sum_{\pi \in \BBc(\restriction{\alpha}{\llbracket |t| \rrbracket\setminus I})} \prod_{V\in\pi}b(\restriction{w}{V}) = \sum_{\pi^{\prime}\in \BBc(t)} \prod_{V\in\pi^{\prime}}b(\restriction{w}{V}).
	\end{equation*}
	This concludes the proof.
\end{proof}

\subsection{preLie algebra II}
The purpose of this section is to exhibit a second sub-preLie algebra $\shalII$ of $\mathfrak{g}$. The preLie product~$\vartriangleleft$ restricts to the zero preLie product on $\shalII$.
Both preLie algebras, $\shalI$ defined in the previous section and $\shalII$ introduced now, are meaningful to the theory of noncommutative probability; they \emph{implement} moment-cumulant relations corresponding to certain notions of independence. For the one at stake here, this is \emph{tensor independence} between random variables in $S$.\footnote{The combinatorial structures we study here only lead to trivial bi-tensor independence, which coincides with usual tensor independence because the full set of partitions of a translucent set is not constrained by colouration. Non-trivial bi-tensor independences, such as appear in \cite{GHU21p} and \cite{Varso-PhD} as deformations are beyond our scope here.}

\begin{Definition}\label{def:M}We introduce the following subsets of linear functionals on incomplete words (the index $\mathsf P$ stands for partitions):
	\begin{itemize}\itemsep=0pt
		\item We denote by $\shgrII$ the set of all linear forms $f\colon \linearspan{\IWc}\rightarrow \CC$ with
		 \begin{itemize}\itemsep=0pt
			 \item $f(u) = 1$ for all $u \in \IWc(\alpha,\mathbf{0})$, $\alpha \in \{\LL,\RR\}^{\star}$,
\item $f(w)=f(\restriction{w}{[t]_1})$ for all $w\in\IWc(t)$, $t\in\SWc$, $b_t\neq \mathbf{0}$.
		 \end{itemize}
		\item We denote by $\shalII$ the set of all linear forms $f\colon \linearspan{\IWc}\to \CC$ such that
		 \begin{itemize}\itemsep=0pt
			 \item $f(u) = 0$ for all $u \in \IWc(\alpha,\mathbf{0})$,
			 \item $f(w)=f(\restriction{w}{[t]_1})$ for all $w\in\IWc(t)$, $t\in\SWc$, $b_t\neq \mathbf{0}$.
		 \end{itemize}
	\end{itemize}
\end{Definition}

\begin{Remark}
As vector spaces, $\shalI$ and $\shalII$ sit in a class of subspaces of $\mathfrak{g}$ parametrized by subsets of \emph{incomplete bipartitions with one block $($besides the translucent block$)$}.
To an incomplete word on random variables $w\in \IWc$, we associate the incomplete bipartition $\pi_{w} = \{[w]_1\}$ with translucent block $[w]_0$.

Pick a subset $U$ of one-block incomplete bipartitions. At the moment, we do not assume anything on $U$. We set
\begin{equation*}
	\mathfrak{m}_U := \{\ell \in \mathfrak{g}\colon \ell(w)=0 \text{ if }\pi_w \not\in U \text{ and }\ell(w)=\ell(\restriction{w}{[w]_1}) \text{ if } \pi_{w} \in U\}
\end{equation*}
for the space $\shalI$, the corresponding set $U$ is the set $\mathsf I$ of incomplete interval bipartitions with only an opaque block (or, what is the same, shaded noncrossing bipartitions with only one opaque block) while for $\shalII$ the corresponding set is $\mathsf P$, the set of all incomplete bipartition with only one opaque block.
\begin{figure}[!ht]\centering
	\includegraphics[scale=0.7]{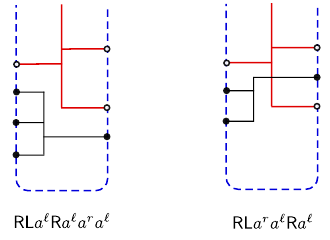}
	\caption{Shaded noncrossing partitions $\pi_w$ and words $w$.}\label{fig:piw}
\end{figure}
A natural question would then address the complete classification of preLie algebras arising as $\mathfrak m_U$ and \emph{meaningful} to the theory of probability. We leave such general considerations for future work.
\end{Remark}

\begin{Proposition}	\label{prop:trivial} The restriction of the preLie product $\vartriangleleft$ to $\shalII$ is trivial,
	\begin{equation*}
		f \vartriangleleft g = 0,\qquad f,g \in \shalII.
	\end{equation*}
	In particular, $\shalII$ is a sub preLie algebra.
\end{Proposition}
\begin{proof}
	Pick two linear forms $f,g \in \shalII$ and a word $w\in \IWc(t)$. We have
	\begin{align*}
		(f\vartriangleleft g)(w) = (f\prec g - g\succ f)(w) = \sum_{\substack{I \subset [t]_1 \\ \bigwedgedots [t]_1 \in I}}f(\restriction{w}{I})g(\translucidation{w}{I}) - \sum_{\substack{I\subset [t]_1 \\ \bigwedgedots [t]_1 \not\in I}}g(\restriction{w}{I})f(\translucidation{w}{I}).
	\end{align*}
	The result is obtained upon changing the summation variable of the second sum to $[t]_1\setminus I$ in the above equality.
\end{proof}
\begin{Theorem}
	Pick a linear form $k \in \shalII$. For any word $w\in \IWc(t)$, $t\in \SWc$, the following formula holds:
	\begin{equation}
		\label{eqn:formulaclassical}
		\exp_{\star}(k)(w)=\exp_{\succ}(k)(w) = \exp_{\prec}(k)(w) = \sum_{\pi \in \BPc(t)} \prod_{V\in\pi} k(\restriction{w}{V}).
	\end{equation}
\end{Theorem}
\begin{proof}
	Pick a linear form $k\in\shalII$. The equality between the left and right half-shuffle exponentials and the full shuffle exponential is a corollary of Proposition~\ref{prop:trivial}. This is a direct application of Theorem 18 in \cite{EbrahimiFardPatras2019applications}: $\Omega^{\prime}={\rm id}$.

 We prove next the last equality in \eqref{eqn:formulaclassical} by using the left half-shuffle fixed point equation. To that end, denote by $G$ the linear functional on $\CC \IWc$ whose values are prescribed by the rightmost terms in \eqref{eqn:formulaclassical}. The proof proceeds from the same lines of arguments we used to prove Theorem~\ref{thm:lefthalfshuffle} and Theorem~\ref{thm:shuffleexpo}.
	The result trivially holds for incomplete words $w$ with $|[w]_1| = 1$. Pick an integer $N\geq 1$ and assume that the result holds for incomplete words $w$ with $|[w]_1| \leq N$. Pick a translucent word $t\in\SWc$ and a word $w \in \IWc(t)$ with $|[w]_1| = N+1$ then
	\begin{equation*}
		(k \prec G)(w) = \sum_{\substack{I\subset [t]_1 \\ \bigwedgedots [t]_1 \in I}} k(\restriction{w}{I}) G(\translucidation{w}{ \hat{I}}) = \sum_{\substack{I \subset [t]_1 \\ \bigwedgedots [t]_1 \in I }} k(\restriction{w}{I})\sum_{\pi\in\BPc(\translucidation{t}{\hat{I}})} \prod_{V\in \pi} G(\restriction{w}{V}),
	\end{equation*}
	where we have used the inductive hypothesis on the word $\translucidation{w}{\hat{I}} $ to derive the last equality. To a partition $\pi \in \BPc(t)$, we associate the pair $(\pi^{\wedgedot}, \translucidation{\pi}{\pi^{\wedgedots}})$ where we recall that $\pi^{\wedgedot}$ is the block of~$\pi$ containing $\wedgedotb[t]_1$. As before, the application
	\[
	\mathcal{P}(t)\ni \pi \mapsto \big(\pi^{\wedgedot}, \translucidation{\pi}{\pi^{\wedgedots}}\big)
	\]
	is injective. Going backward, given a pair $(V,\pi^{\prime})$ of a subset $V\subset [t]_1$ and a partition $\pi^{\prime} \in \BPc(\translucidation{t}{\hat{V}})$, by defining $\pi := 1_{V} \circ \pi^{\prime} $, where \[1_{V} = \{\{1,\dots,|V|\}\} \in \BPc(\restriction{t}{V}),\] one has $\pi^{\wedgedot} = V$ and $\translucidation{\pi}{\pi^{\wedgedots}}=\pi^{\prime}$ (recall that the composition $\circ$ has been defined in Section~\ref{subsec:incompletebinoncrossing}). Hence,
	\begin{align*}
		(k\prec G)(w) = \sum_{\substack{I \subset [t]_1 \\ \bigwedgedots [t]_1 \in I }} k(\restriction{w}{I})\sum_{\pi\in\BPc(\translucidation{t}{ \hat{I}})} \prod_{V\in \pi} G(\restriction{w}{V}) = \sum_{\pi \in \BPc(t)} \prod_{V\in\pi} k(\restriction{w}{V}) = G(w).\tag*{\qed}
	\end{align*}\renewcommand{\qed}{}
\end{proof}

\section{Conclusion and perspectives}

\subsection{Summary}

In this work, we extended the shuffle-algebraic perspective on the free, boolean, and monotone moment-cumulant relations to bifree, biBoolean, and bimonotone moment-cumulant relations. To achieve this, we leverage an idea latent in \cite{ebrahimi2020operads} abridged as follows.
There exists a certain unshuffle structure on noncrossing partitions which, when restricted to ``sticks'' (partitions into singletons) yields the unshuffle coalgebra introduced in \cite{ebrahimi2015cumulants}. That former unshuffle structure stems from a certain composition rule, formalized using the theory of operads, between noncrossing partitions, \emph{the gap insertion operad} in which a partition is inserted into another by ``placing'' blocks in between two consecutive elements of the latter.
In the current work, these gap-insertion operations are extended to noncrossing bipartitions and we chose inputs gaps between two consecutive elements for a necklace order coming with any bipartition. Formalizing, we were led to introduce:
\begin{enumerate}\itemsep=0pt
\item a set of shaded noncrossing bipartitions, following the terminology introduced in \cite{charlesworth2015combinatorics}, that are bipartitions with one-pointed block acting as a placeholder,
\item a M\"obius category $\SWc$ whose morphisms (we call them translucent words) index sets of objects (e.g. shaded noncrossing bipartitions and incomplete words) we compose.
\item To formalize compositions and following \cite{content1980categories}, we defined a certain tensor product $\overt$ on $\icsw$.
\item Additionally, we showed that this category supports an additional semigroupal product $\ominus$ compatible with $\overt$, see Proposition \ref{prop:exchange}.
\item We defined a double monoid $\mathcal{W}$ in $\icsw$ supported by incomplete words also equipped with a certain unshuffle structure (after internalization to the incidence category of $\SWc$ of this notion) compatible with the horizontal composition.
\item We proved finally that inside the convolution monoid of linear functionals on $\mathcal{W}$ sit two groups and their corresponding Lie algebras (defined in Section~\ref{sec:momentscumulantsrelations}) among which a subset of the group of horizontal characters of $\big(\mathcal{W},\hpIWc\big)$. This latter is equipped with three exponentials encoding the bifree, biBoolean, and bimonotone moment-cumulant relations.
\end{enumerate}

\subsection{Outlook}

As alluded to in the introductory section, other independences were introduced in the past few years, including operator-valued counterparts of the
ones investigated in this work. 
Going operator-valued should cause no additional problems,
since the algebra developed here is closely related to the one expounded in \cite{gilliers2020operator}.
An interesting challenge would be to ``blend'' to the shuffle picture developed above also other two-faced independences,
such as the boolean-free one \cite{liu2019}.
This can be approached in at least two different ways:
(1)~including boolean-free independence to a triple of
independences (loosing perhaps positivity), identifying the appropriate sub-preLie algebra of $\mathfrak{g}$ to cast boolean-free moment-cumulant relations as shuffle exponential type relations,
(2)~possibly supplement the shuffle operations $\prec$ and $\succ$ with others, compatible in a certain sense, to reach more moment-cumulant relations.

In this direction, as already pointed out in Section~\ref{sec:bimonoid}, incomplete
words on random variables come with two orders, the natural one and the
necklace one $\stdl$, one can track the first letter for each order to define
half-coproducts.
This yields four of them, among which only two are addressed in this work.

\subsection*{Acknowledgements}
The authors would like to thank Kurusch Ebrahimi-Fard, Michael Sch\"urmann, and Philipp Var\v{s}o for many discussions on moment-cumulant relations and related algebraic questions. We thank the anonymous referees for their valuable feedback. M.G.~was supported by the German Research Foundation (DFG) grant no.~397960675; part of his work was carried out during the tenure of an ERCIM ‘Alain
Bensoussan’ Fellowship Programme.
J.D.~and N.G.~are supported by the trilateral (DFG/ANR/JST) grant ``EnhanceD Data stream Analysis with the Signature Method''. N.G.~is supported by ANR ``STARS''.

\pdfbookmark[1]{References}{ref}
\LastPageEnding

\end{document}